\documentclass[12pt]{article}

\usepackage{tikz}
\usetikzlibrary{calc,through,backgrounds}
\usepackage{multicol}
\usepackage{bm,amssymb,mathtools,amsthm,amsfonts,amsmath,latexsym,cite,color, psfrag,graphicx,ifpdf,ytableau,pgfkeys,pgfopts,xcolor,extarrows}

\newtheorem{theo}{Theorem}[section]

\newtheorem{lemma}[theo]{Lemma}
\newtheorem{coro}[theo]{Corollary}
\newtheorem{prop}[theo]{Proposition}
\newtheorem{prob}[theo]{Problem}
\newtheorem{re}[theo]{Remark}

\allowdisplaybreaks

\makeatletter \@addtoreset{equation}{section}
\@addtoreset{theo}{}\makeatother

\usepackage{hyperref}

\def\qed{\hfill \rule{4pt}{7pt}}
\def\pf{\noindent {\it Proof.} }

\def\S{{\mathfrak{S}}}

\begin{document}
\begin{center}

{\Large\bf Bumpless
Pipedreams, Reduced Word Tableaux  and

 Stanley Symmetric
Functions }\\[12pt]
\begin{center}

\vskip 4mm
{\small  Neil J.Y. Fan$^1$,  Peter L. Guo$^2$, Sophie C.C. Sun$^3$}

\vskip 4mm

$^1$Department of Mathematics\\
Sichuan University, Chengdu,
Sichuan 610064, P.R. China \\[3mm]

$^{2,3}$Center for Combinatorics, LPMC-TJKLC\\
Nankai University,
Tianjin 300071,
P.R. China\\[3mm]

\vskip 4mm
$^1$fan@scu.edu.cn,
$^2$lguo@nankai.edu.cn,
$^3$suncongcong@mail.nankai.edu.cn

\end{center}

\end{center}

\vskip 6mm
\begin{abstract}
Lam, Lee and Shimozono   introduced the structure of
bumpless pipedreams  in their study of back stable Schubert calculus.
They found that a specific family of bumpless pipedreams,
called EG-pipedreams, can be used to interpret the Edelman-Greene coefficients
appearing in the expansion of a Stanley symmetric function in the basis of Schur functions. It is well known that the Edelman-Greene coefficients can  also
be interpreted in terms of reduced word tableaux for  permutations.
Lam, Lee and Shimozono proposed the problem of finding a
shape preserving bijection between   reduced word tableaux
for a permutation $w$ and  EG-pipedreams of $w$.
In this paper, we construct such a bijection.
The key ingredients are two new developed isomorphic
tree structures associated to $w$: the
modified Lascoux-Sch\"{u}tzenberger tree of $w$ and
the Edelman-Greene tree of $w$.
Using the Little map,
we show that the leaves in the modified Lascoux-Sch\"{u}tzenberger of $w$
are in bijection with the reduced word  tableaux for $w$. On the
other hand, applying the droop operation on bumpless pipedreams
also introduced by Lam, Lee and Shimozono, we show that the leaves in the Edelman-Greene
tree of $w$
are in bijection with the EG-pipedreams of $w$.
This allows us to establish a shape preserving one-to-one correspondence
between reduced word tableaux for $w$ and  EG-pipedreams of $w$.
\end{abstract}

\section{Introduction}

The structure of   bumpless pipedreams
was recently introduced by Lam, Lee and Shimozono \cite{LLS}
in their study of backward stable Schubert calculus.
They proved that bumpless pipedreams can generate back stable Schubert polynomials, which are polynomial representatives of the Schubert classes of the infinite flag varieties. Restricting to finite flag varieties, bumpless pipedreams  serve  as a new combinatorial
model for the  Schubert polynomials.
Moreover,   they found  that for a permutation $w$,
the coefficients $c^w_\lambda$
in the expansion of the Stanley symmetric function
\begin{align}\label{expan}
F_w(x)=\sum_{\lambda} c^w_\lambda\, s_{\lambda}(x)
\end{align}
in the basis of Schur functions $s_\lambda(x)$ can be interpreted by a specific
family of bumpless pipedreams, namely, the EG-pipedreams
 of $w$.

Stanley symmetric functions were invented  by
Stanley in his seminal paper \cite{Sta} in order to enumerate
the reduced decompositions of  permutations. The coefficients  $c^w_\lambda$, now called  the Edelman-Greene coefficients or the Stanley coefficients,  were first conjectured by Stanley \cite{Sta} and then proved by Edelman and Greene \cite{EG} to be nonnegative  integers.
More precisely, by developing the Edelman-Greene  algorithm  (also
called the Coxeter-Knuth algorithm), Edelman and Greene \cite{EG}
showed that the coefficient $c^w_\lambda$ is equal to
the number of reduced word tableaux for $w$ with shape $\lambda$,
see also Fomin and Greene \cite{Fomin-Greene}, Lam \cite{Lam}
or Stanley \cite{Sta-2}.

Lam, Lee and Shimozono \cite{LLS}   proposed the following problem.

\begin{prob}[Lam-Lee-Shimozono \mdseries{\cite[Problem 5.19]{LLS}}]\label{PR}
Find a direct shape preserving bijection between EG-pipedreams of $w$ and reduced word tableaux for $w$.
\end{prob}

In this paper, we provide  a desired bijection as asked in Problem \ref{PR}. The construction of our bijection relies  on the Edelman-Greene insertion algorithm \cite{EG}, the Little map \cite{Little} as well as two new structures
developed in this paper: the modified Lascoux-Sch\"utzenberger tree  and the Edelman-Greene tree.

Before outlining our strategy, let us recall the classical Lascoux-Sch\"utzenberger tree (LS-tree), which is an alternative approach  to explain  the Edelman-Greene coefficients. For a permutation $w$, the LS-tree of $w$ was
introduced by Lascoux and Sch\"utzenberger \cite{LS}  by utilizing the maximal transition
formula for Stanley symmetric functions. In the LS-tree of $w$, the children of a node $u$ are the results of applying maximal transitions to $u$. A combinatorial proof of  such  transition
relations was found  by Little \cite{Little}.
In the LS-tree of $w$, each leaf   is labeled by a Grassmannian permutation, whose Stanley symmetric function  is a Schur function.
Thus  the coefficient $c_{\lambda}^w$ is equal to the number of leaves in the LS-tree of $w$ whose Stanley symmetric function  is  $s_\lambda(x)$.

We introduce the modified Lascoux-Sch\"utzenberger tree (modified LS-tree) of a permutation $w$, and show that each leaf in the modified LS-tree of $w$ corresponds to a reduced word tableau for $w$. Unlike the classical LS-tree, the construction of a modified LS-tree  is based on general transition relations satisfied by Stanley symmetric functions. More precisely, in the modified LS-tree of $w$, the children of a node $u$ are the results of  applying some general (not necessarily maximal) transitions to $u$. Each leaf  in the modified LS-tree of $w$ is a dominant permutation.  It is known that for
a dominant permutation $u$, the corresponding Stanley symmetric function
$F_u(x)$    equals
a Schur function $s_{\lambda(u)}(x)$, where $\lambda(u)$ is the Lehmer code of $u$. Therefore,
the modified LS-tree of $w$ also allows us to expand $F_w(x)$ in terms of Schur functions, that is, the Edelman-Greene coefficient $c_{\lambda}^w$ is equal to the number of leaves in the modified LS-tree of $w$ whose Lehmer codes are   $\lambda$. Moreover, employing  the Little map \cite{Little}, we can associate
each leaf  in the modified LS-tree of $w$  to
a reduced word tableau for $w$.

It is worth mentioning  another difference between the modified
LS-tree  and the classical  LS-tree. Let  $w$ be a permutation on $\{1,2,\ldots, n\}$.
Since  maximal transitions may increase the size of  permutations,
there may exist  nodes in the LS-tree of $w$ which are labeled with
permutations on $\{1,2,\ldots, m\}$ with $m>n$. However,
the general transitions used in this paper do not increase the size of permutations,
that is, each node in the modified LS-tree of $w$ is also labeled with
a permutation on $\{1,2,\ldots, n\}$. Therefore, in some sense, the modified LS-tree seems to be more controllable than the LS-tree in the process of expanding a Stanley symmetric function into Schur functions.

On the other hand, we show that the EG-pipedreams of $w$ can also
 be generated as the leaves of a tree associated to $w$,
 which is isomorphic to the modified LS-tree of $w$.
Lam, Lee and Shimozono \cite{LLS} introduced an operation, called droops, on bumpless pipedreams. They showed that any bumpless pipedream of $w$ can be obtained by applying a sequence of droops to the Rothe pipedream of $w$.
By applying the droop operations, we construct a tree of bumpless pipedreams of $w$, called the Edelman-Greene tree (EG-tree) of $w$.
In the EG-tree of $w$, each node is a bumpless pipedream of $w$, and the children of the node $u$ are obtained by applying some specific droops to $u$, which correspond to the general  transitions in the process of constructing the modified LS-tree of $w$. Thus the EG-tree of $w$ and the modified LS-tree of $w$ are isomorphic.
In particular,
the leaves in the EG-tree of $w$ are exactly labeled with
the EG-pipedreams of $w$. Since the leaves in the modified
  LS-tree of $w$ are in one-to-one correspondence
   with the reduced word tableaux for $w$,
  we obtain a bijection  between
reduced word tableaux for $w$ and EG-pipedreams of $w$.

This paper is organized as follows. In Section \ref{preli},
we give  overviews of  the Stanley symmetric function,
 the Edelman-Greene insertion algorithm,
 the Lascoux-Sch\"utzenberger tree and the Little map.
In Section \ref{PP-2}, we describe  the structure of bumpless pipedreams
as well as the droop operation  introduced by Lam, Lee and Shimozono
\cite{LLS}.
In Section \ref{PP-5}, we introduce the structure  of a
modified  LS-tree.
Based on Section \ref{PP-5},  we construct the structure of an EG-tree
from a modified LS-tree in Section \ref{PP-6}.
 In Section \ref{PP-8}, using  the  structures
 developed in Sections \ref{PP-5} and \ref{PP-6}
 together with the Edelman-Greene algorithm and the
 Little map, we establish a shape preserving
 bijection between reduced word tableaux and
 EG-pipedreams.

\section{Preliminaries}\label{preli}

In this section, we collect some notions and structures that we are concerned
 with in this paper, including the Stanley symmetric
function, the Edelman-Greene insertion algorithm,
the Lascoux-Sch\"utzenberger tree  and the Little map.
The reader familiar with these structures could skip this section.

\subsection{Stanley symmetric functions}\label{2-1}

Let $w=w_1w_2\cdots w_n\in S_n$ be a permutation on $\{1,2,\ldots,n\}$. As usual, let $s_i$ denote the simple transposition interchanging the
elements $i$ and $i+1$.
Note that $ws_i$ is the permutation obtained from $w$ by
swapping $w_i$ and $w_{i+1}$.
A decomposition of $w$ as a product of simple transpositions is called reduced if it consists of a minimum number of simple transpositions.
A sequence $(a_1,a_2,\ldots, a_{\ell})$ is called a reduced word
of $w$ if $s_{a_1}s_{a_2}\cdots s_{a_\ell}$ is a reduced decomposition of $w$.
The length of $w$, denoted $\ell(w)$, is the number of simple transpositions in a reduced decomposition of $w$.

 Stanley symmetric functions were introduced
by Stanley \cite{Sta} to enumerate the reduced
decompositions of  a permutation.
For a permutation $w$, the Stanley symmetric function
$F_w(x)$  is defined as
\begin{equation}\label{stanley-sym-fun}
F_w(x)=\sum_{(a_1,a_2,\ldots, a_{\ell})}\sum_{1\le b_1\le b_2\le\cdots\le b_{\ell}\atop a_i<a_{i+1}\Rightarrow b_i<b_{i+1}}x_{b_1}x_{b_2}\cdots x_{b_{\ell}},
\end{equation}
where $(a_1,a_2,\ldots, a_{\ell})$ ranges over the reduced words of $w$.
Note that the definition of $F_w(x)$ in \eqref{stanley-sym-fun} is $F_{w^{-1}}(x)$ in Stanley's original  definition  \cite{Sta}.

It is well known that $F_w(x)$ can be regarded as the stable
limit of the Schubert polynomial $\S_{w}(x)$,
see for example \cite{BJS,BP,Mac,Marberg}. Let us recall the definition of
a double Schubert polynomial $\S_w(x;y)$, which reduces to $\S_{w}(x)$
by setting $y_i=0$ and will also be used in Section \ref{PP-2}.
Double Schubert polynomials were  introduced by
Lascoux and Sch\"{u}tzenberger \cite{LS1} as polynomial
representatives of the $T$-equivariant classes for Schubert varieties in the flag manifold, which can be defined recursively as follows.  If $w$ is  the longest permutation $w_0=n \,(n-1)\cdots  1$, then  set
\[\mathfrak{S}_{w_0}(x;y)=\prod_{i+j\leq n}(x_i-y_j).\]
If $w\neq w_0$, then choose a simple transposition $s_i$   such that $\ell(ws_i)=\ell(w)+1$, and   let
\begin{equation*}
\mathfrak{S}_{w}(x;y)=\partial_i \mathfrak{S}_{ws_i}(x;y).
\end{equation*}
Here,  $\partial_i$ is the  divided difference operator applies only to the $x$ variables. That is,
$\partial_i f=(f-s_if)/(x_i-x_{i+1}),$
where $f$ is a polynomial in $x$ and $s_if$ is obtained by interchanging $x_i$ and $x_{i+1}$ in $f$.   Setting $y_i=0$, $\mathfrak{S}_{w}(x; y)$ reduces to the (single)
Schubert polynomial $\S_{w}(x)$.

There have been several  combinatorial rules to generate $\S_w(x)$ and $\S_w(x;y)$, see, for example \cite{BB,BeSo, BJS,FS,WeYo,Wi, FK2,KnMi,KM}.
To explain the relation between $F_w(x)$ and $\S_{w}(x)$,
let  us recall the  combinatorial construction of $\S_{w}(x)$
due to Billey,   Jockusch and   Stanley \cite{BJS}:
\begin{equation}\label{schu}
\S_w(x)=\sum_{(a_1,a_2,\ldots, a_{\ell})}\sum_{1\leq b_1\le b_2\le\cdots\le b_{\ell}\atop
{ b_i\le a_i\atop{ a_i<a_{i+1}\Rightarrow b_i<b_{i+1}}}}x_{b_1}x_{b_2}\cdots x_{b_{\ell}},
\end{equation}
where $(a_1,a_2,\ldots, a_{\ell})$ ranges over the reduced words of $w$.
In view of \eqref{stanley-sym-fun} and \eqref{schu}, it is clear that
\begin{equation}\label{alnb}
F_w(x)=\lim_{m\rightarrow\infty}\S_{1^m\times w}(x),
\end{equation}
where $1^m\times w$ is the permutation
$12 \cdots m(w_1+m) \cdots (w_n+m).$

\subsection{Edelman-Greene insertion algorithm}


The Edelman-Greene algorithm \cite{EG} inserts a positive integer into
an increasing tableau to yield a new
increasing tableau. An increasing tableau of shape $\lambda$
is a filling of positive integers into the boxes of $\lambda$
such that the entries in each row and each column are  strictly increasing, see for example \cite{Buch,Chen,TY}.
It is worth mentioning that the Hecke insertion algorithm developed by Buch,  Kresch, Shimozono, Tamvakis and Yong  \cite{Buch} specializes to the Edelman-Greene algorithm
when applying to the reduced words of permutations.

Given a reduced word $a=(a_1,a_2,\ldots, a_\ell)$
of a permutation $w$,   the
 Edelman-Greene algorithm transforms $a$  into
a pair $(P(a), Q(a))$ of tableaux of the same shape. The tableaux $P(a)$ and $Q(a)$ are called the insertion tableau and  the recording tableau, respectively.
The tableau $P(a)$ is an increasing tableau which can be
obtained by constructing a sequence of increasing tableaux
$P_1,\ldots, P_\ell$ as follows. Let $P_1$ be the tableau with one box filled
with $a_1$. Suppose that $P_i$ ($1\leq i\leq \ell-1$) has been constructed.
Let us generate $P_{i+1}$ by inserting the integer $a_{i+1}$
into $P_i$. Set  $x=a_{i+1}$.
Let $R$ be the first row of $P_i$.
Roughly speaking, an element in $R$ may be bumped out and then inserted into the
next row. The process is repeated until no element is bumped out. There are two cases.

Case 1: The integer $x$ is strictly larger than  all the entries in $R$.
Let $P_{i+1}$ be the tableau by  adding $x$ as a new
box to the end of $R$, and the process terminates.


Case 2:
The integer $x$ is strictly smaller than some element in $R$. Let $y$ be the leftmost
entry in $R$ that is strictly larger than $x$. If replacing $y$ by $x$ results in an increasing tableau, then $y$ is bumped out by $x$ and $y$ will be inserted into the next row. If replacing
$y$ by $x$ does not result in an increasing tableau, then keep the row $R$ unchanged and the element $y$ will   be inserted into the next row.

Iterating the above procedure, we finally get the tableau $P_{i+1}$.
Set $P(a)=P_\ell$.
It should be noted that in the process of  inserting an integer into a row of $P_i$,  we will not encounter the situation  that this integer  is equal to the largest element in that row.
The tableau $Q(a)$ is the standard Young tableau which records the changes of the shapes of $P_1,P_2,\ldots,P_{\ell}$ as the insertion is performed.

\begin{figure}[h]
\begin{center}
\begin{tabular}{cccc}


\begin{tikzpicture}[scale = 1]
\def\rectanglepath{-- +(4mm,0mm) -- +(4mm,4mm) -- +(0mm,4mm) -- cycle}

\node at (5mm,30mm) {$P_i:$};
\node at (5mm,10mm) {$Q_i:$};
\draw (18mm,28mm)
\rectanglepath;
\draw (18mm,8mm) \rectanglepath;
\node at (20mm,30mm) {2};
\node at (20mm,10mm) {1};
\draw (36mm,28mm)\rectanglepath;
\draw (40mm,28mm)\rectanglepath;
\node at (38mm,30mm) {2};
\node at (42mm,30mm) {3};
\draw (36mm,8mm)\rectanglepath;
\draw (40mm,8mm)\rectanglepath;
\node at (38mm,10mm) {1};
\node at (42mm,10mm) {2};
\draw (56mm,30mm)\rectanglepath;
\draw (56mm,26mm)\rectanglepath;
\draw (60mm,30mm)\rectanglepath;
\node at (58mm,32mm) {1};
\node at (58mm,28mm) {2};
\node at (62mm,32mm) {3};
\draw (56mm,10mm)\rectanglepath;
\draw (56mm,6mm)\rectanglepath;
\draw (60mm,10mm)\rectanglepath;
\node at (58mm,12mm) {1};
\node at (58mm,8mm) {3};
\node at (62mm,12mm) {2};

\draw (74mm,30mm)\rectanglepath;
\draw (74mm,26mm)\rectanglepath;
\draw (78mm,30mm)\rectanglepath;
\draw (82mm,30mm)\rectanglepath;
\node at (76mm,32mm) {1};
\node at (76mm,28mm) {2};
\node at (80mm,32mm) {3};
\node at (84mm,32mm) {6};
\draw (74mm,10mm)\rectanglepath;
\draw (74mm,6mm)\rectanglepath;
\draw (78mm,10mm)\rectanglepath;
\draw (82mm,10mm)\rectanglepath;
\node at (76mm,12mm) {1};
\node at (76mm,8mm) {3};
\node at (80mm,12mm) {2};
\node at (84mm,12mm) {4};

\draw (94mm,30mm)\rectanglepath;
\draw (94mm,26mm)\rectanglepath;
\draw (98mm,30mm)\rectanglepath;
\draw (98mm,26mm)\rectanglepath;
\draw (102mm,30mm)\rectanglepath;
\node at (96mm,32mm) {1};
\node at (96mm,28mm) {2};
\node at (100mm,32mm) {3};
\node at (100mm,28mm) {6};
\node at (104mm,32mm) {4};
\draw (94mm,10mm)\rectanglepath;
\draw (94mm,6mm)\rectanglepath;
\draw (98mm,10mm)\rectanglepath;
\draw (98mm,6mm)\rectanglepath;
\draw (102mm,10mm)\rectanglepath;
\node at (96mm,12mm) {1};
\node at (96mm,8mm) {3};
\node at (100mm,12mm) {2};
\node at (100mm,8mm) {5};
\node at (104mm,12mm) {4};

\draw (114mm,32mm)\rectanglepath;
\draw (114mm,28mm)\rectanglepath;
\draw (114mm,24mm)\rectanglepath;
\draw (118mm,32mm)\rectanglepath;
\draw (118mm,28mm)\rectanglepath;
\draw (122mm,32mm)\rectanglepath;
\node at (116mm,34mm) {\small{1}};
\node at (116mm,30mm) {\small{2}};
\node at (116mm,26mm) {\small{6}};
\node at (120mm,34mm) {\small{3}};
\node at (120mm,30mm) {\small{4}};
\node at (124mm,34mm) {\small{4}};
\draw (114mm,12mm)\rectanglepath;
\draw (114mm,8mm)\rectanglepath;
\draw (114mm,4mm)\rectanglepath;
\draw (118mm,12mm)\rectanglepath;
\draw (118mm,8mm)\rectanglepath;
\draw (122mm,12mm)\rectanglepath;
\node at (116mm,14mm) {\small{1}};
\node at (116mm,10mm) {\small{3}};
\node at (116mm,6mm) {\small{6}};
\node at (120mm,14mm) {\small{2}};
\node at (120mm,10mm) {\small{5}};
\node at (124mm,14mm) {\small{4}};

\draw (134mm,34mm)\rectanglepath;
\draw (134mm,30mm)\rectanglepath;
\draw (134mm,26mm)\rectanglepath;
\draw (134mm,22mm)\rectanglepath;
\draw (138mm,34mm)\rectanglepath;
\draw (138mm,30mm)\rectanglepath;
\draw (142mm,34mm)\rectanglepath;
\node at (136mm,36mm) {\small{1}};
\node at (136mm,32mm) {\small{2}};
\node at (136mm,28mm) {\small{4}};
\node at (136mm,24mm) {\small{6}};
\node at (140mm,36mm) {\small{2}};
\node at (140mm,32mm) {\small{3}};
\node at (144mm,36mm) {\small{4}};
\draw (134mm,14mm)\rectanglepath;
\draw (134mm,10mm)\rectanglepath;
\draw (134mm,6mm)\rectanglepath;
\draw (134mm,2mm)\rectanglepath;
\draw (138mm,14mm)\rectanglepath;
\draw (138mm,10mm)\rectanglepath;
\draw (142mm,14mm)\rectanglepath;
\node at (136mm,16mm) {\small{1}};
\node at (136mm,12mm) {\small{3}};
\node at (136mm,8mm) {\small{6}};
\node at (136mm,4mm) {\small{7}};
\node at (140mm,16mm) {\small{2}};
\node at (140mm,12mm) {\small{5}};
\node at (144mm,16mm) {\small{4}};
\end{tikzpicture}
\end{tabular}
\caption{Construction of $P(a)$ and $Q(a)$ for $a=(2,3,1,6,4,3,2)$.}
\label{EG-insertion}
\end{center}
\end{figure}
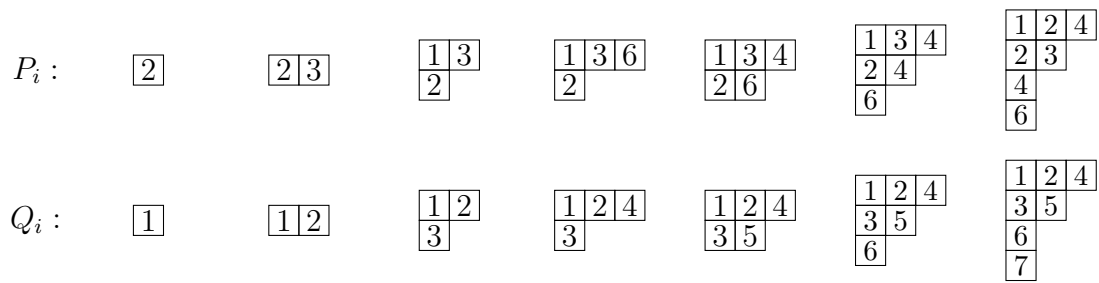

For example, let $a=(2,3,1,6,4,3,2)$ be a reduced word
of the permutation $w=3514276$. Then $P(a)$ and $Q(a)$ can be constructed as in
Figure \ref{EG-insertion}.

To characterize the
insertion tableau $P(a)$, we need the notion of  the row
reading word  of a  tableau $T$,
that is, the word obtained
by reading the entries of $T$ along the rows from left to right,
bottom to top. For example, the tableau $P(a)$
in Figure \ref{EG-insertion} has row reading word $(6,4,2,3,1,2,4)$.

\begin{theo}[Edelman-Greene \mdseries{\cite{EG}}]\label{90}
The Edelman-Greene correspondence is a bijection between
the set of reduced words of $w$ and the set of
pairs $(P,Q)$ of tableaux of the same shape,
where $P$ is an increasing tableaux whose row reading
word is a reduced word of $w$,
and $Q$ is a standard Young tableau.
\end{theo}

The following theorem is attributed to Edelman and Greene \cite{EG}, see also  Fomin and Greene \cite{Fomin-Greene}, Lam \cite{Lam} or Stanley \cite{Sta-2}.

\begin{theo}[Edelman-Greene \mdseries{\cite{EG}}]\label{PO}
The coefficient $c^{w}_\lambda$ equals  the number of increasing
  tableaux of shape $\lambda$ whose row reading words are reduced words of $w^{-1}$.
\end{theo}

Theorem \ref{PO} has an equivalent
statement in terms of column reading words.
The column reading word
 of $T$, denoted $\mathrm{column}(T)$, is obtained by
reading the entries of $T$ along the columns from top to bottom, right to left.

The following theorem follows from  \cite[Theorem 1]{Buch} by restricting a Hecke word to a reduced word.
Here we give a self-contained proof based on
properties of the Edelman-Greene algorithm.

\begin{theo}[Buch-Kresch-Shimozono-Tamvakis-Yong \mdseries{\cite{Buch}}]\label{91}
The coefficient $c^{w}_\lambda$ equals the number of increasing  tableaux of shape $\lambda$ whose column reading words are reduced words of $w$.
\end{theo}

\pf
Let $T$ be an increasing tableau such that the row reading word of $T$ is a reduced word of $w^{-1}$.
We claim that $\mathrm{column}(T)$ is a reduced word of
$w$. Fix a standard Young tableau $Q$ which has the same shape as $T$. Let $a=(a_1,a_2,\ldots,a_\ell)$
be the reduced word of $w^{-1}$
corresponding to the pair $(T, Q)$.  Denote $T^t$ by the transpose of
$T$, and write
\begin{align}\label{rvs}
a^{\mathrm{rev}}=(a_\ell, a_{\ell-1},\ldots, a_1)
\end{align}
for the reverse of $a$. Note that $a^{\mathrm{rev}}$
is a reduced word of $w$.
By Edelman and Greene \cite{EG} (see also
Felsner \cite{Felsner}), the insertion tableau of $a^{\mathrm{rev}}$
is $T^t$. By Theorem \ref{90}, the row reading word of $T^t$ is
a reduced word of $w$, or equivalently, the column reading word of $T$
is a reduced word of $w$. This verifies the claim.

Conversely, we can show that if $T$ is an
increasing tableau whose
column reading word is a reduced word of $w$, then
its row reading word is a reduced word of $w^{-1}$. This
completes the proof.
\qed

Throughout this paper, an increasing tableau $T$ is called
a reduced word tableau for $w$ if $\mathrm{column}(T)$ is a reduced word of $w$.
Equivalently, by the proof of Theorem \ref{91}, an increasing tableau
$T$ is a reduced word tableau for $w$ if its row
reading word is a reduced word of
$w^{-1}$.

\subsection{Lascoux-Sch\"utzenberger tree}\label{lstree}

The LS-tree is built based on maximal transitions
   on Stanley symmetric functions.
For a  permutation $w=w_1w_2\cdots w_n$,
let
\begin{align*}\label{define-r-s}
r&=\max \{i\,|\,w_i>w_{i+1}\},\\[5pt]
s&=\max \{i>r\,|\,w_i<w_{r}\},\\[5pt]
I(w)&=\{i<r\,|\,w_i<w_s\ \ \text{and}\ \ \forall \,j\in(i,r),\,w_j\not\in(w_i,w_s)\},
\end{align*}
where for  two integers $a$ and $b$ with $a<b$, we use $(a,b)$ to denote the
interval $\{a+1,a+2,\ldots, b-1\}$.
In other words, $r$ is the last decent of $w$,  $s$ is the unique position after $r$ such that $\ell(w t_{r,s})=\ell(w)-1$, where $wt_{r,s}$ is the permutation obtained by swapping $w_r$ and $w_s$. And $I(w)$ is the set of positions $i$ before $r$ such that $\ell(w t_{r,s} t_{i,r})=\ell(w)$.
With the above notation, set
\[\Phi(w)=\left\{
            \begin{array}{ll}
              \{w t_{r,s}t_{i,r}\,|\, i\in I(w)\}, & \hbox{if $I(w)\neq\emptyset$}; \\
              \Phi(1\times w), & \hbox{if $I(w)=\emptyset$,}
            \end{array}
          \right.
\]
where $1\times w=1(w_1+1)\cdots (w_n+1)$ as defined  in Section \ref{2-1}.

\begin{figure}[!htb]
\begin{center}
\begin{tikzpicture}

\draw (50mm,72mm) node{$2\,\underline{3}\,\underline{1}\,6\,\boxed{5}\,\bar{4}$};

\draw (30mm,54mm) node{$2\,\underline{4}\,\underline{1}\,\boxed{6}\,3\,\bar{5}$};

\draw (70mm,54mm) node{$2\,3\,\underline{4}\,\boxed{6}\,1\,\bar{5}$};


\draw (15mm,36mm)
node{$\underline{2}\,5\,\underline{1}\,\boxed{4}\,\bar{3}\,6$};

\draw (45mm,36mm)
node{$\underline{2}\,4\,\boxed{5}\,1\,\bar{3}\,6$};

\draw (73mm,36mm) node{$2\,3\,5\,\boxed{4}\,\bar{1}\,6$};
\draw (73mm,30mm) node{$\underline{1}\,3\,4\,6\,\boxed{5}\,\bar{2}\,7$};

%
\draw (-5mm,18mm) node{$3\,5\,1\,2\,4\,6$};

\draw (25mm,18mm) node{$2\,5\,\boxed{3}\,\bar{1}\,4\,6$};
\draw (25mm,12mm) node{$\underline{1}\,3\,6\,\boxed{4}\,\bar{2}\,5\,7$};

\draw (48mm,18mm)
node{$3\,4\,\boxed{2}\,\bar{1}\,5\,6$};
\draw (48mm,12mm)
node{$\underline{1}\,4\,5\,\boxed{3}\,\bar{2}\,6\,7$};

\draw (73mm,18mm) node{$2\,3\,4\,6\,1\,5\,7$};

\draw (25mm,0mm) node{$2\,3\,6\,1\,4\,5\,7$};

\draw (48mm,0mm) node{$2\,4\,5\,1\,3\,6\,7$};

\draw(50mm,67mm)--(70mm,59mm);
\draw(50mm,67mm)--(30mm,59mm);
\draw(30mm,50mm)--(15mm,40mm);
\draw(30mm,50mm)--(45mm,40mm);
\draw(73mm,50mm)--(73mm,40mm);
\draw(10mm,32mm)--(-5mm,22mm);
\draw(10mm,32mm)--(25mm,22mm);
\draw(48mm,32mm)--(48mm,22mm);
\draw(73mm,27mm)--(73mm,22mm);

\draw(25mm,8mm)--(25mm,3mm);

\draw(48mm,8mm)--(48mm,3mm);

\end{tikzpicture}
\caption{The LS-tree for $w=231654$.}\label{LS-tree}
\end{center}
\end{figure}
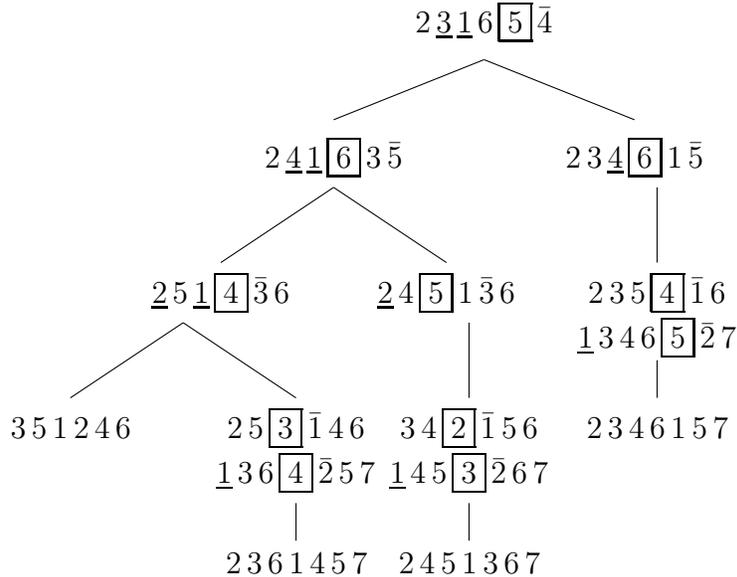

Each permutation in the set $\Phi(w)$ is called a child of $w$.
The  LS-tree of $w$ is obtained by
recursively applying the above operation
until every leaf is a Grassmanian permutation.
A Grassmanian permutation is a permutation
 with at most one descent.  In Figure \ref{LS-tree},
we  illustrate the LS-tree
 of $w=231654$, where,
 as  used  in \cite{Little}, the entries $w_r$,
 $w_s$ and $w_i$ with $i\in I(w)$ are boxed,
 barred and underlined, respectively.

According to the  transition formula for Schubert polynomials
due to Lascoux and Sch\"utzenberger \cite{LS}, one has the
following transition relation for Stanley symmetric functions:
\begin{equation}\label{StanRe}
F_w(x)=\sum_{w'\in \Phi(w)}F_{w'}(x),
\end{equation}
see also \cite[Theorem 1]{Little} or Garsia \cite{Garsia}.
Since each leaf of a LS-tree is labeled
with a Grassmannian permutation,
  $F_w(x)$ can be eventually
written as a sum of Stanley symmetric functions
indexed by Grassmannian permutations.
For a Grassmannian permutation, the Stanley symmetric function
becomes a single
Schur function. More generally,
as proved by Stanley \cite{Sta},
the Stanley symmetric function for
a vexillary permutation is a single Schur function. Since Grassmannian permutations are  vexillary,  the LS-tree implies the Schur positivity of
Stanley symmetric functions.

Recall that a permutation
is called a vexillary (or 2143-avoiding) permutation if it does not contain a  subsequence order-isormorphic to
2143, that is, there are no indices $i_1<i_2<i_3<i_4$ such that
$w_{i_2}<w_{i_1}<w_{i_4}<w_{i_3}$. For  a permutation $w\in S_n$,
let
\[c(w)=(c_1(w),c_2(w),\ldots, c_n(w))\]
be the Lehmer code of $w$, where
$c_i(w)=|\{j\,|\, j>i, w_j<w_i\}|.$
Define a partition $\lambda(w)$ by   rearranging
the Lehmer code of $w$ in weakly decreasing order.
For example, the Lehmer code of  $w=35412$ is $(2,3,2,0,0)$, and so
we have $\lambda(w)=(3,2,2)$. For a vexillary permutation
$w$, Stanley \cite{Sta} proved that
$F_w(x)=s_{\lambda(w)}(x),$
which  also follows from \eqref{alnb}
together with a tableau formula for
Schubert polynomials of   vexillary permutations
due to  Wachs  \cite{Wa} or Knutson, Miller and Yong
\cite{KMY}.

In particular, a permutation $w$ is called a dominant permutation if $w$ is 132-avoiding.  It is clear that a dominant permutation
is also a vexillary permutation.
 The Lehmer code of a dominant permutation $w$
is a weakly decreasing sequence which form a partition shape $\lambda(w)$, see Stanley \cite[Chapter 1]{Sta-5}. Hence, for a dominant permutation $w$, we have
$F_{w}(x)=s_{\lambda(w)}(x),$
where $\lambda(w)=c(w)$ is the partition equal to the Lehmer code of $w$.
Therefore, there is only one reduced word tableau for $w$.
In fact,  the only  reduced word
tableau $T(w)$ for $w$ can be obtained
as follows: fill the entry in the box $(i,j)$ of $\lambda(w)$
 with $i+j-1$.
The tableau
$T(w)$ is also called a frozen tableau, see \cite{Linusson}.

There is a more general relation satisfied by
Stanley symmetric functions,
including \eqref{StanRe} as a special case.
For a permutation  $u=u_1u_2\cdots u_n$ and $1\leq k\leq n$, let
\begin{equation}\label{LP}
I(u,k)=\{i<k\,|\,\ell(ut_{i,k})=\ell(u)+1\},
\end{equation}
and
\begin{align}\label{snks}
S(u,k)=\{j>k\,|\,\ell(ut_{k,j})=\ell(u)+1\}.
\end{align}
Define two sets of permutations by
\begin{align}\label{bigphi}
\Phi(u,k)&=
\left\{
  \begin{array}{ll}
    \{ut_{i,k}\,|\,i\in I(u,k)\}, & \hbox{if $I(u,k)\neq\emptyset$;} \\
    \Phi(1\times u,k+1), & \hbox{otherwise},
  \end{array}
\right.
\end{align}
and
\begin{align}\label{bigpsi}
\Psi(u,k)&=
\left\{
  \begin{array}{ll}
    \{ut_{k,j}\,|\,j\in S(u,k)\}, & \hbox{if $S(u,k)\neq\emptyset$;} \\
    \Psi(u\times 1,k), & \hbox{otherwise},
  \end{array}
\right.
\end{align}
where $u\times 1=u_1u_2\cdots u_n(n+1)$.

By the Monk's rule for Schubert polynomials (see for example
\cite{BB,Mac}) together with \eqref{alnb}, one can easily establish the following relation:
\begin{equation}\label{StanRe2}
    \sum_{w\in \Psi(u,k)}F_{w}(x)=\sum_{w\in\Phi(u,k)}F_{w}(x).
  \end{equation}
Indeed, \eqref{StanRe2} contains  \eqref{StanRe} as a special case,
since 
\[\Phi(wt_{r,s},r)=\Phi(w)\ \ \text{and}\ \  \Psi(wt_{r,s},r)=\{w\}.\]

\subsection{Little map}

Little \cite{Little} developed  a bijection,
known as the Little map, to give a combinatorial proof of \eqref{StanRe2}.
 Let $\mathrm{Red}(w)$ denote  the set of reduced words   of a permutation $w$.
The Little map  is a descent preserving bijection:
\begin{align}\label{thetar}
\theta_k:\ \ \bigcup_{w\in \Psi(u,k)}\mathrm{Red}(w)\ \ \longrightarrow \ \bigcup_{w'\in \Phi(u,k)}\mathrm{Red}(w').
\end{align}
Indeed, in view of  \eqref{stanley-sym-fun},  the Little map yields
 a combinatorial  proof of \eqref{StanRe2}.

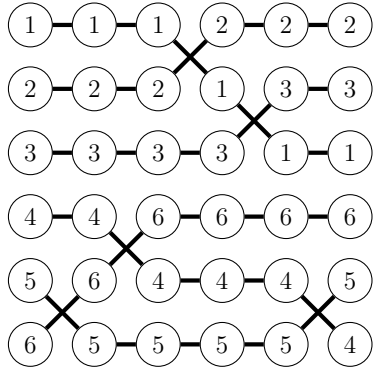
\begin{figure}[h]
\begin{center}
\begin{tabular}{ccc}
\begin{tikzpicture}[scale = .85]
\node[draw,circle,scale=0.8] at
(0,0) (0/1) {1};
\node[draw,circle,scale=0.8] at (0,-1) (0/2) {2};
\node[draw,circle,scale=0.8] at (0,-2) (0/3) {3};
\node[draw,circle,scale=0.8] at (0,-3) (0/4) {4};
\node[draw,circle,scale=0.8] at (0,-4) (0/5) {5};
\node[draw,circle,scale=0.8] at (0,-5) (0/6) {6};

\node[draw,circle,scale=0.8] at
(1,0) (1/1) {1};
\node[draw,circle,scale=0.8] at (1,-1) (1/2) {2};
\node[draw,circle,scale=0.8] at (1,-2) (1/3) {3};
\node[draw,circle,scale=0.8] at (1,-3) (1/4) {4};
\node[draw,circle,scale=0.8] at (1,-4) (1/5) {6};
\node[draw,circle,scale=0.8] at (1,-5) (1/6) {5};

\node[draw,circle,scale=0.8] at
(2,0) (2/1) {1};
\node[draw,circle,scale=0.8] at (2,-1) (2/2) {2};
\node[draw,circle,scale=0.8] at (2,-2) (2/3) {3};
\node[draw,circle,scale=0.8] at (2,-3) (2/4) {6};
\node[draw,circle,scale=0.8] at (2,-4) (2/5) {4};
\node[draw,circle,scale=0.8] at (2,-5) (2/6) {5};

\node[draw,circle,scale=0.8] at
(3,0) (3/1) {2};
\node[draw,circle,scale=0.8] at (3,-1) (3/2) {1};
\node[draw,circle,scale=0.8] at (3,-2) (3/3) {3};
\node[draw,circle,scale=0.8] at (3,-3) (3/4) {6};
\node[draw,circle,scale=0.8] at (3,-4) (3/5) {4};
\node[draw,circle,scale=0.8] at (3,-5) (3/6) {5};

\node[draw,circle,scale=0.8] at (4,0) (4/1) {2};
\node[draw,circle,scale=0.8] at (4,-1) (4/2) {3};
\node[draw,circle,scale=0.8] at (4,-2) (4/3) {1};
\node[draw,circle,scale=0.8] at (4,-3) (4/4) {6};
\node[draw,circle,scale=0.8] at (4,-4) (4/5) {4};
\node[draw,circle,scale=0.8] at (4,-5) (4/6) {5};

\node[draw,circle,scale=0.8] at
(5,0) (5/1) {2};
\node[draw,circle,scale=0.8] at (5,-1) (5/2) {3};
\node[draw,circle,scale=0.8] at (5,-2) (5/3) {1};
\node[draw,circle,scale=0.8] at (5,-3) (5/4) {6};
\node[draw,circle,scale=0.8] at (5,-4) (5/5) {5};
\node[draw,circle,scale=0.8] at (5,-5) (5/6) {4};

\draw[ultra thick] (0/1) -- (1/1) -- (2/1) -- (3/2) -- (4/3) -- (5/3);
\draw[ultra thick] (0/2) -- (1/2) -- (2/2) -- (3/1) -- (4/1) -- (5/1);
\draw[ultra thick] (0/3) -- (1/3) -- (2/3) -- (3/3) -- (4/2) -- (5/2);
\draw[ultra thick] (0/4) -- (1/4) -- (2/5) -- (3/5) -- (4/5)--(5/6);
\draw[ultra thick] (0/5) -- (1/6) -- (2/6) -- (3/6) -- (4/6) -- (5/5);
\draw[ultra thick] (0/6) -- (1/5) -- (2/4) -- (3/4) -- (4/4) -- (5/4);
\end{tikzpicture}

\end{tabular}
\end{center}
\vspace{-.5cm}
\caption{The line diagram of $a=(5,4,1,2,5)$.}\label{fig-line}
\end{figure}

The Little map   is defined based on a bumping algorithm
acting on the line diagrams for permutations.
Assume that $a=(a_1,a_2,\ldots, a_\ell)$ is a   word
(not necessarily reduced) of a permutation $w$.
The line diagram  of   $a$ is the array $\{1,2,
\ldots, \ell\}\times \{1,2,\ldots,n\}$ in the Cartesian coordinates,
which describes the trajectories of the numbers
$1,2,\ldots,n$ as they are arranged into the permutation $w$
by successive simple transpositions. Note that  $a$
is reduced if and only if no two lines cross more than
once.
For example, Figure \ref{fig-line}  is the  line diagram of the
 reduced word $a=(5,4,1,2,5)$ of $w = 231654$.

Given a  word $a=(a_1,a_2,\ldots, a_m)$ and an integer $1\leq t\leq m$,
let
\[a^{(t)}=(a_1,a_2,\ldots, a_{t-1},a_{t+1},\ldots, a_m)\]
 and define
\[\mathclap{a}\uparrow_t\,=\left\{
\begin{array}{ll}
(a_1,a_2,\ldots, a_{t-1},a_t-1,a_{t+1},\ldots, a_m), &\hbox{if $a_t>1$};\\[5pt]
(a_1+1, ,a_2+1,\ldots,a_{t-1}+1, a_t, a_{t+1}+1, \ldots,
a_m+1), &\hbox{if $a_t=1$},
\end{array}
\right.\]
where \  $\mathclap{a}\uparrow_t$ is called  the word obtained from $a$ by bumping at time $t$.
Denote the word obtained after a sequence of  bumps by
\[\mathclap{a}\uparrow_{t_1,t_2,\ldots ,t_i}=(((\ \mathclap{a}\uparrow_{t_1})\uparrow_{t_2})\cdots)\uparrow_{t_i}.\]

The bumping algorithm, called the Little bump,
transforms a reduced word into another reduced word,
which can be sketched as follows. For more detailed information, see
Little \cite{Little} or Hamaker and Young \cite{HY}.
Let  $a=(a_1,a_2,\ldots, a_m)$ be a reduced word.
Assume that $t_1$ is an index such that $a^{(t_1)}$ is also reduced.
Consider the word $\, $  $\mathclap{a}\uparrow_{t_1}$. If  \ $\mathclap{a}\uparrow_{t_1}$ is a reduced word,
then the algorithm terminates.
Otherwise, it can be shown that there is a unique index,
say $t_2$, such that\ $(\ \mathclap{a}\uparrow_{t_1})^{(t_2)}$ is reduced.
Now consider the word $\, $ $\mathclap{a}\uparrow_{t_1, t_2}$. Repeating the above procedure,
one is  eventually left with a reduced word \  $\mathclap{a}\uparrow_{t_1,t_2,\ldots ,t_i}$. The above process is referred to
as the Little bump.

Using the  above bumping algorithm, we can
 define the Little map $\theta_k$.
Given a permutation $w\in \Psi(u,k)$, let $a$ be a reduced word of
$w$. By the definition of the set $\Psi(u,k)$, there is
an index $j>k$ in $S(u,k)$ such that $w=u\, t_{k,j}$ and
$\ell(w)=\ell(u)+1$. Let us first  define a map $\theta_{k, w_j}$ as follows.
Since $k<j$, $w_k>w_j$ and $a$ is a reduced word, there is exactly
 one letter, say $a_{t_1}$, that interchanges   $w_k$ and $w_j$.
Because  $\ell(w)=\ell(u)+1$, the word $a^{(t_1)}$ is reduced.
Thus we can invoke the Little bump beginning at the
position $t_1$. Let $a'$ be the reduced word after applying
the Little bump.  
Define
\[\theta_k(a)=\theta_{k, w_j}(a)=a'.\]
Little \cite{Little} showed  that
\begin{align}\label{usf}
a'\in \mathrm{Red}(w')\  \text{for some} \  w'\in\Phi(u, k).
\end{align}
We remark that although the subscript $w_j$ in $\theta_{k, w_j}$ is determined once
a reduce word   of $w\in \Psi(u,k)$
is given, we would like to use the two parameters $k$ and $w_j$
since this would be more convenient to describe the inverse of the  Little map.

\begin{figure}[h]
\begin{tabular}{cccc}
\begin{tikzpicture}[scale = .62]
\node[draw,circle,scale=0.75] at
(0,0) (0/1) {1};
\node[draw,circle,scale=0.75] at (0,-1) (0/2) {2};
\node[draw,circle,scale=0.75] at (0,-2) (0/3) {3};
\node[draw,circle,scale=0.75] at (0,-3) (0/4) {4};
\node[draw,circle,scale=0.75] at (0,-4) (0/5) {5};
\node[draw,circle,scale=0.75] at (0,-5) (0/6) {6};

\node[draw,circle,scale=0.75] at (1,0) (1/1) {1};
\node[draw,circle,scale=0.75] at (1,-1) (1/2) {2};
\node[draw,circle,scale=0.75] at (1,-2) (1/3) {4};
\node[draw,circle,scale=0.75] at (1,-3) (1/4) {3};
\node[draw,circle,scale=0.75] at (1,-4) (1/5) {5};
\node[draw,circle,scale=0.75] at (1,-5) (1/6) {6};

\node[draw,circle,scale=0.75] at
(2,0) (2/1) {2};
\node[draw,circle,scale=0.75] at (2,-1) (2/2) {1};
\node[draw,circle,scale=0.75] at (2,-2) (2/3) {4};
\node[draw,circle,scale=0.75] at (2,-3) (2/4) {3};
\node[draw,circle,scale=0.75] at (2,-4) (2/5) {5};
\node[draw,circle,scale=0.75] at (2,-5) (2/6) {6};

\node[draw,circle,scale=0.75] at
(3,0) (3/1) {2};
\node[draw,circle,scale=0.75] at (3,-1) (3/2) {1};
\node[draw,circle,scale=0.75] at (3,-2) (3/3) {4};
\node[draw,circle,scale=0.75] at (3,-3) (3/4) {5};
\node[draw,circle,scale=0.75] at (3,-4) (3/5) {3};
\node[draw,circle,scale=0.75] at (3,-5) (3/6) {6};

\node[draw,circle,scale=0.75] at (4,0) (4/1) {2};
\node[draw,circle,scale=0.75] at (4,-1) (4/2) {1};
\node[draw,circle,scale=0.75] at (4,-2) (4/3) {4};
\node[draw,circle,scale=0.75] at (4,-3) (4/4) {5};
\node[draw,circle,scale=0.75] at (4,-4) (4/5) {6};
\node[draw,circle,scale=0.75] at (4,-5) (4/6) {3};

\node[draw,circle,scale=0.75] at
(5,0) (5/1) {2};
\node[draw,circle,scale=0.75] at (5,-1) (5/2) {4};
\node[draw,circle,scale=0.75] at (5,-2) (5/3) {1};
\node[draw,circle,scale=0.75] at (5,-3) (5/4) {5};
\node[draw,circle,scale=0.75] at (5,-4) (5/5) {6};
\node[draw,circle,scale=0.75] at (5,-5) (5/6) {3};

\draw[ultra thick] (0/1) -- (1/1) -- (2/2) -- (3/2) -- (4/2) -- (5/3);
\draw[ultra thick] (0/2) -- (1/2) -- (2/1) -- (3/1) -- (4/1) -- (5/1);
\draw[ultra thick] (0/3) -- (1/4) -- (2/4) -- (3/5) -- (4/6) -- (5/6);
\draw[ultra thick] (0/4) -- (1/3) -- (2/3) -- (3/3) -- (4/3)--(5/2);
\draw[ultra thick] (0/5) -- (1/5) -- (2/5) -- (3/4) -- (4/4) -- (5/4);
\draw[ultra thick] (0/6) -- (1/6) -- (2/6) -- (3/6) -- (4/5) -- (5/5);

\draw[ultra thick,dotted] (3/5) -- (4/4);\draw[ultra thick,dashed] (3/4) -- (4/5);

\end{tikzpicture}

&

\begin{tikzpicture}[scale = .62]
\node[draw,circle,scale=0.75] at
(0,0) (0/1) {1};
\node[draw,circle,scale=0.75] at (0,-1) (0/2) {2};
\node[draw,circle,scale=0.75] at (0,-2) (0/3) {3};
\node[draw,circle,scale=0.75] at (0,-3) (0/4) {4};
\node[draw,circle,scale=0.75] at (0,-4) (0/5) {5};
\node[draw,circle,scale=0.75] at (0,-5) (0/6) {6};

\node[draw,circle,scale=0.75] at
(1,0) (1/1) {1};
\node[draw,circle,scale=0.75] at (1,-1) (1/2) {2};
\node[draw,circle,scale=0.75] at (1,-2) (1/3) {4};
\node[draw,circle,scale=0.75] at (1,-3) (1/4) {3};
\node[draw,circle,scale=0.75] at (1,-4) (1/5) {5};
\node[draw,circle,scale=0.75] at (1,-5) (1/6) {6};

\node[draw,circle,scale=0.75] at
(2,0) (2/1) {2};
\node[draw,circle,scale=0.75] at (2,-1) (2/2) {1};
\node[draw,circle,scale=0.75] at (2,-2) (2/3) {4};
\node[draw,circle,scale=0.75] at (2,-3) (2/4) {3};
\node[draw,circle,scale=0.75] at (2,-4) (2/5) {5};
\node[draw,circle,scale=0.75] at (2,-5) (2/6) {6};

\node[draw,circle,scale=0.75] at (3,0) (3/1) {2};
\node[draw,circle,scale=0.75] at (3,-1) (3/2) {1};
\node[draw,circle,scale=0.75] at (3,-2) (3/3) {4};
\node[draw,circle,scale=0.75] at (3,-3) (3/4) {5};
\node[draw,circle,scale=0.75] at (3,-4) (3/5) {3};
\node[draw,circle,scale=0.75] at (3,-5) (3/6) {6};

\node[draw,circle,scale=0.75] at (4,0) (4/1) {2};
\node[draw,circle,scale=0.75] at (4,-1) (4/2) {1};
\node[draw,circle,scale=0.75] at (4,-2) (4/3) {4};
\node[draw,circle,scale=0.75] at (4,-3) (4/4) {3};
\node[draw,circle,scale=0.75] at (4,-4) (4/5) {5};
\node[draw,circle,scale=0.75] at (4,-5) (4/6) {6};

\node[draw,circle,scale=0.75] at
(5,0) (5/1) {2};
\node[draw,circle,scale=0.75] at (5,-1) (5/2) {4};
\node[draw,circle,scale=0.75] at (5,-2) (5/3) {1};
\node[draw,circle,scale=0.75] at (5,-3) (5/4) {3};
\node[draw,circle,scale=0.75] at (5,-4) (5/5) {5};
\node[draw,circle,scale=0.75] at (5,-5) (5/6) {6};

\draw[ultra thick] (0/1) -- (1/1) -- (2/2) -- (3/2) -- (4/2) -- (5/3);
\draw[ultra thick] (0/2) -- (1/2) -- (2/1) -- (3/1) -- (4/1) -- (5/1);
\draw[ultra thick] (0/3) -- (1/4) -- (2/4) -- (3/5) -- (4/4) -- (5/4);
\draw[ultra thick] (0/4) -- (1/3) -- (2/3) -- (3/3) -- (4/3)--(5/2);
\draw[ultra thick] (0/5) -- (1/5) -- (2/5) -- (3/4) -- (4/5) -- (5/5);
\draw[ultra thick] (0/6) -- (1/6) -- (2/6) -- (3/6) -- (4/6) -- (5/6);

\draw[ultra thick,dotted] (2/4) -- (3/3);\draw[ultra thick,dashed] (2/3) -- (3/4);

\end{tikzpicture}
&

\begin{tikzpicture}[scale = .62]
\node[draw,circle,scale=0.75] at
(0,0) (0/1) {1};
\node[draw,circle,scale=0.75] at (0,-1) (0/2) {2};
\node[draw,circle,scale=0.75] at (0,-2) (0/3) {3};
\node[draw,circle,scale=0.75] at (0,-3) (0/4) {4};
\node[draw,circle,scale=0.75] at (0,-4) (0/5) {5};
\node[draw,circle,scale=0.75] at (0,-5) (0/6) {6};

\node[draw,circle,scale=0.75] at
(1,0) (1/1) {1};
\node[draw,circle,scale=0.75] at (1,-1) (1/2) {2};
\node[draw,circle,scale=0.75] at (1,-2) (1/3) {4};
\node[draw,circle,scale=0.75] at (1,-3) (1/4) {3};
\node[draw,circle,scale=0.75] at (1,-4) (1/5) {5};
\node[draw,circle,scale=0.75] at (1,-5) (1/6) {6};

\node[draw,circle,scale=0.75] at
(2,0) (2/1) {2};
\node[draw,circle,scale=0.75] at (2,-1) (2/2) {1};
\node[draw,circle,scale=0.75] at (2,-2) (2/3) {4};
\node[draw,circle,scale=0.75] at (2,-3) (2/4) {3};
\node[draw,circle,scale=0.75] at (2,-4) (2/5) {5};
\node[draw,circle,scale=0.75] at (2,-5) (2/6) {6};

\node[draw,circle,scale=0.75] at (3,0) (3/1) {2};
\node[draw,circle,scale=0.75] at (3,-1) (3/2) {1};
\node[draw,circle,scale=0.75] at (3,-2) (3/3) {3};
\node[draw,circle,scale=0.75] at (3,-3) (3/4) {4};
\node[draw,circle,scale=0.75] at (3,-4) (3/5) {5};
\node[draw,circle,scale=0.75] at (3,-5) (3/6) {6};

\node[draw,circle,scale=0.75] at (4,0) (4/1) {2};
\node[draw,circle,scale=0.75] at (4,-1) (4/2) {1};
\node[draw,circle,scale=0.75] at (4,-2) (4/3) {3};
\node[draw,circle,scale=0.75] at (4,-3) (4/4) {5};
\node[draw,circle,scale=0.75] at (4,-4) (4/5) {4};
\node[draw,circle,scale=0.75] at (4,-5) (4/6) {6};

\node[draw,circle,scale=0.75] at
(5,0) (5/1) {2};
\node[draw,circle,scale=0.75] at (5,-1) (5/2) {3};
\node[draw,circle,scale=0.75] at (5,-2) (5/3) {1};
\node[draw,circle,scale=0.75] at (5,-3) (5/4) {5};
\node[draw,circle,scale=0.75] at (5,-4) (5/5) {4};
\node[draw,circle,scale=0.75] at (5,-5) (5/6) {6};

\draw[ultra thick] (0/1) -- (1/1) -- (2/2) -- (3/2) -- (4/2) -- (5/3);
\draw[ultra thick] (0/2) -- (1/2) -- (2/1) -- (3/1) -- (4/1) -- (5/1);
\draw[ultra thick] (0/3) -- (1/4) -- (2/4) -- (3/3) -- (4/3) -- (5/2);
\draw[ultra thick] (0/4) -- (1/3) -- (2/3) -- (3/4) -- (4/5)--(5/5);
\draw[ultra thick] (0/5) -- (1/5) -- (2/5) -- (3/5) -- (4/4) -- (5/4);
\draw[ultra thick] (0/6) -- (1/6) -- (2/6) -- (3/6) -- (4/6) -- (5/6);

\draw[ultra thick,dotted] (0/2) -- (1/3);\draw[ultra thick,dashed] (1/2) -- (0/3);
\end{tikzpicture}
&

\begin{tikzpicture}[scale = .62]
\node[draw,circle,scale=0.75] at
(0,0) (0/1) {1};
\node[draw,circle,scale=0.75] at (0,-1) (0/2) {2};
\node[draw,circle,scale=0.75] at (0,-2) (0/3) {3};
\node[draw,circle,scale=0.75] at (0,-3) (0/4) {4};
\node[draw,circle,scale=0.75] at (0,-4) (0/5) {5};
\node[draw,circle,scale=0.75] at (0,-5) (0/6) {6};

\node[draw,circle,scale=0.75] at
(1,0) (1/1) {1};
\node[draw,circle,scale=0.75] at (1,-1) (1/2) {3};
\node[draw,circle,scale=0.75] at (1,-2) (1/3) {2};
\node[draw,circle,scale=0.75] at (1,-3) (1/4) {4};
\node[draw,circle,scale=0.75] at (1,-4) (1/5) {5};
\node[draw,circle,scale=0.75] at (1,-5) (1/6) {6};

\node[draw,circle,scale=0.75] at
(2,0) (2/1) {3};
\node[draw,circle,scale=0.75] at (2,-1) (2/2) {1};
\node[draw,circle,scale=0.75] at (2,-2) (2/3) {2};
\node[draw,circle,scale=0.75] at (2,-3) (2/4) {4};
\node[draw,circle,scale=0.75] at (2,-4) (2/5) {5};
\node[draw,circle,scale=0.75] at (2,-5) (2/6) {6};

\node[draw,circle,scale=0.75] at (3,0) (3/1) {3};
\node[draw,circle,scale=0.75] at (3,-1) (3/2) {1};
\node[draw,circle,scale=0.75] at (3,-2) (3/3) {4};
\node[draw,circle,scale=0.75] at (3,-3) (3/4) {2};
\node[draw,circle,scale=0.75] at (3,-4) (3/5) {5};
\node[draw,circle,scale=0.75] at (3,-5) (3/6) {6};

\node[draw,circle,scale=0.75] at (4,0) (4/1) {3};
\node[draw,circle,scale=0.75] at (4,-1) (4/2) {1};
\node[draw,circle,scale=0.75] at (4,-2) (4/3) {4};
\node[draw,circle,scale=0.75] at (4,-3) (4/4) {5};
\node[draw,circle,scale=0.75] at (4,-4) (4/5) {2};
\node[draw,circle,scale=0.75] at (4,-5) (4/6) {6};

\node[draw,circle,scale=0.75] at
(5,0) (5/1) {3};
\node[draw,circle,scale=0.75] at (5,-1) (5/2) {4};
\node[draw,circle,scale=0.75] at (5,-2) (5/3) {1};
\node[draw,circle,scale=0.75] at (5,-3) (5/4) {5};
\node[draw,circle,scale=0.75] at (5,-4) (5/5) {2};
\node[draw,circle,scale=0.75] at (5,-5) (5/6) {6};

\draw[ultra thick] (0/1) -- (1/1) -- (2/2) -- (3/2) -- (4/2) -- (5/3);
\draw[ultra thick] (0/2) -- (1/3) -- (2/3) -- (3/4) -- (4/5) -- (5/5);
\draw[ultra thick] (0/3) -- (1/2) -- (2/1) -- (3/1) -- (4/1) -- (5/1);
\draw[ultra thick] (0/4) -- (1/4) -- (2/4) -- (3/3) -- (4/3)--(5/2);
\draw[ultra thick] (0/5) -- (1/5) -- (2/5) -- (3/5) -- (4/4) -- (5/4);
\draw[ultra thick] (0/6) -- (1/6) -- (2/6) -- (3/6) -- (4/6) -- (5/6);
\end{tikzpicture}
\end{tabular}
\caption{Applying the Little map $\theta_5$ to $a=(3,1,4,5,2)$.}\label{fig-line-special}
\end{figure}
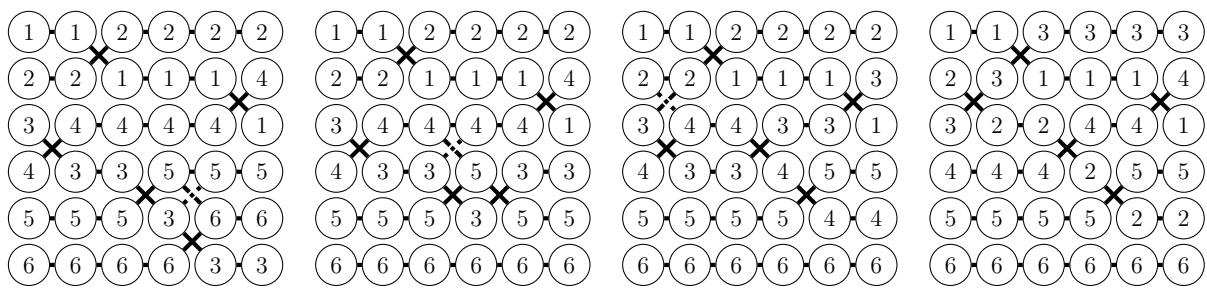

For example, let $u=241536$ and $k=5$. Then $w=u\, t_{5,6}=241563$
is the unique permutation in $\Psi(u, k)$. Let
 $a=(3, 1, 4, 5, 2)$ be a reduced word of $w$. Then $w_k=6,w_j=3$ and $t_1=4$.
 Applying the Little map
 $\theta_5$ to $a$, the resulting word is $(2,1,3,4,2)$,
 which is a reduced word of the permutation $w'=341526$ in $\Phi(u,k)$.
 The bumping process is illustrated in Figure \ref{fig-line-special}, where the dotted crossings indicate the bumped positions.

\begin{theo}[Little \mdseries{\cite{Little}}] \label{LNIU}
The map $\theta_k$   is a descent preserving
bijection.
\end{theo}

As described by Little \cite{Little}, the inverse of  $\theta_k$
can be stated as follows. Assume that  $w'$ is
a permutation in $\Phi(u,k)$ and that $i<k$ is the position in $I(u,k)$ such
that $w'=wt_{i,k}$. Let $a'=(a_1',a_2',\ldots,a_\ell')$ be a reduced word of $w'$ and denote $(a')^c=(n-a_1',\ldots,n-a_{\ell}')$. Note that $(a')^c$ is a reduced
word of the permutation $v=v_1v_2\cdots v_n$ where $v_i=n+1-w'_{n+1-i}$.
Then  one has
\begin{align}\label{invtheta}
\theta_k^{-1}(a')=(\theta_{n+1-k, n+1-w'_{i}}((a')^c))^c.
\end{align}

We conclude with two properties of the Little bump
due to Hamaker and Young  \cite{HY}, which will be  used in Section \ref{PP-8}.
First, the Little bump preserves the recording tableaux of  reduced words
when applying  the Edelman-Greene algorithm, which was conjectured by Lam \cite{Lam} and proved by Hamaker and Young \cite{HY}.

\begin{theo}[Hamaker-Young \mdseries{\cite[Theorem 2]{HY}}]\label{lam-HY}
Let  $a$ and $a'$ be two reduced words such
that there exists a sequence of Little
bumps changing $a$ to $a'$.
Then $Q(a)=Q(a')$.
\end{theo}

The second property   was essentially implied in the proof of Lemma 6 in \cite{HY}.

\begin{theo}[Hamaker-Young \mdseries{\cite[Lemma 6]{HY}}]\label{HY2}
Let $T$ be a reduced word tableau  and $a=\mathrm{column}(T)$. Assume that $b$ is   obtained from $a$ by
 applying a Little bump, and let $T'$ be the
 insertion tableau of the reverse $b^{\mathrm{rev}}$ of $b$ under
 the Edelman-Greene   algorithm. Then
$\mathrm{column}(T')=b$.
\end{theo}

\section{Bumpless pipedreams}\label{PP-2}

In this section, we give an overview of the structure
of bumpless pipedreams as well as  the droop operation
on bumpless pipedreams. The droop operation has a close connection
to the modified  LS-tree, as will be seen in
 Section \ref{PP-5}.

Lam, Lee and Shimozono \cite{LLS} introduced several versions
of bumpless pipedreams for a permutation $w\in S_n$. For the purpose of this paper, we shall only be   concerned  with
the bumpless pipedreams in the region
of a  square grid, which are  called $w$-square
bumpless pipedreams in \cite{LLS}.
Given an $n\times n$ square grid, we use the matrix coordinates for unit squares, that is, the row coordinates increase from top
to bottom and the column coordinates
increase from left to right.
Let us use $(i,j)$ to denote the box in row $i$ and column $j$.
A  bumpless
pipedream for $w$ consists of  $n$ pipes labeled $1,2,\ldots,n$, flowing from
the south boundary of the $n\times n$ square grid  to the east boundary,
such that
\begin{itemize}
\item[(1)] the pipe labeled $i$ enters from the south boundary in column $i$ and exits from the  east boundary  in row $w^{-1}(i)$;

\item [(2)] pipes can only go north or east;

  \item [(3)] no two pipes overlap any step or cross more than once;

  \item [(4)] no two pipes change their directions in a box simultaneously. In other words, each box looks like one of  the first six tiles
      as shown in Figure \ref{PPPs}:
      an empty box, an NW elbow, an SE elbow, a horizontal line, a vertical line,
      and a crossing.

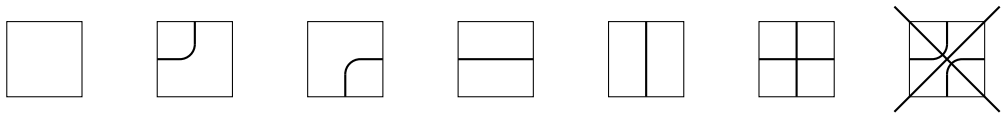
\begin{figure}[h]
\begin{center}
\begin{tikzpicture}

\def\squarepath{-- +(10mm,0mm) -- +(10mm,10mm) -- +(0mm,10mm) -- cycle}

\draw (0mm,0mm)\squarepath;

\draw (20mm,0mm)\squarepath;
\draw [-][thick](23mm,5mm) arc (270:360:0.2)
(20mm,5mm) -- (23mm,5mm);
\draw [-][thick](25mm,7mm) --(25mm,10mm);

\draw (40mm,0mm)\squarepath;
\draw [-][thick](47mm,5mm) arc (270:360:-0.2)
(47mm,5mm) -- (50mm,5mm);
\draw [-][thick](45mm,3mm) --(45mm,0mm);

\draw (60mm,0mm)\squarepath;
\draw [-][thick](60mm,5mm) --(70mm,5mm);

\draw (80mm,0mm)\squarepath;
\draw [-][thick](85mm,0mm) --(85mm,10mm);

\draw (100mm,0mm)\squarepath;
\draw [-][thick](100mm,5mm) --(110mm,5mm);
\draw [-][thick](105mm,0mm) --(105mm,10mm);

\draw (120mm,0mm)\squarepath;
\draw [-][thick](123mm,5mm) arc (270:360:0.2)
(120mm,5mm) -- (123mm,5mm);
\draw [-][thick](125mm,7mm) --(125mm,10mm);
\draw [-][thick](127mm,5mm) arc (270:360:-0.2)
(127mm,5mm) -- (130mm,5mm);
\draw[-][thick](125mm,3mm) --(125mm,0mm);
\draw[-][thick](118mm,-2mm) --(132mm,12mm);
\draw[-][thick](118mm,12mm) --(132mm,-2mm);
\end{tikzpicture}
\end{center}
\vspace{-.6cm}
\caption{Boxes   in a  bumpless pipedream.}\label{PPPs}
\end{figure}
 \end{itemize}
 \vspace{-.3cm}

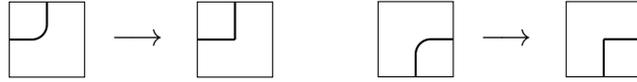
\begin{figure}[h]
\begin{center}
\begin{tabular}{cc}
\begin{tikzpicture}

\def\squarepath{-- +(10mm,0mm) -- +(10mm,10mm) -- +(0mm,10mm) -- cycle}

\draw (0mm,0mm)\squarepath;
\draw [-][thick](3mm,5mm) arc (270:360:0.2)
(0mm,5mm) -- (3mm,5mm);
\draw [-][thick](5mm,7mm) --(5mm,10mm);

\node at (17mm,5mm){$\longrightarrow$};

\draw (25mm,0mm)\squarepath;
\draw [-][thick]
(25mm,5mm) -- (30mm,5mm);
\draw [-][thick](30mm,5mm) --(30mm,10mm);
\end{tikzpicture}
\qquad
&

\qquad
\begin{tikzpicture}

\def\squarepath{-- +(10mm,0mm) -- +(10mm,10mm) -- +(0mm,10mm) -- cycle}

\draw (0mm,0mm)\squarepath;
\draw [-][thick](7mm,5mm) arc (270:360:-0.2)
(7mm,5mm) -- (10mm,5mm);
\draw [-][thick](5mm,3mm) --(5mm,0mm);

\node at (17mm,5mm){$\longrightarrow$};

\draw (25mm,0mm)\squarepath;
\draw [-][thick]
(30mm,5mm) -- (35mm,5mm);
\draw [-][thick](30mm,5mm) --(30mm,0mm);
\end{tikzpicture}
\end{tabular}
\end{center}
\vspace{-.5cm}
\caption{NW elbow and SE elbow.}\label{ON}
\end{figure}

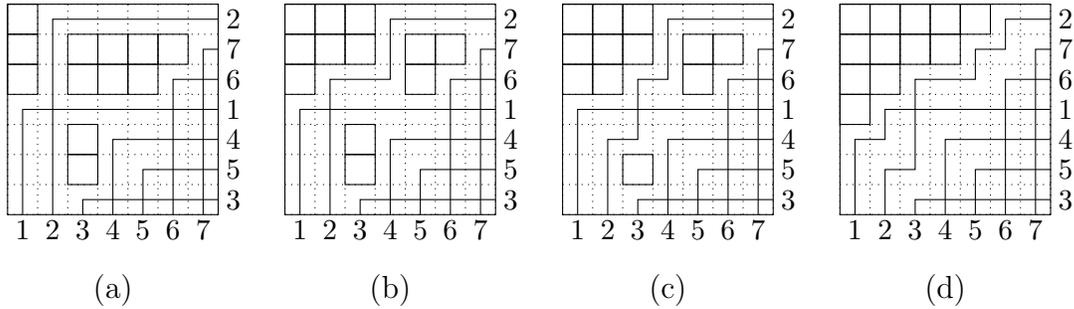
\begin{figure}[h]
\begin{center}
\begin{tabular}{cccc}


\begin{tikzpicture}[scale = 1]

\def\rectanglepath{-- +(4mm,0mm) -- +(4mm,4mm) -- +(0mm,4mm) -- cycle}

\def\squarepath{-- +(28mm,0mm) -- +(28mm,28mm) -- +(0mm,28mm) -- cycle}

\draw (0mm,0mm)\squarepath;
\draw [step=4mm,dotted] (0mm,0mm) grid (28mm,28mm);
\draw (0mm,24mm) \rectanglepath;
\draw (0mm,20mm)\rectanglepath;
\draw (0mm,16mm)\rectanglepath;
\draw (8mm,20mm)\rectanglepath;
\draw (8mm,16mm)\rectanglepath;
\draw (8mm,8mm)\rectanglepath;
\draw (8mm,4mm)\rectanglepath;
\draw (12mm,20mm)\rectanglepath;
\draw (12mm,16mm)\rectanglepath;
\draw (16mm,20mm)\rectanglepath;
\draw (16mm,16mm)\rectanglepath;
\draw (20mm,20mm)\rectanglepath;

\draw(2mm,0mm)--(2mm,14mm)--(28mm,14mm);
\draw(6mm,0mm)--(6mm,26mm)--(28mm,26mm);
\draw(10mm,0mm)--(10mm,2mm)--(28mm,2mm);
\draw(14mm,0mm)--(14mm,10mm)--(28mm,10mm);
\draw(18mm,0mm)--(18mm,6mm)--(28mm,6mm);
\draw(22mm,0mm)--(22mm,18mm)--(28mm,18mm);
\draw(26mm,0mm)--(26mm,22mm)--(28mm,22mm);

\node at (30mm,26mm) {\small{2}};
\node at (30mm,22mm) {\small{7}};
\node at (30mm,18mm) {\small{6}};
\node at (30mm,14mm) {\small{1}};
\node at (30mm,10mm) {\small{4}};
\node at (30mm,6mm) {\small{5}};
\node at (30mm,2mm) {\small{3}};

\node at (2mm,-2mm) {\small{1}};
\node at (6mm,-2mm) {\small{2}};
\node at (10mm,-2mm) {\small{3}};
\node at (14mm,-2mm)  {\small{4}};
\node at (18mm,-2mm) {\small{5}};
\node at (22mm,-2mm) {\small{6}};
\node at (26mm,-2mm) {\small{7}};

\node at (14mm,-10mm) {(a)};
\end{tikzpicture}
&

\begin{tikzpicture}[scale =1]

\def\rectanglepath{-- +(4mm,0mm) -- +(4mm,4mm) -- +(0mm,4mm) -- cycle}

\def\squarepath{-- +(28mm,0mm) -- +(28mm,28mm) -- +(0mm,28mm) -- cycle}

\draw (0mm,0mm)\squarepath;
\draw [step=4mm,dotted] (0mm,0mm) grid (28mm,28mm);
\draw (0mm,24mm)\rectanglepath;
\draw (0mm,20mm)\rectanglepath;
\draw (0mm,16mm)\rectanglepath;
\draw (4mm,24mm)\rectanglepath;
\draw (4mm,20mm)\rectanglepath;
\draw (8mm,24mm)\rectanglepath;
\draw (8mm,20mm)\rectanglepath;
\draw (8mm,8mm)\rectanglepath;
\draw (8mm,4mm)\rectanglepath;
\draw(16mm,20mm)\rectanglepath;
\draw(16mm,16mm)\rectanglepath;
\draw(20mm,20mm)\rectanglepath;

\draw(2mm,0mm)--(2mm,14mm)--(28mm,14mm);
\draw(6mm,0mm)--(6mm,18mm)--(14mm,18mm)--(14mm,26mm)--(28mm,26mm);
\draw(10mm,0mm)--(10mm,2mm)--(28mm,2mm);
\draw(14mm,0mm)--(14mm,10mm)--(28mm,10mm);
\draw(18mm,0mm)--(18mm,6mm)--(28mm,6mm);
\draw(22mm,0mm)--(22mm,18mm)--(28mm,18mm);
\draw(26mm,0mm)--(26mm,22mm)--(28mm,22mm);

\node at (30mm,26mm) {\small{2}};
\node at (30mm,22mm) {\small{7}};
\node at (30mm,18mm) {\small{6}};
\node at (30mm,14mm) {\small{1}};
\node at (30mm,10mm) {\small{4}};
\node at (30mm,6mm) {\small{5}};
\node at (30mm,2mm) {\small{3}};

\node at (2mm,-2mm) {\small{1}};
\node at (6mm,-2mm) {\small{2}};
\node at (10mm,-2mm) {\small{3}};
\node at (14mm,-2mm)  {\small{4}};
\node at (18mm,-2mm) {\small{5}};
\node at (22mm,-2mm) {\small{6}};
\node at (26mm,-2mm) {\small{7}};

\node at (14mm,-10mm) {(b)};
\end{tikzpicture}
&

\begin{tikzpicture}[scale = 1]

\def\rectanglepath{-- +(4mm,0mm) -- +(4mm,4mm) -- +(0mm,4mm) -- cycle}

\def\squarepath{-- +(28mm,0mm) -- +(28mm,28mm) -- +(0mm,28mm) -- cycle}

\draw (0mm,0mm)\squarepath;
\draw [step=4mm,dotted] (0mm,0mm) grid (28mm,28mm);
\draw (0mm,24mm)\rectanglepath;
\draw (0mm,20mm)\rectanglepath;
\draw (0mm,16mm)\rectanglepath;
\draw (4mm,24mm)\rectanglepath;
\draw (4mm,20mm)\rectanglepath;
\draw (4mm,16mm)\rectanglepath;
\draw (8mm,24mm)\rectanglepath;
\draw (8mm,20mm)\rectanglepath;
\draw (8mm,4mm)\rectanglepath;
\draw(16mm,20mm)\rectanglepath;
\draw(16mm,16mm)\rectanglepath;
\draw(20mm,20mm)\rectanglepath;

\draw(2mm,0mm)--(2mm,14mm)--(28mm,14mm);
\draw(6mm,0mm)--(6mm,10mm)--(10mm,10mm)--(10mm,18mm)--(14mm,18mm)--(14mm,26mm)--(28mm,26mm);
\draw(10mm,0mm)--(10mm,2mm)--(28mm,2mm);
\draw(14mm,0mm)--(14mm,10mm)--(28mm,10mm);
\draw(18mm,0mm)--(18mm,6mm)--(28mm,6mm);
\draw(22mm,0mm)--(22mm,18mm)--(28mm,18mm);
\draw(26mm,0mm)--(26mm,22mm)--(28mm,22mm);

\node at (30mm,26mm) {\small{2}};
\node at (30mm,22mm) {\small{7}};
\node at (30mm,18mm) {\small{6}};
\node at (30mm,14mm) {\small{1}};
\node at (30mm,10mm) {\small{4}};
\node at (30mm,6mm) {\small{5}};
\node at (30mm,2mm) {\small{3}};

\node at (2mm,-2mm) {\small{1}};
\node at (6mm,-2mm) {\small{2}};
\node at (10mm,-2mm) {\small{3}};
\node at (14mm,-2mm)  {\small{4}};
\node at (18mm,-2mm) {\small{5}};
\node at (22mm,-2mm) {\small{6}};
\node at (26mm,-2mm) {\small{7}};

\node at (14mm,-10mm) {(c)};
\end{tikzpicture}

&

\begin{tikzpicture}[scale = 1]

\def\rectanglepath{-- +(4mm,0mm) -- +(4mm,4mm) -- +(0mm,4mm) -- cycle}

\def\squarepath{-- +(28mm,0mm) -- +(28mm,28mm) -- +(0mm,28mm) -- cycle}

\draw (0mm,0mm)\squarepath;
\draw [step=4mm,dotted] (0mm,0mm) grid (28mm,28mm);
\draw (0mm,24mm)\rectanglepath;
\draw (0mm,20mm)\rectanglepath;
\draw (0mm,16mm)\rectanglepath;
\draw (0mm,12mm)\rectanglepath;
\draw (4mm,24mm)\rectanglepath;
\draw (4mm,20mm)\rectanglepath;
\draw (4mm,16mm)\rectanglepath;
\draw (8mm,24mm)\rectanglepath;
\draw (8mm,20mm)\rectanglepath;
\draw(12mm,24mm)\rectanglepath;
\draw(12mm,20mm)\rectanglepath;
\draw(16mm,24mm)\rectanglepath;

\draw(2mm,0mm)--(2mm,10mm)--(6mm,10mm)--(6mm,14mm)--(28mm,14mm);
\draw(6mm,0mm)--(6mm,6mm)--(10mm,6mm)--(10mm,18mm)--(18mm,18mm)--(18mm,22mm)--(22mm,22mm)--(22mm,26mm)--(28mm,26mm);
\draw(10mm,0mm)--(10mm,2mm)--(28mm,2mm);
\draw(14mm,0mm)--(14mm,10mm)--(28mm,10mm);
\draw(18mm,0mm)--(18mm,6mm)--(28mm,6mm);
\draw(22mm,0mm)--(22mm,18mm)--(28mm,18mm);
\draw(26mm,0mm)--(26mm,22mm)--(28mm,22mm);

\node at (30mm,26mm) {\small{2}};
\node at (30mm,22mm) {\small{7}};
\node at (30mm,18mm) {\small{6}};
\node at (30mm,14mm) {\small{1}};
\node at (30mm,10mm) {\small{4}};
\node at (30mm,6mm) {\small{5}};
\node at (30mm,2mm) {\small{3}};

\node at (2mm,-2mm) {\small{1}};
\node at (6mm,-2mm) {\small{2}};
\node at (10mm,-2mm) {\small{3}};
\node at (14mm,-2mm)  {\small{4}};
\node at (18mm,-2mm) {\small{5}};
\node at (22mm,-2mm) {\small{6}};
\node at (26mm,-2mm) {\small{7}};

\node at (14mm,-10mm) {(d)};
\end{tikzpicture}
\end{tabular}
\vspace{-.3cm}
\caption{Four  Bumpless pipedreams of $w=2761453$.}
\label{bumpless}
\end{center}
\end{figure}

 For ease of drawing pictures, an NW elbow and an SE elbow
are replaced by the tiles  given in Figure \ref{ON}. Since two pipes cannot bump at the same box, there is no ambiguity for this simplification. We illustrate
four  bumpless pipedreams of $w=2761453$ in Figure \ref{bumpless}.

Lam, Lee and Shimozono \cite{LLS} showed that the bumpless pipedreams of $w$ can generate the Schubert polynomial $\mathfrak{S}_{w}(x)$, or more generally  the double Schubert polynomial $\S_{w}(x;y)$.
For a   bumpless pipedream $P$, define the weight $\mathrm{wt}(P)$ of $P$ to be the product of $x_i-y_j$ over all
empty boxes of $P$ in row $i$ and column $j$.

\begin{theo}[Lam-Lee-Shimozono \mdseries{\cite[Theorem 5.13]{LLS}}]
For any permutation $w$,
\begin{equation}\label{SS}
\S_w(x;y)=\sum_{P}\mathrm{wt}(P),
\end{equation}
where the sum is over the bumpless pipedreams of $w$.
\end{theo}

Lam, Lee and Shimozono \cite{LLS} also discovered
that a specific family of   bumpless pipedreams, called EG-pipedreams,
can be used to interpret  the Edelman-Greene coefficients.
A bumpless pipedream $P$  is called an EG-pipedream if all
 the empty boxes of $P$ are at the northwest corner, where they form a Young diagram $\lambda=\lambda(P)$, called the shape of $P$.
For example, Figure \ref{bumpless}(d) is an EG-pipedream with shape $\lambda=(5,4,2,1)$.

\begin{theo}[Lam-Lee-Shimozono \mdseries{\cite[Theorem 5.14]{LLS}}]\label{dsb}
The Edelman-Greene coefficient $c^w_\lambda$
is equal to the number of EG-pipedreams of $w$ with shape $\lambda$.
\end{theo}

\begin{figure}[h]
\begin{center}
\begin{tabular}{ccc}
\begin{tikzpicture}

\def\rectanglepath{-- +(4mm,0mm) -- +(4mm,4mm) -- +(0mm,4mm) -- cycle}

\def\squarepath{-- +(24mm,0mm) -- +(24mm,24mm) -- +(0mm,24mm) -- cycle}

\draw (0mm,0mm)\squarepath;
\draw [step=4mm,dotted] (0mm,0mm) grid (24mm,24mm);
\draw (0mm,20mm)\rectanglepath;
\draw (0mm,16mm)\rectanglepath;
\draw (4mm,20mm)\rectanglepath;
\draw (4mm,16mm)\rectanglepath;
\draw (8mm,20mm)\rectanglepath;
\draw (12mm,20mm)\rectanglepath;

\draw(2mm,0mm)--(2mm,14mm)--(24mm,14mm);
\draw(6mm,0mm)--(6mm,10mm)--(10mm,10mm)--(10mm,18mm)--(24mm,18mm);
\draw(10mm,0mm)--(10mm,6mm)--(14mm,6mm)--(14mm,10mm)--(18mm,10mm)--(18mm,22mm)--(24mm,22mm);
\draw(14mm,0mm)--(14mm,2mm)--(24mm,2mm);
\draw(18mm,0mm)--(18mm,6mm)--(24mm,6mm);
\draw(22mm,0mm)--(22mm,10mm)--(24mm,10mm);
\node at (26mm,22mm) {\small{3}};
\node at (26mm,18mm) {\small{2}};
\node at (26mm,14mm) {\small{1}};
\node at (26mm,10mm) {\small{6}};
\node at (26mm,6mm) {\small{5}};
\node at (26mm,2mm) {\small{4}};
%
\node at (2mm,-2mm) {\small{1}};
\node at (6mm,-2mm) {\small{2}};
\node at (10mm,-2mm) {\small{3}};
\node at (14mm,-2mm)  {\small{4}};
\node at (18mm,-2mm) {\small{5}};
\node at (22mm,-2mm) {\small{6}};
\end{tikzpicture}&

\quad
\begin{tikzpicture}

\def\rectanglepath{-- +(4mm,0mm) -- +(4mm,4mm) -- +(0mm,4mm) -- cycle}

\def\squarepath{-- +(24mm,0mm) -- +(24mm,24mm) -- +(0mm,24mm) -- cycle}

\draw (0mm,0mm)\squarepath;
\draw [step=4mm,dotted] (0mm,0mm) grid (24mm,24mm);
\draw (0mm,20mm)\rectanglepath;
\draw (0mm,16mm)\rectanglepath;
\draw (0mm,12mm)\rectanglepath;
\draw (4mm,20mm)\rectanglepath;
\draw (4mm,16mm)\rectanglepath;
\draw (8mm,20mm)\rectanglepath;

\draw(2mm,0mm)--(2mm,10mm)--(18mm,10mm)--(18mm,14mm)--(24mm,14mm);
\draw(6mm,0mm)--(6mm,14mm)--(10mm,14mm)--(10mm,18mm)--(24mm,18mm);
\draw(10mm,0mm)--(10mm,6mm)--(14mm,6mm)--(14mm,22mm)--(24mm,22mm);
\draw(14mm,0mm)--(14mm,2mm)--(24mm,2mm);
\draw(18mm,0mm)--(18mm,6mm)--(24mm,6mm);
\draw(22mm,0mm)--(22mm,10mm)--(24mm,10mm);
\node at (26mm,22mm) {\small{3}};
\node at (26mm,18mm) {\small{2}};
\node at (26mm,14mm) {\small{1}};
\node at (26mm,10mm) {\small{6}};
\node at (26mm,6mm) {\small{5}};
\node at (26mm,2mm) {\small{4}};
%
\node at (2mm,-2mm) {\small{1}};
\node at (6mm,-2mm) {\small{2}};
\node at (10mm,-2mm) {\small{3}};
\node at (14mm,-2mm)  {\small{4}};
\node at (18mm,-2mm) {\small{5}};
\node at (22mm,-2mm) {\small{6}};
\end{tikzpicture}&

\quad
\begin{tikzpicture}

\def\rectanglepath{-- +(4mm,0mm) -- +(4mm,4mm) -- +(0mm,4mm) -- cycle}

\def\squarepath{-- +(24mm,0mm) -- +(24mm,24mm) -- +(0mm,24mm) -- cycle}

\draw (0mm,0mm)\squarepath;
\draw [step=4mm,dotted] (0mm,0mm) grid (24mm,24mm);
\draw (0mm,20mm)\rectanglepath;
\draw (0mm,16mm)\rectanglepath;
\draw (0mm,12mm)\rectanglepath;
\draw (4mm,20mm)\rectanglepath;
\draw (4mm,16mm)\rectanglepath;
\draw (8mm,20mm)\rectanglepath;

\draw(2mm,0mm)--(2mm,10mm)--(6mm,10mm)--(6mm,14mm)--(24mm,14mm);
\draw(6mm,0mm)--(6mm,6mm)--(14mm,6mm)--(14mm,10mm)--(18mm,10mm)--(18mm,18mm)--(24mm,18mm);
\draw(10mm,0mm)--(10mm,18mm)--(14mm,18mm)--(14mm,22mm)--(24mm,22mm);
\draw(14mm,0mm)--(14mm,2mm)--(24mm,2mm);
\draw(18mm,0mm)--(18mm,6mm)--(24mm,6mm);
\draw(22mm,0mm)--(22mm,10mm)--(24mm,10mm);
\node at (26mm,22mm) {\small{3}};
\node at (26mm,18mm) {\small{2}};
\node at (26mm,14mm) {\small{1}};
\node at (26mm,10mm) {\small{6}};
\node at (26mm,6mm) {\small{5}};
\node at (26mm,2mm) {\small{4}};
%
\node at (2mm,-2mm) {\small{1}};
\node at (6mm,-2mm) {\small{2}};
\node at (10mm,-2mm) {\small{3}};
\node at (14mm,-2mm)  {\small{4}};
\node at (18mm,-2mm) {\small{5}};
\node at (22mm,-2mm) {\small{6}};
\end{tikzpicture}

\end{tabular}

\caption{The EG-pipedreams of $w=321654$.}\label{EGs}
\end{center}
\end{figure}
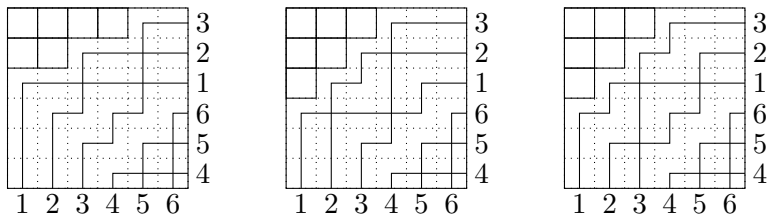

For example, there are three  EG-pipedreams for  $w=321654$,
as illustrated in Figure \ref{EGs}.
By Theorem \ref{dsb}, we see that
$F_w(x)=s_{(4,2)}(x)+2\,s_{(3,2,1)}(x).$

In the remaining of this section, we recall the droop
operation on bumpless pipedreams introduced by
Lam, Lee and Shimozono \cite{LLS}.
As will be seen in Corollary \ref{droop-app}, the droop operation
is closely related  to the modified LS-tree.

Let $P$ be a bumpless pipedream of $w$. The pipes in $P$ are determined by the locations of the NW elbows and the SE elbows.
Generally speaking, a  droop is a local move which swaps an SE elbow $e$ with
 an empty  box $t$, when the SE elbow $e$ lies strictly to the
northwest of the empty box $t$.
To be more specific, let $R$ be the
rectangle with northwest corner $e$ and southeast corner $t$,
 and let $L$ be the pipe passing through the SE elbow
$e$. A droop is allowed only if
\begin{itemize}
  \item [(1)] the pipe $L$ passes through the westmost column and northmost row of $R$;

  \item [(2)]  the rectangle $R$ contains only one elbow: the SE elbow which is at $e$;

  \item [(3)] after the droop we obtain another bumpless pipedream.
\end{itemize}

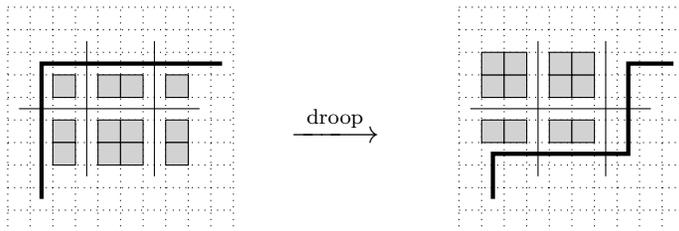
\begin{figure}[h]
\begin{center}
\begin{tikzpicture}

\def\rectanglepath{-- +(3mm,0mm) -- +(3mm,3mm) -- +(0mm,3mm) -- cycle}

\def\squarepath{-- +(27mm,0mm) -- +(27mm,27mm) -- +(0mm,27mm) -- cycle}

\def\srectanglepath{[fill=gray!35!white,draw=black!20!black]-- +(3mm,0mm) -- +(3mm,3mm) -- +(0mm,3mm) -- cycle}

\draw [step=3mm,dotted] (0mm,0mm)grid (30mm,30mm);

\draw (6mm,18mm) \srectanglepath;
\draw (6mm,12mm) \srectanglepath;
\draw (6mm,9mm) \srectanglepath;
\draw (12mm,18mm) \srectanglepath;
\draw (12mm,12mm) \srectanglepath;
\draw (12mm,9mm) \srectanglepath;
\draw (15mm,18mm) \srectanglepath;
\draw (15mm,12mm) \srectanglepath;
\draw (15mm,9mm) \srectanglepath;
\draw (21mm,18mm) \srectanglepath;
\draw (21mm,12mm) \srectanglepath;
\draw (21mm,9mm) \srectanglepath;

\draw[-][ultra thick](4.5mm,4.5mm)--(4.5mm,22.5mm)--(28.5mm,22.5mm);
\draw(10.5mm,7.5mm)--(10.5mm,25.5mm);
\draw(19.5mm,7.5mm)--(19.5mm,25.5mm);
\draw(1.5mm,16.5mm)--(25.5mm,16.5mm);
\draw [step=3mm,dotted] (60mm,0mm)grid (90mm,30mm);
\draw [dotted] (60mm,0mm)--(60mm,30mm);

\draw (63mm,18mm) \srectanglepath;
\draw (63mm,12mm) \srectanglepath;
\draw (63mm,21mm) \srectanglepath;
\draw (66mm,18mm) \srectanglepath;
\draw (66mm,12mm) \srectanglepath;
\draw (66mm,21mm) \srectanglepath;
\draw (72mm,18mm) \srectanglepath;
\draw (72mm,12mm) \srectanglepath;
\draw (72mm,21mm) \srectanglepath;
\draw (75mm,18mm) \srectanglepath;
\draw (75mm,12mm) \srectanglepath;
\draw (75mm,21mm) \srectanglepath;

\draw[-][ultra thick](64.5mm,4.5mm)--(64.5mm,10.5mm)--(82.5mm,10.5mm)
--(82.5mm,22.5mm)--(88.5mm,22.5mm);
\draw(70.5mm,7.5mm)--(70.5mm,25.5mm);
\draw(79.5mm,7.5mm)--(79.5mm,25.5mm);
\draw(61.5mm,16.5mm)--(85.5mm,16.5mm);
\node at (43.5mm,13mm) {$\xlongrightarrow[\,\qquad]{\text{droop}}$};

\end{tikzpicture}
\caption{The droop operation.}\label{DROOP}
\end{center}
\end{figure}

After a droop, the pipe $L$ travels along the southmost row and
eastmost column of $R$, and  an NW elbow   occupies the box that used
to be empty while the box that contained an SE elbow becomes an empty box.
Pipes in $P$ except the pipe $L$ do not change after the droop.
Figure \ref{DROOP} is an illustration of the local move of the droop operation.
For more  examples,  Figure \ref{bumpless}(b) is obtained from
Figure \ref{bumpless}(a) by a  droop, and  Figure \ref{bumpless}(c)
is  obtained from Figure \ref{bumpless}(b) by a droop.

Lam, Lee and Shimozono \cite{LLS} proved that each  bumpless
pipedream of $w$ can be generated from the Rothe pipedream of $w$ by applying a sequence of
droops.
The  Rothe  pipedream   of $w$, denoted  $D(w)$, is the unique bumpless pipedream of $w$ that does not  contain any NW elbows,
 see Figure \ref{bumpless}(a) for an example.
Clearly, for $1\leq i\leq n$,    pipe $i$ passes through
exactly the SE elbow $(i, w_i)$ in row $i$ and column $w_i$.

The  Rothe  pipedream $D(w)$ can also be constructed
 as follows. From the center of each box $(i, w_i)$,
draw a horizonal line to the right and a vertical line to the bottom.
This forms $n$ hooks with turning points at the center
of the boxes $(i, w_i)$.
Thinking of each hook as a pipe,  the pipes together with the
$n\times n$ grid form the   Rothe  pipedream $D(w)$.
The empty boxes of $D(w)$ are known as the Rothe diagram of $w$, denoted $\mathrm{Rothe}(w)$,
which encode  the positions of  inversions of $w$. That is, there is an empty box of $\mathrm{Rothe}(w)$ at $(i,j)$  if and only if $w_i>j$ and the number $j$ appears in $w$ after the position $i$.

\begin{prop}[Lam-Lee-Shimozono \mdseries{\cite[Proposition 5.3]{LLS}}]
For a permutation $w$, every  bumpless pipedream of $w$
can be obtained from the Rothe pipedream $D(w)$ by a sequence of droops.
\end{prop}

\section{Modified Lascoux-Sch\"utzenberger  tree}\label{PP-5}

In this section, we introduce  the structure of the
 modified   LS-tree of a permutation.
We also explain the relation between the  modified  LS-trees
and the droop operations.

For a permutation  $w=w_1w_2\cdots w_n$,
let $p$ be  the largest index in $w$ such that
there are indices $i$ and $j$ satisfying
\begin{equation}\label{Condi}
 i<p<j\ \ \ \text{and}\ \ \  w_i<w_j<w_p,
\end{equation}
namely,
\begin{align}\label{PPPP}
p=\max\{t\,|\,1\leq t\leq n-1,\ \exists\ i<t<j,\ \text{s.t.}\  w_i<w_j<w_t\}.
\end{align}
In other words, $p$ is the largest position such that there is a subsequence $w_iw_pw_j$ which is order-isomorphic to  132.
For example, for the permutation $w=2431$ we have $p=2$,
whereas for the permutation $w=3421$ we cannot find a position satisfying
the condition in \eqref{Condi}.

By definition, there does not exist an index $p$
satisfying \eqref{Condi} if and only if $w$ is 132-avoiding.
If a permutation $w$ contains a 132 pattern, we also call
  $w$ a non-dominant permutation.
%

%
%

Assume that $w$ is a non-dominant permutation. Let
$q$ be the largest index  after $p$ such that $w_p>w_q$ and  there exists an index  $i<p$ such that $w_i<w_q$, namely,
\begin{align}\label{QQQQ}
q=\max \{j\,|\,j>p,\,w_j<w_{p},\,\exists\,i<p,\, \text{s.t.}\,w_i<w_j\}.
\end{align}
For example, for   $w=645978321$,
we have $p=4$ and $q=6$.
We have the following equivalent description of the index $q$.

\begin{lemma}\label{XXX}
Let $w$ be a non-dominant permutation. Then the index $q$ defined in \eqref{QQQQ} is the index of the largest number following  $w_p$ which is less than $w_p$.
\end{lemma}

\pf
Suppose otherwise  that $q'\neq q$ is   the index of the largest number
following  $w_p$ which is less than $w_p$. By the definition of $q$ in \eqref{QQQQ},
we see that $p<q'<q$.
On the other hand, by the definition of $p$, there is an index $i<p$ such that $w_i<w_q<w_p$.
Thus we have $w_i<w_{q}<w_{q'}$, implying that $q'$ satisfies
\eqref{Condi}. This  contradicts the choice of
the index $p$. So the  proof is complete.
\qed

The following proposition implies that the indices $p$ and $q$
 play  similar roles to the indices $r$ and $s$ as defined  in
Section \ref{lstree}.

\begin{prop}\label{CPU}
Let $w$ be  a non-dominant permutation.
Then we have
\begin{align}
\Psi(wt_{p,q},p)=\{w\}.
\end{align}
\end{prop}

\pf Write $u=wt_{p,q}=u_1u_2\cdots u_n$, where $u_p=w_q$,
$u_q=w_p$, and  $u_i=w_i$ for $i\neq p, q$.
By Lemma \ref{XXX}, we see that
$\ell(wt_{p,q})=\ell(w)-1,$
which implies
$\ell(u t_{p,q})=\ell(w)=\ell(u)+1.$
So we have $w\in \Psi(u, p)$.

Suppose that there is another  permutation $w'\in\Psi(u, p)$ which is
not equal to $w$. By the definition of $\Psi(u, p)$ in
\eqref{bigpsi}, there exists an index $q'\neq q$ such that
$q'>p$, $w'=u \,t_{p, q'}$ and
$\ell(w')=\ell(u)+1.$
Then we have $u_{q'}>u_p$, that is,
\begin{equation}\label{AA-2}
w_{q'}>w_q=u_p.
\end{equation}
There are two cases.

Case 1: $p<q'<q$.  By the definition of $q$,  there is an index $i<p$ such that $w_i<w_q<w_p$. In view of \eqref{AA-2}, we see that $w_i<w_{q}<w_{q'}$,
which implies that the index $q'$ satisfies the condition in \eqref{Condi}.
This contradicts the choice of
the index $p$.

Case 2: $q'>q$. In this case, since $\ell(w')=\ell(u)+1$, we must have $w_{q'}<w_{p}$, which,
together with \eqref{AA-2}, implies that $w_{q}<w_{q'}<w_{p}$.
This is   contrary to Lemma \ref{XXX}. So the proof is complete.
\qed

Combining \eqref{StanRe2} and  Proposition \ref{CPU},
we arrive at the following transition relation
satisfied by Stanley symmetric functions.

\begin{theo}\label{CC-1}
Let $w\in S_n$ be a non-dominant permutation. Then we have
\begin{align}\label{BB-1}
F_w(x)=\sum_{w'\in \Phi(wt_{p,q}, p)}F_{w'}(x).
\end{align}
Moreover, by the choices of the indices   $p$ and $q$, the set $I(wt_{p,q},p)$ is
not empty, and
hence each permutation in $ \Phi(wt_{p,q},p)$ is still a permutation on
$\{1,2,\ldots,n\}$.
\end{theo}

For example, for  $w=645978321$, we see that $p=4$,
$q=6$ and so we have
$wt_{p,q}=645879321.$
Moreover, one can check that $ {I}(wt_{p,q}, p)=\{1,3\}$. Therefore,
\[\Phi(wt_{p,q},p)=\{845679321,648579321\}.\]

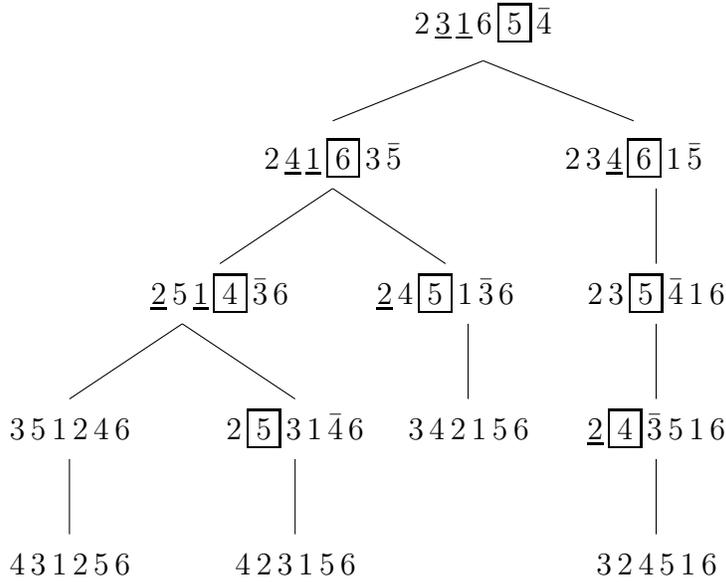
\begin{figure}[!htb]
\begin{center}
\begin{tikzpicture}

\draw (50mm,72mm) node{$2\,\underline{3}\,\underline{1}\,6\,\boxed{5}\,\bar{4}$};

\draw (30mm,54mm) node{$2\,\underline{4}\,\underline{1}\,\boxed{6}\,3\,\bar{5}$};

\draw (70mm,54mm) node{$2\,3\,\underline{4}\,\boxed{6}\,1\,\bar{5}$};


\draw (15mm,36mm)
node{$\underline{2}\,5\,\underline{1}\,\boxed{4}\,\bar{3}\,6$};

\draw (45mm,36mm)
node{$\underline{2}\,4\,\boxed{5}\,1\,\bar{3}\,6$};

\draw (73mm,36mm) node{$2\,3\,\boxed{5}\,\bar{4}\,1\,6$};

%
\draw (-5mm,18mm) node{$3\,5\,1\,2\,4\,6$};

\draw (25mm,18mm) node{$2\,\boxed{5}\,3\,1\,\bar{4}\,6$};

\draw (48mm,18mm)
node{$3\,4\,2\,1\,5\,6$};

\draw (73mm,18mm) node{$\underline{2}\,\boxed{4}\,\bar{3}\,5\,1\,6$};

\draw (-5mm,0mm) node{$4\,3\,1\,2\,5\,6$};
\draw (25mm,0mm) node{$4\,2\,3\,1\,5\,6$};

\draw (73mm,0mm) node{$3\,2\,4\,5\,1\,6$};

\draw(50mm,67mm)--(70mm,59mm);
\draw(50mm,67mm)--(30mm,59mm);
\draw(30mm,50mm)--(15mm,40mm);
\draw(30mm,50mm)--(45mm,40mm);
\draw(73mm,50mm)--(73mm,40mm);
\draw(10mm,32mm)--(-5mm,22mm);
\draw(10mm,32mm)--(25mm,22mm);
\draw(48mm,32mm)--(48mm,22mm);
\draw(73mm,32mm)--(73mm,22mm);

\draw(-5mm,14mm)--(-5mm,4mm);
\draw(25mm,14mm)--(25mm,4mm);

\draw(73mm,14mm)--(73mm,4mm);

\end{tikzpicture}
\caption{The modified LS-tree of $w=231654$.}\label{MLS-tree}
\end{center}
\end{figure}

The construction of the modified LS-tree of $w$ relies
on Theorem \ref{CC-1}.  To be  specific, for a permutation
 $w\in S_n$, iterate
the relation in \eqref{BB-1} until each leaf is
a dominant permutation on $\{1,2,\ldots,n\}$.
The resulting  tree is called
the  modified LS-tree of $w$.
Figure \ref{MLS-tree}
illustrates the   modified LS-tree of $w=231654$.

\begin{re}

The   construction of a modified LS-tree is  feasible because
the process
eventually terminates. This can be seen as follows.
For two sequences
$c=(c_1,c_2,\ldots,c_n)$ and $c'=(c'_1,c'_2,\ldots,c'_n)$ of nonnegative integers of length $n$,
write  $c\geq c'$ if for   $1\le i\le n$,
\[c_1+\cdots+ c_i\geq c_1'+\cdots+ c_i'.\]
This defines a partial   order    on the sequences of
nonnegative integers of length $n$.
Assume that  $w\in S_n$ is a non-dominant permutation. It is easy to check that
for  $w'\in\Phi(wt_{p,q},p)$, the Lehmer code of $w$ is smaller than the Lehmer code of $w'$  under the above order. Hence,
the nodes in any path from the root to a leaf in
 the modified LS-tree
are labeled by distinct permutations
on $\{1,2,\ldots,n\}$. So the process
  terminates.
\end{re}

The  differences between
a modified LS-tree and an ordinary LS-tree are obvious.
First, each node in a modified LS-tree is labeled with
 a permutaton on $\{1,2,\ldots, n\}$, whereas in an ordinary  LS-tree,
 there may exist nodes which are labeled with  permutations
 on $\{1,2,\ldots, m\}$ with $m>n$.
 Second, each leaf in a modified LS-tree is labeled with
  a dominant permutation, whereas each leaf
  in an ordinary LS-tree is labeled with
  a Grassmannian permutation.
  To make a comparison, see Figure \ref{LS-tree} and
Figure \ref{MLS-tree}.

Since each leaf in the modified LS-tree of $w$ is a dominant permutation, according to Theorem \ref{CC-1}, we
obtain the following interpretation for the
Edelman-Greene coefficients.

\begin{theo}\label{EG-3}
The Edelman-Greene coefficient $c^w_\lambda$ equals the number of leaves in the modified LS-tree of $w$ whose labels are the dominant permutations with Lehmer code $\lambda$.
\end{theo}

\begin{figure}[h]
\begin{center}
\begin{tikzpicture}

\def\rectanglepath{-- +(2mm,0mm) -- +(2mm,2mm) -- +(0mm,2mm) -- cycle}

\def\squarepath{-- +(12mm,0mm) -- +(12mm,12mm) -- +(0mm,12mm) -- cycle}

\draw (44mm,72mm) \squarepath;

\draw (44mm,82mm) \rectanglepath;
\draw (44mm,80mm) \rectanglepath;
\draw (50mm,76mm) \rectanglepath;
\draw (52mm,76mm) \rectanglepath;
\draw (50mm,74mm) \rectanglepath;

\draw(45mm,72mm)--(45mm,79mm)--(56mm,79mm);
\draw(47mm,72mm)--(47mm,83mm)--(56mm,83mm);
\draw(49mm,72mm)--(49mm,81mm)--(56mm,81mm);
\draw(51mm,72mm)--(51mm,73mm)--(56mm,73mm);
\draw(53mm,72mm)--(53mm,75mm)--(56mm,75mm);
\draw(55mm,72mm)--(55mm,77mm)--(56mm,77mm);

%
\draw (32mm,54mm) \squarepath;
\draw (32mm,64mm) \rectanglepath;
\draw (32mm,62mm) \rectanglepath;
\draw (36mm,62mm) \rectanglepath;
\draw (40mm,58mm) \rectanglepath;
\draw (36mm,58mm) \rectanglepath;

\draw(33mm,54mm)--(33mm,61mm)--(44mm,61mm);
\draw(35mm,54mm)--(35mm,65mm)--(44mm,65mm);
\draw(37mm,54mm)--(37mm,57mm)--(44mm,57mm);
\draw(39mm,54mm)--(39mm,63mm)--(44mm,63mm);
\draw(41mm,54mm)--(41mm,55mm)--(44mm,55mm);
\draw(43mm,54mm)--(43mm,59mm)--(44mm,59mm);

\draw (56mm,54mm) \squarepath;
\draw (56mm,58mm) \rectanglepath;
\draw (56mm,60mm) \rectanglepath;
\draw (56mm,62mm) \rectanglepath;
\draw (56mm,64mm) \rectanglepath;
\draw (64mm,58mm) \rectanglepath;

\draw(57mm,54mm)--(57mm,57mm)--(68mm,57mm);
\draw(59mm,54mm)--(59mm,65mm)--(68mm,65mm);
\draw(61mm,54mm)--(61mm,63mm)--(68mm,63mm);
\draw(63mm,54mm)--(63mm,61mm)--(68mm,61mm);
\draw(65mm,54mm)--(65mm,55mm)--(68mm,55mm);
\draw(67mm,54mm)--(67mm,59mm)--(68mm,59mm);


\draw (20mm,36mm) \squarepath;
\draw (20mm,46mm) \rectanglepath;
\draw (20mm,44mm) \rectanglepath;
\draw (24mm,44mm) \rectanglepath;
\draw (26mm,44mm) \rectanglepath;
\draw (24mm,40mm) \rectanglepath;

\draw(21mm,36mm)--(21mm,43mm)--(32mm,43mm);
\draw(23mm,36mm)--(23mm,47mm)--(32mm,47mm);
\draw(25mm,36mm)--(25mm,39mm)--(32mm,39mm);
\draw(27mm,36mm)--(27mm,41mm)--(32mm,41mm);
\draw(29mm,36mm)--(29mm,45mm)--(32mm,45mm);
\draw(31mm,36mm)--(31mm,37mm)--(32mm,37mm);

\draw (44mm,36mm) \squarepath;
\draw (44mm,46mm) \rectanglepath;
\draw (44mm,42mm) \rectanglepath;
\draw (44mm,44mm) \rectanglepath;
\draw (48mm,42mm) \rectanglepath;
\draw (48mm,44mm) \rectanglepath;

\draw(45mm,36mm)--(45mm,41mm)--(56mm,41mm);
\draw(47mm,36mm)--(47mm,47mm)--(56mm,47mm);
\draw(49mm,36mm)--(49mm,39mm)--(56mm,39mm);
\draw(51mm,36mm)--(51mm,45mm)--(56mm,45mm);
\draw(53mm,36mm)--(53mm,43mm)--(56mm,43mm);
\draw(55mm,36mm)--(55mm,37mm)--(56mm,37mm);

\draw (68mm,36mm) \squarepath;
\draw (68mm,40mm) \rectanglepath;
\draw (68mm,42mm) \rectanglepath;
\draw (68mm,44mm) \rectanglepath;
\draw (68mm,46mm) \rectanglepath;
\draw (74mm,42mm) \rectanglepath;

\draw(69mm,36mm)--(69mm,39mm)--(80mm,39mm);
\draw(71mm,36mm)--(71mm,47mm)--(80mm,47mm);
\draw(73mm,36mm)--(73mm,45mm)--(80mm,45mm);
\draw(75mm,36mm)--(75mm,41mm)--(80mm,41mm);
\draw(77mm,36mm)--(77mm,43mm)--(80mm,43mm);
\draw(79mm,36mm)--(79mm,37mm)--(80mm,37mm);

%
\draw (8mm,18mm) \squarepath;
\draw (8mm,26mm) \rectanglepath;
\draw (8mm,28mm) \rectanglepath;
\draw (10mm,26mm) \rectanglepath;
\draw (10mm,28mm) \rectanglepath;
\draw (14mm,26mm) \rectanglepath;

\draw(9mm,18mm)--(9mm,25mm)--(20mm,25mm);
\draw(11mm,18mm)--(11mm,23mm)--(20mm,23mm);
\draw(13mm,18mm)--(13mm,29mm)--(20mm,29mm);
\draw(15mm,18mm)--(15mm,21mm)--(20mm,21mm);
\draw(17mm,18mm)--(17mm,27mm)--(20mm,27mm);
\draw(19mm,18mm)--(19mm,19mm)--(20mm,19mm);

\draw (32mm,18mm) \squarepath;
\draw (32mm,28mm) \rectanglepath;
\draw (32mm,24mm) \rectanglepath;
\draw (32mm,26mm) \rectanglepath;
\draw (36mm,26mm) \rectanglepath;
\draw (38mm,26mm) \rectanglepath;

\draw(33mm,18mm)--(33mm,23mm)--(44mm,23mm);
\draw(35mm,18mm)--(35mm,29mm)--(44mm,29mm);
\draw(37mm,18mm)--(37mm,25mm)--(44mm,25mm);
\draw(39mm,18mm)--(39mm,21mm)--(44mm,21mm);
\draw(41mm,18mm)--(41mm,27mm)--(44mm,27mm);
\draw(43mm,18mm)--(43mm,19mm)--(44mm,19mm);


\draw (56mm,18mm) \squarepath;
\draw [thick](56mm,24mm) \rectanglepath;
\draw[thick] (56mm,26mm) \rectanglepath;
\draw [thick](56mm,28mm) \rectanglepath;
\draw[thick] (58mm,26mm) \rectanglepath;
\draw[thick] (58mm,28mm) \rectanglepath;

\draw(57mm,18mm)--(57mm,23mm)--(68mm,23mm);
\draw(59mm,18mm)--(59mm,25mm)--(68mm,25mm);
\draw(61mm,18mm)--(61mm,29mm)--(68mm,29mm);
\draw(63mm,18mm)--(63mm,27mm)--(68mm,27mm);
\draw(65mm,18mm)--(65mm,21mm)--(68mm,21mm);
\draw(67mm,18mm)--(67mm,19mm)--(68mm,19mm);

\draw (80mm,18mm) \squarepath;
\draw (80mm,22mm) \rectanglepath;
\draw (80mm,24mm) \rectanglepath;
\draw (80mm,26mm) \rectanglepath;
\draw (80mm,28mm) \rectanglepath;
\draw (84mm,26mm) \rectanglepath;

\draw(81mm,18mm)--(81mm,21mm)--(92mm,21mm);
\draw(83mm,18mm)--(83mm,29mm)--(92mm,29mm);
\draw(85mm,18mm)--(85mm,25mm)--(92mm,25mm);
\draw(87mm,18mm)--(87mm,27mm)--(92mm,27mm);
\draw(89mm,18mm)--(89mm,23mm)--(92mm,23mm);
\draw(91mm,18mm)--(91mm,19mm)--(92mm,19mm);

\draw (8mm,0mm) \squarepath;
\draw [thick](8mm,8mm) \rectanglepath;\draw [thick](8mm,10mm) \rectanglepath;
\draw [thick](12mm,10mm) \rectanglepath;
\draw [thick](10mm,8mm) \rectanglepath;
\draw [thick](10mm,10mm) \rectanglepath;

\draw(9mm,0mm)--(9mm,7mm)--(20mm,7mm);
\draw(11mm,0mm)--(11mm,5mm)--(20mm,5mm);
\draw(13mm,0mm)--(13mm,9mm)--(20mm,9mm);
\draw(15mm,0mm)--(15mm,11mm)--(20mm,11mm);
\draw(17mm,0mm)--(17mm,3mm)--(20mm,3mm);
\draw(19mm,0mm)--(19mm,1mm)--(20mm,1mm);

\draw (32mm,0mm) \squarepath;
\draw [thick](32mm,6mm) \rectanglepath;
\draw [thick](32mm,8mm) \rectanglepath;
\draw [thick](32mm,10mm) \rectanglepath;
\draw [thick](34mm,10mm) \rectanglepath;
\draw [thick](36mm,10mm) \rectanglepath;

\draw(33mm,0mm)--(33mm,5mm)--(44mm,5mm);
\draw(35mm,0mm)--(35mm,9mm)--(44mm,9mm);
\draw(37mm,0mm)--(37mm,7mm)--(44mm,7mm);
\draw(39mm,0mm)--(39mm,11mm)--(44mm,11mm);
\draw(41mm,0mm)--(41mm,3mm)--(44mm,3mm);
\draw(43mm,0mm)--(43mm,1mm)--(44mm,1mm);

\draw (80mm,0mm) \squarepath;
\draw[thick] (80mm,4mm) \rectanglepath;
\draw [thick](80mm,6mm) \rectanglepath;
\draw [thick](80mm,8mm) \rectanglepath;
\draw [thick](80mm,10mm) \rectanglepath;
\draw[thick] (82mm,10mm) \rectanglepath;

\draw(81mm,0mm)--(81mm,3mm)--(92mm,3mm);
\draw(83mm,0mm)--(83mm,9mm)--(92mm,9mm);
\draw(85mm,0mm)--(85mm,11mm)--(92mm,11mm);
\draw(87mm,0mm)--(87mm,7mm)--(92mm,7mm);
\draw(89mm,0mm)--(89mm,5mm)--(92mm,5mm);
\draw(91mm,0mm)--(91mm,1mm)--(92mm,1mm);

\draw(50mm,71mm)--(62mm,67mm);
\draw(50mm,71mm)--(38mm,67mm);
\draw(38mm,53mm)--(26mm,49mm);
\draw(38mm,53mm)--(50mm,49mm);
\draw(62mm,53mm)--(74mm,49mm);

\draw(26mm,35mm)--(38mm,31mm);
\draw(26mm,35mm)--(14mm,31mm);
\draw(50mm,35mm)--(62mm,31mm);
\draw(74mm,35mm)--(86mm,31mm);

\draw(14mm,17mm)--(14mm,13mm);
\draw(38mm,17mm)--(38mm,13mm);
\draw(86mm,17mm)--(86mm,13mm);

\end{tikzpicture}
\caption{The Rothe pipedream version of Figure \ref{MLS-tree}.}\label{mMLS}
\end{center}
\end{figure}
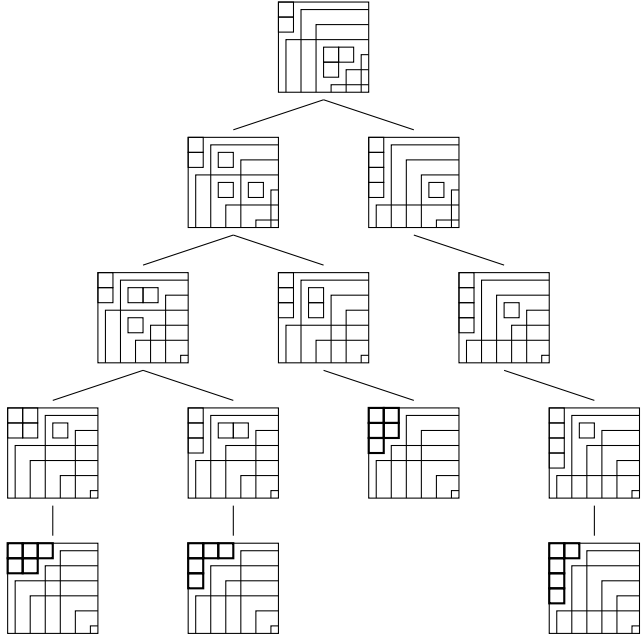

In the modified LS-tree of $w$, if replacing  the label of each node with its corresponding
Rothe pipedream, then we obtain a tree labeled with   Rothe pipedreams.
Figure \ref{mMLS} is the  Rothe pipedream version of Figure \ref{MLS-tree}.
 By using Rothe pipedreams to label the nodes, we  can directly read off the Schur functions in the expansion of
a Stanley symmetric function. Moreover,
it will be more convenient to use the    Rothe pipedream version
of a modified LS-tree  in the construction of
an EG-tree in Section \ref{PP-6}.

The remaining of this section are devoted to several  propositions
 that will be used in Section \ref{PP-6} to construct
 the EG-tree of $w$ from the modified  LS-tree of $w$.

Propositions \ref{LAR} and \ref{droop-app} tell us how to use
the droop operation to generate
the children of a node in a modified LS-tree. To describe these two
propositions, we define the notion of a pivot  of an empty box
in the Rothe pipedream of a permutation $w$. Note that for the specific empty box
$(r, w_s)$, where $r$ and $s$ are defined in Section \ref{lstree},
the notion of pivots  has been defined by
Knutson and Yong \cite{KY}.
Let $(i,j)$ be an empty box in the Rothe pipedream $D(w)$ of $w$.
We say that
an SE elbow $e$ of  $D(w)$  is
 a pivot of   $(i,j)$ if $e$ is northwest of   $(i,j)$ and there are no   other elbows contained in the rectangle with northwest corner $e$   and southeast corner $(i, j)$. In other words,
  among the SE elbows that are northwest of $(i,j)$,
the pivots are   maximally southeast.
For example, the Rothe pipedream   of
$w= 2761453$ is displayed in Figure \ref{bumpless}(a).
The empty box $(3,1)$ has no pivots, whereas the empty box
$(6,3)$ has  two pivots: $(1, 2)$ and $(4,1)$.

For a non-dominant permutation $w$, it is easily seen  that
$(p, w_q)$ is an empty box of  $D(w)$.
The following proposition gives a characterization of  the pivots of
this specific empty box.

\begin{prop}\label{LAR}
Let $w$ be a non-dominant permutation.
 Then the set of pivots of the empty box $(p, w_q)$ is
\begin{equation}\label{PUE}
\{(i, w_i)\,|\, i\in I(wt_{p,q}, p)\},
\end{equation}
where $I(wt_{p,q}, p)$ is the index set as defined in \eqref{LP}.
\end{prop}

\pf
 The assertion follows from  the definition of
 pivots as well as  the definition  of the index set
$I(wt_{p,q}, p)$.
\qed

Clearly, we can apply a droop
to the Rothe pipedream  of $w$
with respect to the empty box $(p, w_q)$ and a pivot
of $(p, w_q)$.

\begin{prop}\label{droop-app}
Assume that $w$ is a non-dominant permutation.
Let $P$ be the bumpless pipedream of $w$ obtained from the Rothe pipedream  $D(w)$ by applying a droop swapping the empty box $(p, w_q)$ and a pivot $(i, w_i)$, where $i\in I(wt_{p,q}, p)$. Then the Rothe diagram
of $w'=wt_{p,q}t_{i,p}\in \Phi(wt_{p,q}, p)$
equals the set of empty boxes of $P$.
\end{prop}

\pf
Let us first explain that the Rothe diagram  $\mathrm{Rothe}(w')$ of $w'$
can be  constructed  from the Rothe diagram  $\mathrm{Rothe}(w)$ of $w$
as follows.
First, note that $\mathrm{Rothe}(wt_{p,q})$  can be obtained from $\mathrm{Rothe}(w)$ by deleting the box $(p, w_q)$, see Figure \ref{Thm-111} for an illustration.
\begin{figure}[h]
\begin{center}

\begin{tabular}{cc}
\begin{tikzpicture}[scale = 1.1]

\def\rectanglepath{-- +(3mm,0mm) -- +(3mm,3mm) -- +(0mm,3mm) -- cycle}

\def\squarepath{-- +(27mm,0mm) -- +(27mm,27mm) -- +(0mm,27mm) -- cycle}

\def\srectanglepath{[fill=gray!35!white,draw=black!20!black]-- +(3mm,0mm) -- +(3mm,3mm) -- +(0mm,3mm) -- cycle}
\def\ssrectanglepath{[fill=gray!80!white,draw=black!20!black]-- +(3mm,0mm) -- +(3mm,3mm) -- +(0mm,3mm) -- cycle}

\draw (0mm,0mm) \squarepath;
\draw [step=3mm,dotted] (0mm,0mm)grid (27mm,27mm);

\draw (0mm,24mm) \srectanglepath;
\draw (0mm,21mm) \srectanglepath;
\draw (0mm,18mm) \srectanglepath;
\draw (0mm,15mm) \srectanglepath;
\draw (0mm,12mm) \srectanglepath;
\draw (0mm,9mm) \srectanglepath;
\draw (3mm,24mm) \srectanglepath;
\draw (6mm,24mm) \srectanglepath;
\draw (6mm,18mm) \srectanglepath;
\draw (6mm,15mm) \srectanglepath;
\draw (6mm,12mm) \ssrectanglepath;
\draw (12mm,18mm) \srectanglepath;
\draw (12mm,15mm) \srectanglepath;
\draw (12mm,12mm) \ssrectanglepath;
\draw (15mm,18mm) \ssrectanglepath;
\draw (15mm,15mm) \ssrectanglepath;
\draw (15mm,12mm) \ssrectanglepath;

\draw(1.5mm,0mm)--(1.5mm,7.5mm)--(27mm,7.5mm);
\draw(4.5mm,0mm)--(4.5mm,22.5mm)--(27mm,22.5mm);
\draw(7.5mm,0mm)--(7.5mm,10.5mm)--(27mm,10.5mm);
\draw(10.5mm,0mm)--(10.5mm,25.5mm)--(27mm,25.5mm);
\draw(13.5mm,0mm)--(13.5mm,4.5mm)--(27mm,4.5mm);
\draw(16.5mm,0mm)--(16.5mm,1.5mm)--(27mm,1.5mm);
\draw(19.5mm,0mm)--(19.5mm,19.5mm)--(27mm,19.5mm);
\draw(22.5mm,0mm)--(22.5mm,16.5mm)--(27mm,16.5mm);
\draw(25.5mm,0mm)--(25.5mm,13.5mm)--(27mm,13.5mm);

\node at (38.5mm,12.5mm) {$\xlongrightarrow[\qquad]{}$};
\end{tikzpicture}&

\begin{tikzpicture}[scale = 1.1]

\def\rectanglepath{-- +(3mm,0mm) -- +(3mm,3mm) -- +(0mm,3mm) -- cycle}

\def\squarepath{-- +(27mm,0mm) -- +(27mm,27mm) -- +(0mm,27mm) -- cycle}

\def\srectanglepath{[fill=gray!35!white,draw=black!20!black]-- +(3mm,0mm) -- +(3mm,3mm) -- +(0mm,3mm) -- cycle}

\def\ssrectanglepath{[fill=gray!80!white,draw=black!20!black]-- +(3mm,0mm) -- +(3mm,3mm) -- +(0mm,3mm) -- cycle}

\draw (0mm,0mm) \squarepath;
\draw [step=3mm,dotted] (0mm,0mm)grid (27mm,27mm);

\draw (0mm,24mm) \srectanglepath;
\draw (0mm,21mm) \srectanglepath;
\draw (0mm,18mm) \srectanglepath;
\draw (0mm,15mm) \srectanglepath;
\draw (0mm,12mm) \srectanglepath;
\draw (0mm,9mm) \srectanglepath;
\draw (3mm,24mm) \srectanglepath;
\draw (6mm,24mm) \srectanglepath;
\draw (6mm,18mm) \srectanglepath;
\draw (6mm,15mm) \srectanglepath;
\draw (6mm,12mm) \ssrectanglepath;
\draw (12mm,18mm) \srectanglepath;
\draw (12mm,15mm) \srectanglepath;
\draw (12mm,12mm) \ssrectanglepath;
\draw (15mm,18mm) \ssrectanglepath;
\draw (15mm,15mm) \ssrectanglepath;

\draw(1.5mm,0mm)--(1.5mm,7.5mm)--(27mm,7.5mm);
\draw(4.5mm,0mm)--(4.5mm,22.5mm)--(27mm,22.5mm);
\draw(7.5mm,0mm)--(7.5mm,10.5mm)--(27mm,10.5mm);
\draw(10.5mm,0mm)--(10.5mm,25.5mm)--(27mm,25.5mm);
\draw(13.5mm,0mm)--(13.5mm,4.5mm)--(27mm,4.5mm);
\draw(16.5mm,0mm)--(16.5mm,13.5mm)--(27mm,13.5mm);
\draw(19.5mm,0mm)--(19.5mm,19.5mm)--(27mm,19.5mm);
\draw(22.5mm,0mm)--(22.5mm,16.5mm)--(27mm,16.5mm);
\draw(25.5mm,0mm)--(25.5mm,1.5mm)--(27mm,1.5mm);


\end{tikzpicture}
\end{tabular}

\caption{Construction of $\mathrm{Rothe}(wt_{p,q})$ from $\mathrm{Rothe}(w)$.}\label{Thm-111}
\end{center}
\end{figure}
Let us proceed to determine $\mathrm{Rothe}(w')$ from $\mathrm{Rothe}(wt_{p,q})$.
Locate the rectangle $R$  in the $n\times n$ square grid,
such that the  northwest corner is   $(i, w_i)$ and the southeast corner is   $(p, w_q)$. Then  $\mathrm{Rothe}(w')$ can be obtained from $\mathrm{Rothe}(wt_{p,q})$ as follows:
\begin{itemize}

\item[(1)]
Move each box of $\mathrm{Rothe}(wt_{p,q})$ in the bottom row of $R$   to the top row of $R$;

 \item[(2)]
  Move each box of $\mathrm{Rothe}(wt_{p,q})$ in the rightmost column of $R$ to the leftmost column of $R$;

 \item[(3)] Add a new box to the position $(i, w_i)$.

\end{itemize}
Figure \ref{Thm-122} illustrates
the above construction of $\mathrm{Rothe}(w')$ from $\mathrm{Rothe}(wt_{p,q})$.

\begin{figure}[h]
\begin{center}

\begin{tabular}{cccc}

\begin{picture}(105,100)
\begin{tikzpicture}[scale = 1]

\def\rectanglepath{-- +(3mm,0mm) -- +(3mm,3mm) -- +(0mm,3mm) -- cycle}

\def\squarepath{-- +(27mm,0mm) -- +(27mm,27mm) -- +(0mm,27mm) -- cycle}

\def\srectanglepath{[fill=gray!35!white,draw=black!20!black]-- +(3mm,0mm) -- +(3mm,3mm) -- +(0mm,3mm) -- cycle}

\def\ssrectanglepath{[fill=gray!80!white,draw=black!20!black]-- +(3mm,0mm) -- +(3mm,3mm) -- +(0mm,3mm) -- cycle}

\draw (0mm,0mm) \squarepath;
\draw [step=3mm,dotted] (0mm,0mm)grid (27mm,27mm);

\draw (0mm,24mm) \srectanglepath;
\draw (0mm,21mm) \srectanglepath;
\draw (0mm,18mm) \srectanglepath;
\draw (0mm,15mm) \srectanglepath;
\draw (0mm,12mm) \srectanglepath;
\draw (0mm,9mm) \srectanglepath;
\draw (3mm,24mm) \srectanglepath;
\draw (6mm,24mm) \srectanglepath;
\draw (6mm,18mm) \srectanglepath;
\draw (6mm,15mm) \srectanglepath;
\draw (6mm,12mm) \ssrectanglepath;
\draw (12mm,18mm) \srectanglepath;
\draw (12mm,15mm) \srectanglepath;
\draw (12mm,12mm) \ssrectanglepath;
\draw (15mm,18mm) \ssrectanglepath;
\draw (15mm,15mm) \ssrectanglepath;

\draw(1.5mm,0mm)--(1.5mm,7.5mm)--(27mm,7.5mm);
\draw(4.5mm,0mm)--(4.5mm,22.5mm)--(27mm,22.5mm);
\draw(7.5mm,0mm)--(7.5mm,10.5mm)--(27mm,10.5mm);
\draw(10.5mm,0mm)--(10.5mm,25.5mm)--(27mm,25.5mm);
\draw(13.5mm,0mm)--(13.5mm,4.5mm)--(27mm,4.5mm);
\draw(16.5mm,0mm)--(16.5mm,13.5mm)--(27mm,13.5mm);
\draw(19.5mm,0mm)--(19.5mm,19.5mm)--(27mm,19.5mm);
\draw(22.5mm,0mm)--(22.5mm,16.5mm)--(27mm,16.5mm);
\draw(25.5mm,0mm)--(25.5mm,1.5mm)--(27mm,1.5mm);

\node at (34.5mm,12.5mm) {$\longrightarrow$};

\end{tikzpicture}
\end{picture}
&

\begin{picture}(100,99)
\begin{tikzpicture}[scale = 1]

\def\rectanglepath{-- +(3mm,0mm) -- +(3mm,3mm) -- +(0mm,3mm) -- cycle}

\def\squarepath{-- +(27mm,0mm) -- +(27mm,27mm) -- +(0mm,27mm) -- cycle}

\def\srectanglepath{[fill=gray!35!white,draw=black!20!black]-- +(3mm,0mm) -- +(3mm,3mm) -- +(0mm,3mm) -- cycle}

\def\ssrectanglepath{[fill=gray!80!white,draw=black!20!black]-- +(3mm,0mm) -- +(3mm,3mm) -- +(0mm,3mm) -- cycle}

\draw (0mm,0mm) \squarepath;
\draw [step=3mm,dotted] (0mm,0mm)grid (27mm,27mm);

\draw (0mm,24mm) \srectanglepath;
\draw (0mm,21mm) \srectanglepath;
\draw (0mm,18mm) \srectanglepath;
\draw (0mm,15mm) \srectanglepath;
\draw (0mm,12mm) \srectanglepath;
\draw (0mm,9mm) \srectanglepath;
\draw (3mm,24mm) \srectanglepath;
\draw (6mm,24mm) \srectanglepath;
\draw (6mm,21mm) \ssrectanglepath;
\draw (6mm,18mm) \srectanglepath;
\draw (6mm,15mm) \srectanglepath;
\draw (12mm,21mm) \ssrectanglepath;
\draw (12mm,18mm) \srectanglepath;
\draw (12mm,15mm) \srectanglepath;
\draw (15mm,18mm) \ssrectanglepath;
\draw (15mm,15mm) \ssrectanglepath;

\node at (32.5mm,12.5mm) {$\longrightarrow$};
\end{tikzpicture}
\end{picture}&

\begin{picture}(100,99)
\begin{tikzpicture}[scale = 1]

\def\rectanglepath{-- +(3mm,0mm) -- +(3mm,3mm) -- +(0mm,3mm) -- cycle}

\def\squarepath{-- +(27mm,0mm) -- +(27mm,27mm) -- +(0mm,27mm) -- cycle}

\def\srectanglepath{[fill=gray!35!white,draw=black!20!black]-- +(3mm,0mm) -- +(3mm,3mm) -- +(0mm,3mm) -- cycle}

\def\ssrectanglepath{[fill=gray!80!white,draw=black!20!black]-- +(3mm,0mm) -- +(3mm,3mm) -- +(0mm,3mm) -- cycle}

\draw (0mm,0mm) \squarepath;
\draw [step=3mm,dotted] (0mm,0mm)grid (27mm,27mm);

\draw (0mm,24mm) \srectanglepath;
\draw (0mm,21mm) \srectanglepath;
\draw (0mm,18mm) \srectanglepath;
\draw (0mm,15mm) \srectanglepath;
\draw (0mm,12mm) \srectanglepath;
\draw (0mm,9mm) \srectanglepath;
\draw (3mm,24mm) \srectanglepath;
\draw (3mm,18mm) \ssrectanglepath;
\draw (3mm,15mm) \ssrectanglepath;
\draw (6mm,24mm) \srectanglepath;
\draw (6mm,21mm) \ssrectanglepath;
\draw (6mm,18mm) \srectanglepath;
\draw (6mm,15mm) \srectanglepath;

\draw (12mm,21mm) \ssrectanglepath;
\draw (12mm,18mm) \srectanglepath;
\draw (12mm,15mm) \srectanglepath;

\node at (32.5mm,12.5mm) {$\longrightarrow$};
\end{tikzpicture}
\end{picture}&

\begin{tikzpicture}[scale = 1]

\def\rectanglepath{-- +(3mm,0mm) -- +(3mm,3mm) -- +(0mm,3mm) -- cycle}

\def\squarepath{-- +(27mm,0mm) -- +(27mm,27mm) -- +(0mm,27mm) -- cycle}

\def\srectanglepath{[fill=gray!35!white,draw=black!20!black]-- +(3mm,0mm) -- +(3mm,3mm) -- +(0mm,3mm) -- cycle}

\def\ssrectanglepath{[fill=gray!80!white,draw=black!20!black]-- +(3mm,0mm) -- +(3mm,3mm) -- +(0mm,3mm) -- cycle}

\draw (0mm,0mm) \squarepath;
\draw [step=3mm,dotted] (0mm,0mm)grid (27mm,27mm);

\draw (0mm,24mm) \srectanglepath;
\draw (0mm,21mm) \srectanglepath;
\draw (0mm,18mm) \srectanglepath;
\draw (0mm,15mm) \srectanglepath;
\draw (0mm,12mm) \srectanglepath;
\draw (0mm,9mm) \srectanglepath;
\draw (3mm,24mm) \srectanglepath;
\draw (3mm,21mm) \ssrectanglepath;
\draw (3mm,18mm) \ssrectanglepath;
\draw (3mm,15mm) \ssrectanglepath;
\draw (6mm,24mm) \srectanglepath;
\draw (6mm,21mm) \ssrectanglepath;
\draw (6mm,18mm) \srectanglepath;
\draw (6mm,15mm) \srectanglepath;

\draw (12mm,21mm) \ssrectanglepath;
\draw (12mm,18mm) \srectanglepath;
\draw (12mm,15mm) \srectanglepath;

\draw(1.5mm,0mm)--(1.5mm,7.5mm)--(27mm,7.5mm);
\draw(4.5mm,0mm)--(4.5mm,13.5mm)--(27mm,13.5mm);
\draw(7.5mm,0mm)--(7.5mm,10.5mm)--(27mm,10.5mm);
\draw(10.5mm,0mm)--(10.5mm,25.5mm)--(27mm,25.5mm);
\draw(13.5mm,0mm)--(13.5mm,4.5mm)--(27mm,4.5mm);
\draw(16.5mm,0mm)--(16.5mm,22.5mm)--(27mm,22.5mm);
\draw(19.5mm,0mm)--(19.5mm,19.5mm)--(27mm,19.5mm);
\draw(22.5mm,0mm)--(22.5mm,16.5mm)--(27mm,16.5mm);
\draw(25.5mm,0mm)--(25.5mm,1.5mm)--(27mm,1.5mm);

\end{tikzpicture}
\end{tabular}

\caption{Construction of $\mathrm{Rothe}(w')$ from $\mathrm{Rothe}(wt_{p,q})$.}\label{Thm-122}
\end{center}
\end{figure}

One the other hand, one can apply a droop on $D(w)$  swapping
$(i, w_i)$ and $(p, w_q)$ to obtain a bumpless pipedream $P$ of $w$, see Figure \ref{fig-mar-gen} for an illustration.

\begin{figure}[h]
\begin{center}

\begin{tabular}{cc}
\begin{tikzpicture}[scale = 1.1]

\def\rectanglepath{-- +(3mm,0mm) -- +(3mm,3mm) -- +(0mm,3mm) -- cycle}

\def\squarepath{-- +(27mm,0mm) -- +(27mm,27mm) -- +(0mm,27mm) -- cycle}

\def\srectanglepath{[fill=gray!35!white,draw=black!20!black]-- +(3mm,0mm) -- +(3mm,3mm) -- +(0mm,3mm) -- cycle}

\draw (0mm,0mm) \squarepath;
\draw [step=3mm,dotted] (0mm,0mm)grid (27mm,27mm);

\draw (0mm,24mm) \srectanglepath;
\draw (0mm,21mm) \srectanglepath;
\draw (0mm,18mm) \srectanglepath;
\draw (0mm,15mm) \srectanglepath;
\draw (0mm,12mm) \srectanglepath;
\draw (0mm,9mm) \srectanglepath;
\draw (3mm,24mm) \srectanglepath;
\draw (6mm,24mm) \srectanglepath;
\draw (6mm,18mm) \srectanglepath;
\draw (6mm,15mm) \srectanglepath;
\draw (6mm,12mm) \srectanglepath;
\draw (12mm,18mm) \srectanglepath;
\draw (12mm,15mm) \srectanglepath;
\draw (12mm,12mm) \srectanglepath;
\draw (15mm,18mm) \srectanglepath;
\draw (15mm,15mm) \srectanglepath;
\draw (15mm,12mm) \srectanglepath;

\draw(1.5mm,0mm)--(1.5mm,7.5mm)--(27mm,7.5mm);
\draw(4.5mm,0mm)--(4.5mm,22.5mm)--(27mm,22.5mm);
\draw[-][ultra thick](4.5mm,0mm)--(4.5mm,22.5mm)--(27mm,22.5mm);
\draw(7.5mm,0mm)--(7.5mm,10.5mm)--(27mm,10.5mm);
\draw(10.5mm,0mm)--(10.5mm,25.5mm)--(27mm,25.5mm);
\draw(13.5mm,0mm)--(13.5mm,4.5mm)--(27mm,4.5mm);
\draw(16.5mm,0mm)--(16.5mm,1.5mm)--(27mm,1.5mm);
\draw(19.5mm,0mm)--(19.5mm,19.5mm)--(27mm,19.5mm);
\draw(22.5mm,0mm)--(22.5mm,16.5mm)--(27mm,16.5mm);
\draw(25.5mm,0mm)--(25.5mm,13.5mm)--(27mm,13.5mm);

\node at (38.5mm,12.5mm) {$\xlongrightarrow[\qquad]{}$};
\end{tikzpicture}&

\begin{tikzpicture}[scale = 1.1]

\def\rectanglepath{-- +(3mm,0mm) -- +(3mm,3mm) -- +(0mm,3mm) -- cycle}

\def\squarepath{-- +(27mm,0mm) -- +(27mm,27mm) -- +(0mm,27mm) -- cycle}

\def\srectanglepath{[fill=gray!35!white,draw=black!20!black]-- +(3mm,0mm) -- +(3mm,3mm) -- +(0mm,3mm) -- cycle}

\draw (0mm,0mm) \squarepath;
\draw [step=3mm,dotted] (0mm,0mm)grid (27mm,27mm);

\draw (0mm,24mm) \srectanglepath;
\draw (0mm,21mm) \srectanglepath;
\draw (0mm,18mm) \srectanglepath;
\draw (0mm,15mm) \srectanglepath;
\draw (0mm,12mm) \srectanglepath;
\draw (0mm,9mm) \srectanglepath;
\draw (3mm,24mm) \srectanglepath;
\draw (3mm,21mm) \srectanglepath;
\draw (6mm,24mm) \srectanglepath;
\draw (6mm,18mm) \srectanglepath;
\draw (6mm,15mm) \srectanglepath;
\draw (6mm,21mm) \srectanglepath;
\draw (12mm,18mm) \srectanglepath;
\draw (12mm,15mm) \srectanglepath;
\draw (12mm,21mm) \srectanglepath;
\draw (3mm,18mm) \srectanglepath;
\draw (3mm,15mm) \srectanglepath;

\draw(1.5mm,0mm)--(1.5mm,7.5mm)--(27mm,7.5mm);
\draw(4.5mm,0mm)--(4.5mm,13.5mm)--(16.5mm,13.5mm)--(16.5mm,22.5mm)--(27mm,22.5mm);
\draw[-][ultra thick](4.5mm,0mm)--(4.5mm,13.5mm)--(16.5mm,13.5mm)--(16.5mm,22.5mm)--(27mm,22.5mm);

\draw(7.5mm,0mm)--(7.5mm,10.5mm)--(27mm,10.5mm);
\draw(10.5mm,0mm)--(10.5mm,25.5mm)--(27mm,25.5mm);
\draw(13.5mm,0mm)--(13.5mm,4.5mm)--(27mm,4.5mm);
\draw(16.5mm,0mm)--(16.5mm,1.5mm)--(27mm,1.5mm);
\draw(19.5mm,0mm)--(19.5mm,19.5mm)--(27mm,19.5mm);
\draw(22.5mm,0mm)--(22.5mm,16.5mm)--(27mm,16.5mm);
\draw(25.5mm,0mm)--(25.5mm,13.5mm)--(27mm,13.5mm);

\end{tikzpicture}
\end{tabular}

\caption{A  droop operation on $D(w)$.}\label{fig-mar-gen}
\end{center}
\end{figure}
Evidently,  the droop
has the same effect on the empty boxes of $D(w)$ (that is,
the Rothe diagram of $w$) as the
operation illustrated  in Figures \ref{Thm-111} and \ref{Thm-122}.
This competes the proof.
\qed

\begin{re}
When $p$ is the last descent $r$ of $w$,
Knutson and Yong \cite{KY} introduced  the marching
operation to generate the Rothe diagram of a child of $w$
from the Rothe diagram $\mathrm{Rothe}(w)$ of $w$. In this specific case, the marching operation
has the same effect on $\mathrm{Rothe}(w)$ as the droop operation.
\end{re}

Propositions \ref{LAR} and \ref{droop-app} give an explicit way
to generate the children of
a non-dominant permutation $w$ in the modified LS-tree of $w$. Assume that $e_1, \ldots, e_m$
are the pivots of $(p, w_q)$. For $1\leq i\leq m$, let $P_i$ be the
bumpless pipedream of $w$ obtained from  $D(w)$
by applying a droop
with respect to  $(p, w_q)$ and  $e_i$. Denote by $\mathrm{Box}(P_i)$
the empty boxes of $P_i$.
Then $\mathrm{Box}(P_1), \ldots, \mathrm{Box}(P_m)$ are the Rothe diagrams of
the children of $w$. Equivalently,  $\mathrm{Box}(P_1), \ldots, \mathrm{Box}(P_m)$
are the empty
boxes of Rothe pipedreams of the children of $w$.

The next two propositions investigate  the properties concerning
the specific empty box $(p, w_q)$.
Define  a total order on the boxes of the $n\times n$ grid
by letting $(i,j)<(k, \ell)$ if
\begin{equation}\label{totall}
\text{either $i<k$, or $i=k$ and $j<\ell$}.
\end{equation}
Let $\mathrm{pivot}(w)$ denote the set of empty boxes
in the Rothe pipedream  of $w$ which has at least one pivot.

\begin{prop}\label{LAR-2}
Let $w$ be a non-dominant permutation.
Then  $(p, w_q)$ is the largest box
in the set $\mathrm{pivot}(w)$ under the total order  defined
in \eqref{totall}.
\end{prop}

\pf Suppose otherwise that
 $(p', j)\neq (p, w_q)$  is
the largest box in $\mathrm{pivot}(w)$.
Then there is an index $q'>p'$ such that $w_{q'}=j$.
Since $(p, w_q)$ is the
 rightmost empty box in the $p$-th
row of $D(w)$, we must have $p'>p$. Let $(i', w_{i'})$
be a pivot of $(p', w_{q'})$, where $i'<p'$. By the definition of
a pivot, we have $w_{i'}<j=w_{q'}$.  Moreover, notice that
$w_{q'}=j<w_{p'}$. So, for the indices  $i'<p'<q'$,  we find that
$w_{i'}<w_{q'}<w_{p'}$. This implies
that $p'$ satisfies the condition \eqref{Condi},
 contrary to the choice of $p$. This concludes the proof.
\qed

Using Proposition \ref{LAR-2}, we can prove  the following assertion.

\begin{prop}\label{FORW}
Let $w$ be a non-dominant permutation, and
let   $w'=wt_{p,q}t_{i,p}\in \Phi(wt_{p,q}, p)$ be a child  of $w$.
Suppose that $w'$ is also a non-dominant permutation.
Write  $p'$ and $q'$ for the indices as defined in \eqref{PPPP}
and \eqref{QQQQ} for $w'$, respectively. Then
\begin{equation}\label{III}
(p', w_{q'}')< (p,w_q)
\end{equation}
under the total order defined in \eqref{totall}.
\end{prop}

\begin{figure}[h]
\begin{center}

\begin{tikzpicture}[scale = .85]
\def\rectanglepath{-- +(4mm,0mm) -- +(4mm,4mm) -- +(0mm,4mm) -- cycle}

\def\squarepath{-- +(32mm,0mm) -- +(32mm,32mm) -- +(0mm,32mm) -- cycle}

\draw (0mm,0mm) \squarepath;
\draw (12mm,12mm)\rectanglepath;
\draw(6mm,0mm)--(6mm,28mm)--(32mm,28mm);
\draw(14mm,0mm)--(14mm,6mm)--(32mm,6mm);
\draw(22mm,0mm)--(22mm,14mm)--(32mm,14mm);
\draw[dashed,line width=1.2pt](4mm,12mm)--(16mm,12mm)--(16mm,30mm)--(4mm,30mm)--cycle;
\node at (35mm,28mm) {\small{$w_i$}};
\node at (35mm,14mm) {\small{$w_p$}};
\node at (35mm,6mm) {\small{$w_q$}};
\node at (6mm,-3mm) {\small{$w_i$}};
\node at (14mm,-3mm) {\small{$w_q$}};
\node at (22mm,-3mm) {\small{$w_p$}};


\filldraw[fill=gray!35!white,draw=gray!20!white]
(80mm,12mm)--(84mm,12mm)--(84mm,16mm)--(112mm,16mm)--(112mm,32mm)--(80mm,32mm)--cycle;
\draw (80mm,0mm) \squarepath;
\draw(86mm,0mm)--(86mm,14mm)--(112mm,14mm);
\draw(94mm,0mm)--(94mm,28mm)--(112mm,28mm);
\draw(102mm,0mm)--(102mm,6mm)--(112mm,6mm);
\draw[dashed,line width=1.2pt](84mm,12mm)--(96mm,12mm)--(96mm,30mm)--(84mm,30mm)--cycle;

\node at (115mm,28mm) {\small{$w_q$}};
\node at (115mm,14mm) {\small{$w_i$}};
\node at (115mm,6mm) {\small{$w_p$}};
\node at (86mm,-3mm) {\small{$w_i$}};
\node at (94mm,-3mm) {\small{$w_q$}};
\node at (102mm,-3mm) {\small{$w_p$}};

%
\end{tikzpicture}

\caption{Rothe pipedreams of $w$ and $w'$.}\label{Sha-1-1}
\end{center}
\end{figure}
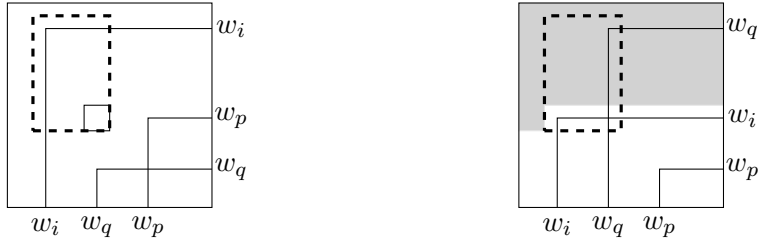

\pf
The proof is best understood by means of the pictures of the
Rothe pipedreams of $w$ and $w'$, as illustrated in
Figure \ref{Sha-1-1}.
The dashed rectangle in Figure \ref{Sha-1-1}
 signifies the rectangle
 where the droop operation
 on  $D(w)$, as described  in Proposition \ref{droop-app},  took place.
The shaded area of $D(w')$ in Figure \ref{Sha-1-1} contains the boxes of $D(w')$
smaller than the SE elbow $(p,w_i)$ under
the total order  \eqref{totall}.
Since $w'_p=w_i$,   there are no empty boxes in the $p$-th row
of $D(w')$ after the box $(p, w_i)$. So \eqref{III} is equivalent to saying
that the empty box $(p', w_{q'}')$ of $D(w')$ lies in the shaded area.

Suppose otherwise that $(p', w_{q'}')> (p,w_q)$.
Then $(p', w_{q'}')$
lies strictly below   row $p$.
Note  that $(p', w_{q'}')$
is also an empty box of $D(w)$.
There are two cases to discuss.

\noindent
Case 1:  $w_i<w_{q'}'$. In this case, in the Rothe
pipedream of $w$, $(i,w_i)$ is an SE elbow that is northwest of
$(p', w_{q'}')$. Hence  $(p', w_{q'}')$ has a pivot in  the Rothe
pipedream of $w$, which  is contrary to Proposition \ref{LAR-2}.

\noindent
Case 2:  $w_i>w_{q'}'$. Let $e$ be an SE elbow in the Rothe pipedream
of $w'$ which is a pivot of $(p', w_{q'}')$. Still, from Figure \ref{Sha-1-1},
we see that $e$ is also an
SE elbow in the Rothe pipedream
of $w$. This implies that  $(p', w_{q'}')$ has a pivot in  the Rothe
pipedream of $w$, which   again
contradicts  Proposition \ref{LAR-2}.
\qed

\section{Edelman-Greene tree}\label{PP-6}

In this section, we introduce  the structure of
the Edelman-Greene tree (EG-tree) of a permutation $w$.
Each node in the EG-tree of $w$ is labeled with a bumpless pipedream of $w$.
In particular, we  show that the leaves in the
EG-tree of $w$ are  exactly  labeled with  the EG-pipedreams of $w$.

Let us proceed with the construction of the EG-tree of  $w$. The idea behind is that for each path in the modified LS-tree of $w$, we construct a corresponding path of bumpless pipedreams of $w$.
Assume that
\[w=w^{(0)}\rightarrow w^{(1)}\rightarrow \cdots \rightarrow w^{(m)}\]
is a path in the modified LS-tree of $w$ from the root $w^{(0)}=w$ to a leaf $w^{(m)}$.  We  aim to construct a sequence
\begin{equation}\label{EG-22}
P_0\rightarrow P_1 \rightarrow \cdots \rightarrow P_m
\end{equation}
of  bumpless pipedreams of $w$ such that
\begin{itemize}
\item[(1)]
$P_0$ is the Rothe pipedream of $w$ and $P_m$ is an EG-pipedream of $w$;

\item[(2)]
For $0\leq j\leq m$, the empty boxes of $P_j$ are the same as
the empty boxes of the Rothe pipedream $D(w^{(j+1)})$  of $w^{(j)}$.
\end{itemize}

We  now describe the construction of the pipedreams in \eqref{EG-22}.
For $0\le j\le m-1$, let $p_j$ and $q_j$ be  the indices of the permutation $w^{(j)}$ as defined in \eqref{PPPP} and
\eqref{QQQQ}. By the construction of the modified LS-tree of $w$, there exists an index $i_j\in I(w^{(j)}t_{p_j, q_j}, p_j)$
such that $w^{(j+1)}=w^{(j)} t_{p_j, q_j} t_{i_j, p_j}$.
By Proposition \ref{droop-app}, we can generate the empty boxes of
$D(w^{(j+1)})$ by applying a
droop  to   $ D(w^{(j)})$ with respect to the
empty box $(p_j,w^{(j)}_{q_j})$ and the
pivot $(i_j,w^{(j)}_{i_j})$. As will be seen in the proof of  Proposition
\ref{LSB},  we can also  apply a droop to $P_j$ with respect to
 $(p_j,w^{(j)}_{q_j})$ and
$(i_j,w^{(j)}_{i_j})$.
Let $P_{j+1}$ be the bumpless pipedream of $w$
obtained by applying a droop to $P_{j}$ with respect to
$(p_j,w^{(j)}_{q_j})$ and  $(i_j,w^{(j)}_{i_j})$.

\begin{prop}\label{LSB}
The above construction  from  $P_j$ to $P_{j+1}$ is feasible.
\end{prop}

\pf
Let us first consider the case $j=0$.  Since $P_0=D(w)$,
by Proposition  \ref{droop-app},
we can apply a droop
to $P_0$ by swapping the empty box $(p_0,w^{(0)}_{q_0})$ and its pivot $(i_0,w^{(0)}_{i_0})$, resulting in the bumpless pipedream
 $P_1$.

We proceed to consider the case $j=1$.
By Proposition  \ref{droop-app},
the empty boxes of  $D(w^{(2)})$ can be obtained
from  $D(w^{(1)})$ by applying a droop
with respect to the empty box $(p_1,w^{(1)}_{q_1})$
and its pivot $(i_1,w^{(1)}_{i_1})$. By Proposition \ref{FORW},
there holds that $(p_1,w^{(1)}_{p_1})<(p_0, w^{(0)}_{p_0})$.
By the construction
of a droop operation, it is easy to check that $P_1$ and $D(w^{(1)})$
have the same tile at each   position smaller than  $(p_0,w^{(0)}_{q_0})$.
Thus we can apply a droop to $P_1$ with respect to   $(p_1,w^{(1)}_{q_1})$ and $(i_1,w^{(1)}_{i_1})$, yielding  the bumpless pipedream $P_2$ of $w$.

For the same reason as above, in the general case for $j\geq 2$,
we can check that $P_j$ and $D(w^{(j)})$
have the same tile at each   position
smaller than  $(p_{j-1},w^{(j-1)}_{q_{j-1}})$, and so that
we can apply a droop to $P_j$ with respect to   $(p_j,w^{(1)}_{q_j})$ and $(i_j,w^{(j)}_{i_j})$. This generates the bumpless pipedream  $P_{j+1}$ of $w$.
\qed

By the proof of  Proposition
\ref{LSB},
the empty boxes of $P_j$ are the same as
the empty boxes of the Rothe pipedream
of $w^{(j)}$.
Notice that the Rothe diagram of a dominant permutation
is  a partition shape at the northwest corner of the $n\times n$
grid, see Fulton \cite{Ful} or Stanley \cite[Chapter 1]{Sta-5}.
Thus the bumpless pipedream $P_m$  is
an EG-pipedream of $w$.

\begin{figure}[!htb]
\begin{center}
\begin{tikzpicture}

\def\rectanglepath{-- +(2mm,0mm) -- +(2mm,2mm) -- +(0mm,2mm) -- cycle}

\def\squarepath{-- +(12mm,0mm) -- +(12mm,12mm) -- +(0mm,12mm) -- cycle}

\draw (44mm,72mm) \squarepath;

\draw (44mm,82mm) \rectanglepath;
\draw (44mm,80mm) \rectanglepath;
\draw (50mm,76mm) \rectanglepath;
\draw (52mm,76mm) \rectanglepath;
\draw (50mm,74mm) \rectanglepath;

\draw(45mm,72mm)--(45mm,79mm)--(56mm,79mm);
\draw(47mm,72mm)--(47mm,83mm)--(56mm,83mm);
\draw(49mm,72mm)--(49mm,81mm)--(56mm,81mm);
\draw(51mm,72mm)--(51mm,73mm)--(56mm,73mm);
\draw(53mm,72mm)--(53mm,75mm)--(56mm,75mm);
\draw(55mm,72mm)--(55mm,77mm)--(56mm,77mm);

%
\draw (32mm,54mm) \squarepath;
\draw (32mm,64mm) \rectanglepath;
\draw (32mm,62mm) \rectanglepath;
\draw (36mm,62mm) \rectanglepath;
\draw (40mm,58mm) \rectanglepath;
\draw (36mm,58mm) \rectanglepath;

\draw(33mm,54mm)--(33mm,61mm)--(44mm,61mm);
\draw(35mm,54mm)--(35mm,65mm)--(44mm,65mm);
\draw(37mm,54mm)--(37mm,57mm)--(39mm,57mm)--(39mm,63mm)--(44mm,63mm);
\draw(39mm,54mm)--(39mm,55mm)--(44mm,55mm);
\draw(41mm,54mm)--(41mm,57mm)--(44mm,57mm);
\draw(43mm,54mm)--(43mm,59mm)--(44mm,59mm);

\draw (56mm,54mm) \squarepath;
\draw (56mm,58mm) \rectanglepath;
\draw (56mm,60mm) \rectanglepath;
\draw (56mm,62mm) \rectanglepath;
\draw (56mm,64mm) \rectanglepath;
\draw (64mm,58mm) \rectanglepath;

\draw(57mm,54mm)--(57mm,57mm)--(63mm,57mm)--(63mm,61mm)--(68mm,61mm);
\draw(59mm,54mm)--(59mm,65mm)--(68mm,65mm);
\draw(61mm,54mm)--(61mm,63mm)--(68mm,63mm);
\draw(63mm,54mm)--(63mm,55mm)--(68mm,55mm);
\draw(65mm,54mm)--(65mm,57mm)--(68mm,57mm);
\draw(67mm,54mm)--(67mm,59mm)--(68mm,59mm);


\draw (20mm,36mm) \squarepath;
\draw (20mm,46mm) \rectanglepath;
\draw (20mm,44mm) \rectanglepath;
\draw (24mm,44mm) \rectanglepath;
\draw (26mm,44mm) \rectanglepath;
\draw (24mm,40mm) \rectanglepath;

\draw(21mm,36mm)--(21mm,43mm)--(32mm,43mm);
\draw(23mm,36mm)--(23mm,47mm)--(32mm,47mm);
\draw(25mm,36mm)--(25mm,39mm)--(27mm,39mm)--(27mm,41mm)--(29mm,41mm)--(29mm,45mm)--(32mm,45mm);

\draw(27mm,36mm)--(27mm,37mm)--(32mm,37mm);
\draw(29mm,36mm)--(29mm,39mm)--(32mm,39mm);
\draw(31mm,36mm)--(31mm,41mm)--(32mm,41mm);

\draw (44mm,36mm) \squarepath;
\draw (44mm,46mm) \rectanglepath;
\draw (44mm,42mm) \rectanglepath;
\draw (44mm,44mm) \rectanglepath;
\draw (48mm,42mm) \rectanglepath;
\draw (48mm,44mm) \rectanglepath;

\draw(45mm,36mm)--(45mm,41mm)--(53mm,41mm)--(53mm,43mm)--(56mm,43mm);
\draw(47mm,36mm)--(47mm,47mm)--(56mm,47mm);
\draw(49mm,36mm)--(49mm,39mm)--(51mm,39mm)--(51mm,45mm)--(56mm,45mm);

\draw(51mm,36mm)--(51mm,37mm)--(56mm,37mm);
\draw(53mm,36mm)--(53mm,39mm)--(56mm,39mm);
\draw(55mm,36mm)--(55mm,41mm)--(56mm,41mm);

\draw (68mm,36mm) \squarepath;
\draw (68mm,40mm) \rectanglepath;
\draw (68mm,42mm) \rectanglepath;
\draw (68mm,44mm) \rectanglepath;
\draw (68mm,46mm) \rectanglepath;
\draw (74mm,42mm) \rectanglepath;

\draw(69mm,36mm)--(69mm,39mm)--(75mm,39mm)--(75mm,41mm)--(77mm,41mm)--(77mm,43mm)--(80mm,43mm);
\draw(71mm,36mm)--(71mm,47mm)--(80mm,47mm);
\draw(73mm,36mm)--(73mm,45mm)--(80mm,45mm);

\draw(75mm,36mm)--(75mm,37mm)--(80mm,37mm);
\draw(77mm,36mm)--(77mm,39mm)--(80mm,39mm);
\draw(79mm,36mm)--(79mm,41mm)--(80mm,41mm);

%
\draw (8mm,18mm) \squarepath;
\draw (8mm,26mm) \rectanglepath;
\draw (8mm,28mm) \rectanglepath;
\draw (10mm,26mm) \rectanglepath;
\draw (10mm,28mm) \rectanglepath;
\draw (14mm,26mm) \rectanglepath;

\draw(9mm,18mm)--(9mm,25mm)--(20mm,25mm);
\draw(11mm,18mm)--(11mm,23mm)--(13mm,23mm)--(13mm,29mm)--(20mm,29mm);
\draw(13mm,18mm)--(13mm,21mm)--(15mm,21mm)--(15mm,23mm)--(17mm,23mm)--(17mm,27mm)--(20mm,27mm);
\draw(15mm,18mm)--(15mm,19mm)--(20mm,19mm);
\draw(17mm,18mm)--(17mm,21mm)--(20mm,21mm);
\draw(19mm,18mm)--(19mm,23mm)--(20mm,23mm);

\draw (32mm,18mm) \squarepath;
\draw (32mm,28mm) \rectanglepath;
\draw (32mm,24mm) \rectanglepath;
\draw (32mm,26mm) \rectanglepath;
\draw (36mm,26mm) \rectanglepath;
\draw (38mm,26mm) \rectanglepath;

\draw(33mm,18mm)--(33mm,23mm)--(37mm,23mm)--(37mm,25mm)--(44mm,25mm);
\draw(35mm,18mm)--(35mm,29mm)--(44mm,29mm);
\draw(37mm,18mm)--(37mm,21mm)--(39mm,21mm)--(39mm,23mm)--(41mm,23mm)--(41mm,27mm)--(44mm,27mm);
\draw(39mm,18mm)--(39mm,19mm)--(44mm,19mm);
\draw(41mm,18mm)--(41mm,21mm)--(44mm,21mm);
\draw(43mm,18mm)--(43mm,23mm)--(44mm,23mm);


\draw (56mm,18mm) \squarepath;

\draw [thick](56mm,24mm) \rectanglepath;
\draw[thick] (56mm,26mm) \rectanglepath;
\draw[thick] (56mm,28mm) \rectanglepath;
\draw[thick] (58mm,26mm) \rectanglepath;
\draw[thick] (58mm,28mm) \rectanglepath;

\draw(57mm,18mm)--(57mm,23mm)--(65mm,23mm)--(65mm,25mm)--(68mm,25mm);
\draw(59mm,18mm)--(59mm,25mm)--(61mm,25mm)--(61mm,29mm)--(68mm,29mm);
\draw(61mm,18mm)--(61mm,21mm)--(63mm,21mm)--(63mm,27mm)--(68mm,27mm);
\draw(63mm,18mm)--(63mm,19mm)--(68mm,19mm);
\draw(65mm,18mm)--(65mm,21mm)--(68mm,21mm);
\draw(67mm,18mm)--(67mm,23mm)--(68mm,23mm);

\draw (80mm,18mm) \squarepath;
\draw (80mm,22mm) \rectanglepath;
\draw (80mm,24mm) \rectanglepath;
\draw (80mm,26mm) \rectanglepath;
\draw (80mm,28mm) \rectanglepath;
\draw (84mm,26mm) \rectanglepath;

\draw(81mm,18mm)--(81mm,21mm)--(87mm,21mm)--(87mm,23mm)--(89mm,23mm)--(89mm,25mm)--(92mm,25mm);
\draw(83mm,18mm)--(83mm,29mm)--(92mm,29mm);
\draw(85mm,18mm)--(85mm,25mm)--(87mm,25mm)--(87mm,27mm)--(92mm,27mm);
\draw(87mm,18mm)--(87mm,19mm)--(92mm,19mm);
\draw(89mm,18mm)--(89mm,21mm)--(92mm,21mm);
\draw(91mm,18mm)--(91mm,23mm)--(92mm,23mm);

\draw (8mm,0mm) \squarepath;
\draw[thick] (8mm,8mm) \rectanglepath;
\draw[thick] (8mm,10mm) \rectanglepath;
\draw[thick] (8mm,10mm) \rectanglepath;
\draw[thick] (10mm,8mm) \rectanglepath;
\draw[thick] (10mm,10mm) \rectanglepath;
\draw[thick] (12mm,10mm) \rectanglepath;

\draw(9mm,0mm)--(9mm,7mm)--(20mm,7mm);
\draw(11mm,0mm)--(11mm,5mm)--(13mm,5mm)--(13mm,9mm)--(15mm,9mm)--(15mm,11mm)--(20mm,11mm);
\draw(13mm,0mm)--(13mm,3mm)--(15mm,3mm)--(15mm,5mm)--(17mm,5mm)--(17mm,9mm)--(20mm,9mm);
\draw(15mm,0mm)--(15mm,1mm)--(20mm,1mm);
\draw(17mm,0mm)--(17mm,3mm)--(20mm,3mm);
\draw(19mm,0mm)--(19mm,5mm)--(20mm,5mm);

\draw (32mm,0mm) \squarepath;
\draw[thick] (32mm,6mm) \rectanglepath;
\draw[thick] (32mm,8mm) \rectanglepath;
\draw[thick] (32mm,10mm) \rectanglepath;
\draw[thick] (34mm,10mm) \rectanglepath;
\draw[thick] (36mm,10mm) \rectanglepath;

\draw(33mm,0mm)--(33mm,5mm)--(37mm,5mm)--(37mm,7mm)--(44mm,7mm);
\draw(35mm,0mm)--(35mm,9mm)--(39mm,9mm)--(39mm,11mm)--(44mm,11mm);
\draw(37mm,0mm)--(37mm,3mm)--(39mm,3mm)--(39mm,5mm)--(41mm,5mm)--(41mm,9mm)--(44mm,9mm);
\draw(39mm,0mm)--(39mm,1mm)--(44mm,1mm);
\draw(41mm,0mm)--(41mm,3mm)--(44mm,3mm);
\draw(43mm,0mm)--(43mm,5mm)--(44mm,5mm);

\draw (80mm,0mm) \squarepath;
\draw [thick](80mm,4mm) \rectanglepath;
\draw [thick](80mm,6mm) \rectanglepath;
\draw [thick](80mm,8mm) \rectanglepath;
\draw [thick](80mm,10mm) \rectanglepath;
\draw [thick](82mm,10mm) \rectanglepath;

\draw(83mm,0mm)--(83mm,9mm)--(85mm,9mm)--(85mm,11mm)--(92mm,11mm);
\draw(81mm,0mm)--(81mm,3mm)--(87mm,3mm)--(87mm,5mm)--(89mm,5mm)--(89mm,7mm)--(92mm,7mm);
\draw(85mm,0mm)--(85mm,7mm)--(87mm,7mm)--(87mm,9mm)--(92mm,9mm);
\draw(87mm,0mm)--(87mm,1mm)--(92mm,1mm);
\draw(89mm,0mm)--(89mm,3mm)--(92mm,3mm);
\draw(91mm,0mm)--(91mm,5mm)--(92mm,5mm);

\draw(50mm,71mm)--(62mm,67mm);
\draw(50mm,71mm)--(38mm,67mm);
\draw(38mm,53mm)--(26mm,49mm);
\draw(38mm,53mm)--(50mm,49mm);
\draw(62mm,53mm)--(74mm,49mm);

\draw(26mm,35mm)--(38mm,31mm);
\draw(26mm,35mm)--(14mm,31mm);
\draw(50mm,35mm)--(62mm,31mm);
\draw(74mm,35mm)--(86mm,31mm);

\draw(14mm,17mm)--(14mm,13mm);
\draw(38mm,17mm)--(38mm,13mm);
\draw(86mm,17mm)--(86mm,13mm);

\end{tikzpicture}
\caption{The EG-tree of $w=231654$.}\label{EGtree}
\end{center}
\end{figure}
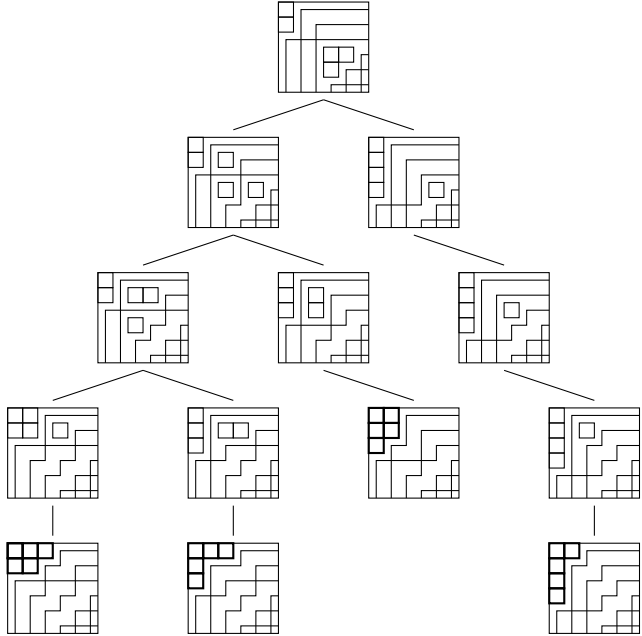

We  apply the above procedure to each path from the root to
a leaf in the modified LS-tree of $w$.
The resulting tree
is called the Edelman-Greene tree (EG-tree)
of $w$.  Figure \ref{EGtree} is the EG-tree
of  $w=231654$.

By the  construction of an EG-tree of $w$, we summarize
its relation to a modified LS-tree of $w$ in the following theorem.

\begin{theo}\label{lseg}
For  permutation $w$, the modified LS-tree of $w$ and the
EG-tree of $w$ are isomorphic. For any given  node labeled with
a permutation $u$ in the modified LS-tree of $w$, let $P$ be the corresponding
bumpless  pipedream of $w$ in the EG-tree of $w$. Then
 the empty boxes  of the Rothe pipedream of $u$
 are the same as the empty boxes of $P$.
\end{theo}

We have obtained a map from the leaves of an EG-tree
of $w$ to the EG-pipedreams of $w$. In the following theorem,
we shall give the reverse procedure.

\begin{theo}\label{EG-2}
Let $w\in S_n$ be a permutation and  $P$ be an EG-pipedream of $w$. Then
there is a leaf in the   EG-tree of $w$ whose label
is $P$.
\end{theo}

\pf
Assume that $P$ has $m$ NW elbows, say,
\[(i_m,j_m)<(i_{m-1},j_{m-1})<\cdots<(i_1, j_1),\]
 which are listed in
the total order as defined in  \eqref{totall}.
We construct a sequence
\begin{align}\label{invp}
P=P_m\rightarrow P_{m-1}\rightarrow\cdots\rightarrow P_0
\end{align}
of bumpless pipedreams of $w$ such that $P_m=P$ and $P_0=D(w)$.

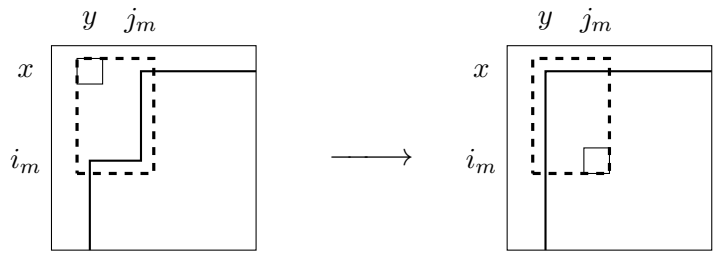
\begin{figure}[h]
\begin{center}
\begin{tabular}{cc}

\begin{tikzpicture}[scale = .85]

\def\rectanglepath{-- +(4mm,0mm) -- +(4mm,4mm) -- +(0mm,4mm) -- cycle}

\def\squarepath{-- +(32mm,0mm) -- +(32mm,32mm) -- +(0mm,32mm) -- cycle}

\draw[dashed,line width=1.2pt]
(4mm,12mm)--(16mm,12mm)--(16mm,30mm)--(4mm,30mm)--cycle;

\draw (0mm,0mm) \squarepath;
\draw (4mm,26mm)\rectanglepath;
\draw[thick](6mm,0mm)--(6mm,14mm)--(14mm,14mm)--(14mm,28mm)--(32mm,28mm);

\node at (-4mm,28mm) {\small{$x$}};
\node at (-4mm,14mm) {\small{$i_m$}};

\node at (6mm,36mm) {\small{$y$}};
\node at (14mm,36mm) {\small{$j_m$}};

\node at (50mm,15mm)
{$\xlongrightarrow{\qquad}$};
\end{tikzpicture}
&
\begin{tikzpicture}[scale = .85]

\def\rectanglepath{-- +(4mm,0mm) -- +(4mm,4mm) -- +(0mm,4mm) -- cycle}

\def\squarepath{-- +(32mm,0mm) -- +(32mm,32mm) -- +(0mm,32mm) -- cycle}

\draw (0mm,0mm) \squarepath;

\draw[dashed,line width=1.2pt]
(4mm,12mm)--(16mm,12mm)--(16mm,30mm)--(4mm,30mm)--cycle;

\draw (12mm,12mm)\rectanglepath;
\draw[thick](6mm,0mm)--(6mm,28mm)--(32mm,28mm);

\node at (-4mm,28mm) {\small{$x$}};
\node at (-4mm,14mm) {\small{$i_m$}};

\node at (6mm,36mm) {\small{$y$}};
\node at (14mm,36mm) {\small{$j_m$}};
\end{tikzpicture}
\end{tabular}
\caption{Reverse droop operation.}\label{fig-gendroop}
\end{center}
\end{figure}

First, we construct the  pipedream  $P_{m-1}$ from $P_m=P$.
The construction  is the same as that
 in the proof of \cite[Proposition 5.3]{LLS}
and is sketched below.
 Let $L$ be the pipe in $P_m$ passing through the NW elbow $(i_m,j_m)$.
Then   $L$ passes through an SE elbow $(i_m,y)$ (respectively,
$(x,j_m)$) in the same row (respectively, column).
Let  $R$ be the rectangle with corners $(i_m,y)$ and $(x,j_m)$.
It is easy to check that the northwest corner of $R$ is an empty box, and that
there are no any other elbows in $R$. Let $P_{m-1}$ be the pipedream
obtained from $P_m$ by a ``reverse droop'', that is,
change  the pipe $L$ to
travel  along the westmost column and
northmost row of $R$, see Figure \ref{fig-gendroop} for
an illustration.

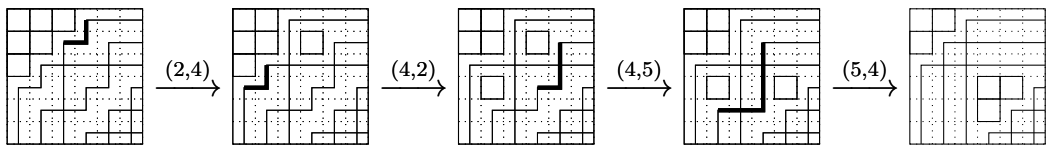
\begin{figure}[h]
\begin{center}
\begin{tikzpicture}

\def\rectanglepath{-- +(3mm,0mm) -- +(3mm,3mm) -- +(0mm,3mm) -- cycle}

\def\squarepath{-- +(18mm,0mm) -- +(18mm,18mm) -- +(0mm,18mm) -- cycle}

\draw (0mm,0mm) \squarepath;\draw [step=3mm,dotted] (0mm,0mm) grid (18mm,18mm);
\draw (0mm,9mm) \rectanglepath;\draw (0mm,12mm) \rectanglepath;
\draw (0mm,15mm) \rectanglepath;\draw (3mm,12mm) \rectanglepath;
\draw (3mm,15mm) \rectanglepath;\draw (6mm,15mm) \rectanglepath;

\draw(1.5mm,0mm)--(1.5mm,7.5mm)--(4.5mm,7.5mm)--(4.5mm,10.5mm)--(18mm,10.5mm);
\draw(4.5mm,0mm)--(4.5mm,4.5mm)--(10.5mm,4.5mm)--(10.5mm,7.5mm)--(13.5mm,7.5mm)--(13.5mm,13.5mm)--(18mm,13.5mm);
\draw(7.5mm,0mm)--(7.5mm,13.5mm);
\draw[-][ultra thick](7.5mm,13.5mm)--(10.5mm,13.5mm)--(10.5mm,16.5mm);
\draw(10.5mm,16.5mm)--(18mm,16.5mm);

\draw(10.5mm,0mm)--(10.5mm,1.5mm)--(18mm,1.5mm);
\draw(13.5mm,0mm)--(13.5mm,4.5mm)--(18mm,4.5mm);
\draw(16.5mm,0mm)--(16.5mm,7.5mm)--(18mm,7.5mm);


\node at (24mm,9mm) {$\xlongrightarrow[]{(2,4)}$};
\node at (54mm,9mm) {$\xlongrightarrow[]{(4,2)}$};
\node at (84mm,9mm) {$\xlongrightarrow[]{(4,5)}$};
\node at (114mm,9mm) {$\xlongrightarrow[]{(5,4)}$};

\draw (30mm,0mm) \squarepath;
\draw [step=3mm,dotted] (30mm,0mm) grid (48mm,18mm);
\draw (30mm,9mm) \rectanglepath;\draw (30mm,12mm) \rectanglepath;
\draw (30mm,15mm) \rectanglepath;
\draw (33mm,12mm) \rectanglepath;
\draw (33mm,15mm) \rectanglepath;
\draw (39mm,12mm) \rectanglepath;

\draw(31.5mm,0mm)--(31.5mm,7.5mm)--(34.5mm,7.5mm)--(34.5mm,10.5mm)--(48mm,10.5mm);
\draw[-][ultra thick](31.5mm,7.5mm)--(34.5mm,7.5mm)--(34.5mm,10.5mm);

\draw(34.5mm,0mm)--(34.5mm,4.5mm)--(40.5mm,4.5mm)--(40.5mm,7.5mm)--(43.5mm,7.5mm)--(43.5mm,13.5mm)--(48mm,13.5mm);
\draw(37.5mm,0mm)--(37.5mm,16.5mm)--(48mm,16.5mm);
\draw(40.5mm,0mm)--(40.5mm,1.5mm)--(48mm,1.5mm);
\draw(43.5mm,0mm)--(43.5mm,4.5mm)--(48mm,4.5mm);
\draw(46.5mm,0mm)--(46.5mm,7.5mm)--(48mm,7.5mm);

\draw (60mm,0mm) \squarepath;\draw [step=3mm,dotted] (60mm,0mm) grid (78mm,18mm);
\draw (63mm,6mm) \rectanglepath;\draw (60mm,12mm) \rectanglepath;
\draw (60mm,15mm) \rectanglepath;\draw (63mm,12mm) \rectanglepath;
\draw (63mm,15mm) \rectanglepath;\draw (69mm,12mm) \rectanglepath;

\draw(61.5mm,0mm)--(61.5mm,10.5mm)--(78mm,10.5mm);
\draw(64.5mm,0mm)--(64.5mm,4.5mm)--(70.5mm,4.5mm)--(70.5mm,7.5mm)--(73.5mm,7.5mm)--(73.5mm,13.5mm)--(78mm,13.5mm);
\draw[-][ultra thick](70.5mm,7.5mm)--(73.5mm,7.5mm)--(73.5mm,13.5mm);

\draw(67.5mm,0mm)--(67.5mm,16.5mm)--(78mm,16.5mm);
\draw(70.5mm,0mm)--(70.5mm,1.5mm)--(78mm,1.5mm);
\draw(73.5mm,0mm)--(73.5mm,4.5mm)--(78mm,4.5mm);
\draw(76.5mm,0mm)--(76.5mm,7.5mm)--(78mm,7.5mm);

\draw (90mm,0mm) \squarepath;\draw [step=3mm,dotted] (90mm,0mm) grid (108mm,18mm);
\draw (93mm,6mm) \rectanglepath;\draw (90mm,12mm) \rectanglepath;
\draw (90mm,15mm) \rectanglepath;\draw (93mm,12mm) \rectanglepath;
\draw (93mm,15mm) \rectanglepath;\draw (102mm,6mm) \rectanglepath;

\draw(91.5mm,0mm)--(91.5mm,10.5mm)--(108mm,10.5mm);
\draw(94.5mm,0mm)--(94.5mm,4.5mm)--(100.5mm,4.5mm)--(100.5mm,13.5mm)--(103.5mm,13.5mm)--(108mm,13.5mm);
\draw[-][ultra thick](94.5mm,4.5mm)--(100.5mm,4.5mm)--(100.5mm,13.5mm);

\draw(97.5mm,0mm)--(97.5mm,16.5mm)--(108mm,16.5mm);
\draw(100.5mm,0mm)--(100.5mm,1.5mm)--(108mm,1.5mm);
\draw(103.5mm,0mm)--(103.5mm,4.5mm)--(108mm,4.5mm);
\draw(106.5mm,0mm)--(106.5mm,7.5mm)--(108mm,7.5mm);

\draw (120mm,0mm) \squarepath;\draw[step=3mm,dotted] (120mm,0mm) grid (138mm,18mm);
\draw (129mm,3mm) \rectanglepath;\draw (120mm,12mm) \rectanglepath;
\draw (129mm,6mm) \rectanglepath;
\draw (123mm,15mm) \rectanglepath;\draw (132mm,6mm) \rectanglepath;

\draw(121.5mm,0mm)--(121.5mm,10.5mm)--(138mm,10.5mm);
\draw(124.5mm,0mm)--(124.5mm,13.5mm)--(138mm,13.5mm);
\draw(127.5mm,0mm)--(127.5mm,16.5mm)--(138mm,16.5mm);
\draw(130.5mm,0mm)--(130.5mm,1.5mm)--(138mm,1.5mm);
\draw(133.5mm,0mm)--(133.5mm,4.5mm)--(138mm,4.5mm);
\draw(136.5mm,0mm)--(136.5mm,7.5mm)--(138mm,7.5mm);

\draw (0mm,0mm) \squarepath;\draw [step=3mm,dotted] (0mm,0mm) grid (18mm,18mm);
\draw (0mm,9mm) \rectanglepath;\draw (0mm,12mm) \rectanglepath;
\draw (0mm,15mm) \rectanglepath;\draw (3mm,12mm) \rectanglepath;
\draw (3mm,15mm) \rectanglepath;\draw (6mm,15mm) \rectanglepath;

\draw(1.5mm,0mm)--(1.5mm,7.5mm)--(4.5mm,7.5mm)--(4.5mm,10.5mm)--(18mm,10.5mm);
\draw(4.5mm,0mm)--(4.5mm,4.5mm)--(10.5mm,4.5mm)--(10.5mm,7.5mm)--(13.5mm,7.5mm)--(13.5mm,13.5mm)--(18mm,13.5mm);
\draw(7.5mm,0mm)--(7.5mm,13.5mm);
\draw[-][ultra thick](7.5mm,13.5mm)--(10.5mm,13.5mm)--(10.5mm,16.5mm);
\draw(10.5mm,16.5mm)--(18mm,16.5mm);

\draw(10.5mm,0mm)--(10.5mm,1.5mm)--(18mm,1.5mm);
\draw(13.5mm,0mm)--(13.5mm,4.5mm)--(18mm,4.5mm);
\draw(16.5mm,0mm)--(16.5mm,7.5mm)--(18mm,7.5mm);


\node at (24mm,9mm) {$\xlongrightarrow[]{(2,4)}$};
\node at (54mm,9mm) {$\xlongrightarrow[]{(4,2)}$};
\node at (84mm,9mm) {$\xlongrightarrow[]{(4,5)}$};
\node at (114mm,9mm) {$\xlongrightarrow[]{(5,4)}$};

\draw (30mm,0mm) \squarepath;
\draw [step=3mm,dotted] (30mm,0mm) grid (48mm,18mm);
\draw (30mm,9mm) \rectanglepath;\draw (30mm,12mm) \rectanglepath;
\draw (30mm,15mm) \rectanglepath;
\draw (33mm,12mm) \rectanglepath;
\draw (33mm,15mm) \rectanglepath;
\draw (39mm,12mm) \rectanglepath;

\draw(31.5mm,0mm)--(31.5mm,7.5mm)--(34.5mm,7.5mm)--(34.5mm,10.5mm)--(48mm,10.5mm);
\draw[-][ultra thick](31.5mm,7.5mm)--(34.5mm,7.5mm)--(34.5mm,10.5mm);

\draw(34.5mm,0mm)--(34.5mm,4.5mm)--(40.5mm,4.5mm)--(40.5mm,7.5mm)--(43.5mm,7.5mm)--(43.5mm,13.5mm)--(48mm,13.5mm);
\draw(37.5mm,0mm)--(37.5mm,16.5mm)--(48mm,16.5mm);
\draw(40.5mm,0mm)--(40.5mm,1.5mm)--(48mm,1.5mm);
\draw(43.5mm,0mm)--(43.5mm,4.5mm)--(48mm,4.5mm);
\draw(46.5mm,0mm)--(46.5mm,7.5mm)--(48mm,7.5mm);

\draw (60mm,0mm) \squarepath;\draw [step=3mm,dotted] (60mm,0mm) grid (78mm,18mm);
\draw (63mm,6mm) \rectanglepath;\draw (60mm,12mm) \rectanglepath;
\draw (60mm,15mm) \rectanglepath;\draw (63mm,12mm) \rectanglepath;
\draw (63mm,15mm) \rectanglepath;\draw (69mm,12mm) \rectanglepath;

\draw(61.5mm,0mm)--(61.5mm,10.5mm)--(78mm,10.5mm);
\draw(64.5mm,0mm)--(64.5mm,4.5mm)--(70.5mm,4.5mm)--(70.5mm,7.5mm)--(73.5mm,7.5mm)--(73.5mm,13.5mm)--(78mm,13.5mm);
\draw[-][ultra thick](70.5mm,7.5mm)--(73.5mm,7.5mm)--(73.5mm,13.5mm);

\draw(67.5mm,0mm)--(67.5mm,16.5mm)--(78mm,16.5mm);
\draw(70.5mm,0mm)--(70.5mm,1.5mm)--(78mm,1.5mm);
\draw(73.5mm,0mm)--(73.5mm,4.5mm)--(78mm,4.5mm);
\draw(76.5mm,0mm)--(76.5mm,7.5mm)--(78mm,7.5mm);

\draw (90mm,0mm) \squarepath;\draw [step=3mm,dotted] (90mm,0mm) grid (108mm,18mm);
\draw (93mm,6mm) \rectanglepath;\draw (90mm,12mm) \rectanglepath;
\draw (90mm,15mm) \rectanglepath;\draw (93mm,12mm) \rectanglepath;
\draw (93mm,15mm) \rectanglepath;\draw (102mm,6mm) \rectanglepath;

\draw(91.5mm,0mm)--(91.5mm,10.5mm)--(108mm,10.5mm);
\draw(94.5mm,0mm)--(94.5mm,4.5mm)--(100.5mm,4.5mm)--(100.5mm,13.5mm)--(103.5mm,13.5mm)--(108mm,13.5mm);
\draw[-][ultra thick](94.5mm,4.5mm)--(100.5mm,4.5mm)--(100.5mm,13.5mm);

\draw(97.5mm,0mm)--(97.5mm,16.5mm)--(108mm,16.5mm);
\draw(100.5mm,0mm)--(100.5mm,1.5mm)--(108mm,1.5mm);
\draw(103.5mm,0mm)--(103.5mm,4.5mm)--(108mm,4.5mm);
\draw(106.5mm,0mm)--(106.5mm,7.5mm)--(108mm,7.5mm);

\end{tikzpicture}
\caption{Example of the construction of the sequence
 in \eqref{invp}.}\label{p2d}
\end{center}
\end{figure}

Using the same procedure as above, we can construct $P_{k-1}$ from $P_{k}$
$(1\leq k\leq m)$ by applying  a  ``reverse droop'' corresponding
 to
the NW elbow $(i_k,j_k)$.
Figure \ref{p2d} gives  an example to illustrate  the generation of
the chain from an EG-pipedream.

We show that the sequence in \eqref{invp} is a path
from a leaf to the root in the EG-tree of $w$.
Since $P_{k-1}$ has one fewer NW elbows than    $P_{k}$, the bumpless
pipedream $P_0$ has no NW elbows and hence is the Rothe pipedream   of $w$.
Set $u^{(0)}=w$.
By the construction of the sequence \eqref{invp} together with the fact that
the empty boxes of $P_m$ form a partition
at the northwest corner of the $n\times n$ grid, it is not hard to check that
$(i_1, j_1)$ is the largest box in the set $\mathrm{pivot}(w)$
under the total order  defined
in \eqref{totall}. In view of  Proposition  \ref{droop-app},
  the empty boxes of $P_{1}$ are also the empty boxes of the  Rothe pipedream of
  $u^{(1)}$  for some $u^{(1)}\in\Phi(wt_{p,q}, p)$.
Along  the same line, we can deduce that for $2\leq k\leq m$,
the empty boxes of $P_{k}$ form the empty boxes of the   Rothe pipedream of
some  $u^{(k)}$. In particular,
for $1\leq k\leq m$, $u^{(k)}$
is a child of $u^{(k-1)}$ in the modified LS-tree of $w$.
Hence, the sequence
\[w=u^{(0)}\rightarrow u^{(1)}\rightarrow \cdots \rightarrow u^{(m)}\]
forms a path in the modified LS-tree of $w$ from the root  to a leaf $u^{(m)}$. This shows that \eqref{invp} is a path
from a leaf to the root in the EG-tree of $w$.
Now we see that $P$ is the label of a leaf in the EG-tree of $w$.
This completes the proof.
\qed

It is easy to verify that the construction of \eqref{invp}
is the reverse process of the construction of \eqref{EG-22}.
Hence we arrive at the following conclusion.

\begin{coro}\label{COO}
For a permutation $w$, the labels of leaves
in the EG-tree of  $w$  are in bijection with
the EG-pipedreams of $w$.
\end{coro}

By Theorem \ref{EG-3}  and Corollary \ref{COO},
we obtain an alternative  proof of Theorem \ref{dsb}.

\section{The bijection}\label{PP-8}

In this section, we establish
the promised bijection
between reduced word  tableaux  and  EG-pipedreams.
For a permutation $w$, let $\mathrm{RT}(w)$ denote the set of
reduced word tableaux for $w$, namely, the set of increasing tableaux whose column reading words are reduced words of $w$. Let $\mathrm{EG}(w)$ denote the set of EG-pipedreams of $w$.

\begin{theo}\label{EG}
There is a shape preserving bijection  between  $\mathrm{RT}(w)$ and $\mathrm{EG}(w)$.
\end{theo}

By Theorem \ref{lseg} and Corollary \ref{COO},
we need only to establish a bijection between
the set $\mathrm{RT}(w)$ and the set of leaves in the
modified LS-tree of $w$. However,
we shall directly  construct
a  shape preserving bijection between
$\mathrm{RT}(w)$ and $\mathrm{EG}(w)$.
Of course, such a construction implies a bijection between
the set $\mathrm{RT}(w)$ and the set of leaves in the
modified LS-tree of $w$.

We first define two maps
\[\Gamma:\ \  \mathrm{RT}(w)\longrightarrow \mathrm{EG}(w).\]
and
\[\widetilde{\Gamma}:\ \  \mathrm{EG}(w)\longrightarrow \mathrm{RT}(w).\]
Then we show that they are the inverses of each other.

We need to employ the Little map for the transition in Theorem
\ref{CC-1}. In this case, the Little map   can be
written as $\theta_{p, w_q}$, since $q$ is the unique element of $S(wt_{p,q},p)$ defined in \eqref{snks}. Note that $(p, w_q)$ is the maximum  empty box in the Rothe pipedream
of $w$ which has a pivot. Moreover, for a reduced word $a$ of a permutation
in $\Phi(wt_{p,q}, p)$, we have
\[\theta_{p, w_q}^{-1}(a)=(\theta_{n+1-p, n+1-w_q}(a^c))^c,\]
which is a reduced word of $w$.

Let $T\in \mathrm{RT}(w)$ be a reduced word tableau for $w$ with shape $\lambda$.
We first construct a path of permutations
\begin{equation}\label{BUT}
w=w^{(0)}\rightarrow w^{(1)}\rightarrow \cdots \rightarrow w^{(m)}
\end{equation}
in the modified LS-tree of $w$ from the root
$w^{(0)}=w$ to a leaf $w^{(m)}$. Then there is a path in the EG-tree of $w$ which corresponds to the path in \eqref{BUT}. Let $P$ be the EG-pipedream corresponding to $w^{(m)}$ in the EG-tree of $w$. Define $\Gamma(T)=P$.

To construct \eqref{BUT}, we need a path of reduced words
\begin{equation}\label{FI}
\tau=\tau^{(0)}\rightarrow \tau^{(1)}\rightarrow \cdots \rightarrow \tau^{(m)}
\end{equation} such that $\tau^{(i)}$ is a reduced word of $w^{(i)}$ for $0\le i\le m$.
Let $\tau^{(0)}=\mathrm{column}(T)$. By definition,
$\mathrm{column}(T)$ is a reduced word of $w$. Note that $T$ is the insertion
tableau of $\mathrm{column}(T)^{\mathrm{rev}}$ under the Edelman-Greene algorithm.
For $0\le i\le m-1$, let $\tau^{(i+1)}$ be obtained from $\tau^{(i)}$ by applying the Little map $\theta_{p_i, w^{(i)}_{q_i}}$, where $(p_i,w^{(i)}_{q_i})$ is
the maximum empty box in the Rothe pipedream  of $w^{(i)}$ which has a pivot.

Conversely, let $P$ be an EG-pipedream of $w$ with $m$ NW elbows, say
\[ (i_{m},j_{m})<\cdots<(i_1, j_1)\]
in the order \eqref{totall}. Assume that $w'$ is the permutation in the modified LS-tree of $w$ which corresponds to $P$. Then $w'$ is a dominant permutation with a unique reduced word tableau, say $T'$, namely, the frozen tableau of $w'$ as mentioned in Section \ref{lstree}. Take the column reading word  $\mathrm{column}(T')$
of $w'$.
 Let  $w(P)$ be the reduced word of $w$ defined by
\begin{align}\label{gamma}
w(P)=\theta_{i_1,j_1}^{-1}\circ\cdots\circ\theta_{i_{m}, j_{m}}^{-1}(\mathrm{column}(T')).
\end{align}
Let $\widetilde{T}$  be the insertion tableau of $w(P)^{\mathrm{rev}}$ by the Edelman-Greene  algorithm. Define $\widetilde{\Gamma}(P)=\widetilde{T}$.
By Theorem \ref{90}, the row reading word of $\widetilde{T}$
is a reduced word of $w^{-1}$.
By the proof of Theorem \ref{91},
the column reading word of  $\widetilde{T}$ is a reduced word of $w$, and thus $\widetilde{T}$
is a reduced word tableau for $w$.

\noindent{\it Proof of Theorem 6.1}.
We show that   $\widetilde{\Gamma}$ is
the inverse of $\Gamma$. Let $P=\Gamma(T)$ where $T\in\mathrm{RT}(w)$.
We need to show that $\widetilde{T}=\widetilde{\Gamma}(P)=T$.
For $0\leq i\leq m$, let $T_i$ be the insertion
tableau of $(\tau^{(i)})^{\mathrm{rev}}$ for
the reduced words in \eqref{FI}. Note that $T_0=T$.
 By Theorem \ref{lam-HY}, the tableaux $T_i$ have the same shape.
 Moreover, by Theorem \ref{HY2}, the column reading word of $T_i$
is $\tau^{(i)}$.
By the construction of the  EG-tree of $w$ and in view of
the construction of $\widetilde{T}$, it is easy to check that $\widetilde{T}=T$.
In the same manner, it is also easy to verify that $\Gamma$ is
the inverse of $\widetilde{\Gamma}$.
So the proof is complete. \qed

For example, let $w=231654$ and let $T$ be the following reduced word
tableau of $w$:
\[T=\begin{ytableau}
\scriptstyle1 &  \scriptstyle4  & \scriptstyle5\\
\scriptstyle2 \\
 \scriptstyle5
\end{ytableau}.\]
We see that $\mathrm{column}(T)=(5,4,1,2,5)$.
Let us construct the path of permutations in \eqref{BUT}.
Since the maximum box in $\mathrm{pivot}(w)$  is $(5,4)$, we should apply $\theta_{5,4}$ to $\mathrm{column}(T)=(5,4,1,2,5)$. That is,
\[\theta_{5,4}(5,4,1,2,5)=(5,3,1,2,4),\]
which is a reduced word  for $w^{(1)}=241635$.
The maximum box in $\mathrm{pivot}(w^{(1)})$  is $(4,5)$. Then we should apply $\theta_{4,5}$ to $53124$ to yield
\[\theta_{4,5}(5,3,1,2,4)=(4,3,1,2,4),\]
which is a reduced word of $w^{(2)}=251436$.
The maximum box in  $\mathrm{pivot}(w^{(2)})$ is $(4,3)$. Thus we apply $\theta_{4,3}$ to $43124$ to obtain
\[\theta_{4,3}(4,3,1,2,4)=(4,3,1,2,3),\]
which is a reduced word of $w^{(3)}=253146$.
The maximum box in  $\mathrm{pivot}(w^{(3)})$   is $(2,4)$. Then we apply $\theta_{2,4}$ to $(4,3,1,2,3)$ to yield
\[\theta_{2,4}(4,3,1,2,3)=(3,2,1,2,3),\]
which is a reduced word of $w^{(4)}=423156$. Since $w^{(4)}$ is a dominant permutation, the Rothe diagram $\mathrm{Rothe}(w^{(4)})$ is a partition, and  we stop. Therefore, the   path of permutations
in the  modified  LS-tree of $w$ is \[231654\rightarrow241635\rightarrow251436\rightarrow253146\rightarrow423156.\]
By the EG-tree of $w$ displayed in Figure \ref{EGtree} , the EG-pipedream corresponding to the leaf $423156$ is
\begin{figure}[h]
\begin{center}
\begin{tikzpicture}

\def\rectanglepath{-- +(3mm,0mm) -- +(3mm,3mm) -- +(0mm,3mm) -- cycle}

\def\squarepath{-- +(18mm,0mm) -- +(18mm,18mm) -- +(0mm,18mm) -- cycle}

\draw (0mm,0mm)\squarepath;
\draw [step=3mm,dotted] (0mm,0mm) grid (18mm,18mm);
\draw (0mm,9mm) \rectanglepath;
\draw (0mm,12mm)\rectanglepath;
\draw (0mm,15mm)\rectanglepath;
\draw (3mm,15mm)\rectanglepath;
\draw (6mm,15mm)\rectanglepath;

\draw(1.5mm,0mm)--(1.5mm,7.5mm)--(7.5mm,7.5mm)--(7.5mm,10.5mm)--(18mm,10.5mm);
\draw(4.5mm,0mm)--(4.5mm,13.5mm)--(10.5mm,13.5mm)--(10.5mm,16.5mm)--(18mm,16.5mm);
\draw(7.5mm,0mm)--(7.5mm,4.5mm)--(10.5mm,4.5mm)--(10.5mm,7.5mm)--(13.5mm,7.5mm)--(13.5mm,13.5mm)--(18mm,13.5mm);
\draw(10.5mm,0mm)--(10.5mm,1.5mm)--(18mm,1.5mm);
\draw(13.5mm,0mm)--(13.5mm,4.5mm)--(18mm,4.5mm);
\draw(16.5mm,0mm)--(16.5mm,7.5mm)--(18mm,7.5mm);


\node  [left] at (0mm,9mm) {$\Gamma(T)=$};
\node   at (19mm,8mm) {.};

\end{tikzpicture}
\end{center}
\end{figure}

Conversely, let $P$ be the second EG-pipedream in the bottom row of Figure \ref{EGtree}. The set of NW elbows of $P$ is $\{(2,4),(4,3),(4,5),(5,4)\}_<$. The leaf in the modified LS-tree of $w$  corresponding to $P$ is $w'=423156$, which has a unique reduced word tableau
\[
T'=\begin{ytableau}
\scriptstyle1 &  \scriptstyle2  & \scriptstyle3\\
\scriptstyle2 \\
 \scriptstyle3
\end{ytableau}.
\]
Then $\mathrm{column}(T')=(3,2,1,2,3)$ and \[w(P)=\theta_{5,4}^{-1}\circ\theta_{4,5}^{-1}\circ\theta_{4,3}^{-1}\circ\theta_{2,4}^{-1} (3,2,1,2,3).\]
We have the following calculations:
\begin{align*}
  &\theta_{2,4}^{-1}(3,2,1,2,3)=(\theta_{5,3}^{-1}(3,4,5,4,3))^c=(2,3,5,4,3)^c=(4,3,1,2,3);\\[5pt]
  &\theta_{4,3}^{-1}(4,3,1,2,3)=(\theta_{3,4}^{-1}(2,3,5,4,3))^c=(2,3,5,4,2)^c=(4,3,1,2,4);\\[5pt]
  &\theta_{4,5}^{-1}(4,3,1,2,4)=(\theta_{3,2}^{-1}(2,3,5,4,2))^c=(1,3,5,4,2)^c=(5,3,1,2,4);\\[5pt]
  &\theta_{5,4}^{-1}(5,3,1,2,4)=(\theta_{2,3}^{-1}(1,3,5,4,2))^c=(1,2,5,4,1)^c=(5,4,1,2,5).
\end{align*}
Therefore, $w(P)=(5,4,1,2,5)$. Finally,   insert  $w(P)^{\mathrm{rev}}=(5,2,1,4,5)$ by the Edelman-Greene algorithm to obtain
\[T=\begin{ytableau}
\scriptstyle1 &  \scriptstyle4  & \scriptstyle5\\
\scriptstyle2 \\
 \scriptstyle5
\end{ytableau},
\]
which is a reduced word tableau for $w=231654$.

\vskip 3mm \noindent {\bf Acknowledgments.}
Part of this work was completed during Neil Fan was visiting the Department of Mathematics at the University of Illinois at Urbana-Champaign, he wishes to thank the department for its hospitality and thank Alexander Yong for helpful conversations.
This work was supported
by the 973 Project, the PCSIRT Project of
the Ministry of Education, the National Science Foundation of
China.


\begin{thebibliography}{99}

\bibitem{BB}
N. Bergeron and S. Billey, RC-graphs and Schubert polynomials, Experiment. Math. 2 (1993), 257--269.


\bibitem{BeSo}
N. Bergeron and F. Sottile, Schubert polynomials, the Bruhat order, and the geometry of flag manifolds, Duke Math. J. 95 (1998),  373--423.

\bibitem{BJS}
S. Billey, W. Jockusch and R.P. Stanley, Some combinatorial properties of Schubert
polynomials, J. Algebraic Combin. 2 (1993), 345--374.

\bibitem{BP}
S. Billey and B. Pawlowski, Permutation patterns, Stanley symmetric functions, and generalized Specht modules, J. Combin. Theory Ser. A  127 (2014),  85--120.


\bibitem{Buch}
A.S. Buch, A. Kresch, M. Shimozono, H. Tamvakis and A. Yong, Stable
Grothendieck polynomials and K-theoretic factor sequences, Math. Ann. 340
(2008), 359--382.


\bibitem{Chen}
W.Y.C. Chen, P.L. Guo and S.X.M. Pang,
Vacillating Hecke tableaux and linked partitions,
Math. Z. 281 (2015), 661--672.

\bibitem{Felsner}
S. Felsner, The Skeleton of a reduced word and a
correspondence of Edelman and Greene, Electron. J.
Combin. 8 (2001), \#R10.


\bibitem{Fomin-Greene}
 S. Fomin and C. Greene, Noncommutative Schur functions and their applications, Discrete
Math. 193 (1998), 179--200.

%

\bibitem{EG}
P. Edelman and C. Greene, Balanced tableaux, Adv. Math., 63   (1987), 42--99.


\bibitem{FK2}
S. Fomin and A. N. Kirillov, The Yang-Baxter equation, symmetric functions, and Schubert polynomials,
Proceedings of the 5th Conference on Formal Power Series and Algebraic Combinatorics (Florence 1993), Discrete Math. 153 (1996),  123--143.

\bibitem{FS}
S. Fomin and R.P. Stanley, Schubert polynomials and the NilCoxeter algebra, Adv.
Math. 103 (1994), 196--207.

\bibitem{Garsia}
A.M. Garsia, The Saga of Reduced Factorizations of Elements of the
Symmetric Group, Publications du
LaCIM, Universit\'e du Qu\'ebec \'a Montr\'eal, Canada, Vol. 29, 2002.

\bibitem{Ful}
W. Fulton, Flags, Schubert polynomials, degeneracy loci, and determinantal formulae, Duke Math. J. 65 (1992), 381--420.




\bibitem{KnMi}
A. Knutson and E. Miller, Gr\"obner geometry of Schubert polynomials, Ann. Math.
161 (2005), 1245--1318.

\bibitem{KM}
A. Knutson and E. Miller, Subword complexes in Coxeter groups, Adv. Math. 184
(2004), 161--176.


\bibitem{KMY}
A. Knutson, E. Miller and A. Yong, Gr\"obner geometry of vertex decompositions and of flagged tableaux, J. Reine Angew. Math. 63 (2009), 1--31.

\bibitem{KY}
A. Knutson and A. Yong, A formula for $K$-theory truncation Schubert calculus, Int. Math. Res. Not.  70 (2004), 3741--3756.

\bibitem{Lam}
 T. Lam, Stanley symmetric functions and Peterson algebras. In: Lam, T., Lapointe, L., Morse, J., Schilling, A., Shimozono, M., Zabrocki, M. (eds) $k$-Schur Functions and Affine Schubert Calculus, Fields Institute Monographs, Vol. 33, pp. 133--168. Springer, New York, 2014.


\bibitem{LLS}
T. Lam, S. Lee and M. Shimozono, Back stable Schubert calculus, arXiv:1806.11233v1.


\bibitem{LS1}
A. Lascoux and M. Sch\"{u}tzenberger,
Polyn\^{o}mes de Schubert, C.R. Acad. Sci. Paris 294 (1982), 447--450.



\bibitem{LS}
A. Lascoux and M. Sch\"{u}tzenberger, Schubert polynomials and the Littlewood-Richardson rule, Lett. Math. Phys. 10 (2-3) (1985), 111--124.



\bibitem{Linusson}
S. Linusson and S. Potka,
New properties of the Edelman-Greene bijection,
arXiv:1804.10034v1.

\bibitem{Little}
D. P. Little, Combinatorial aspects of the
Lascoux-Sch\"{u}tzenberger tree, Adv. Math. 174  (2003), 236--253.


%
\bibitem{Marberg}
 Z. Hamaker, E. Marberg and B. Pawlowski,
Schur $P$-positivity and involution Stanley symmetric functions,
Int. Math. Res. Not. (2017), rnx274.

\bibitem{HY}
Z. Hamaker and B. Young, Relating Edelman-Greene insertion to the Little map,
J. Algebr. Comb.,  40 (2014), 693--710.

\bibitem{Mac}
I.G. Macdonald, Notes on Schubert Polynomials, Montr\'{e}al: D\'{e}p. de math\'{e}matique et d'informatique, Universit\'{e} du Qu\'{e}bec \`{a} Montr\'{e}al, 1991.





\bibitem{Sta-5}
R.P. Stanley, Enumerative Combinatorics, Vol. 1, Second edition, Cambridge Studies in Advanced Mathematics,  Cambridge University Press, Cambridge, 2012.


\bibitem{Sta}
R. Stanley, On the number of reduced decompositions of elements of Coxeter groups, European J. Combin. 5 (1984), 359--372.

\bibitem{Sta-2}
 R.P. Stanley,
Reduced decompositions, http://www-math.mit.edu/~rstan/\\
transparencies/redec-ams.pdf

\bibitem{TY}
H. Thomas and A. Yong, A jeu de taquin theory for increasing tableaux, with applications to $K$-theoretic Schubert calculus, Algebra Number Theory  3  (2009), 121--148.

\bibitem{Wa}
M. Wachs, Flagged Schur functions, Schubert polynomials, and symmetrizing
operators, J. Combin. Theory  Ser. A. 40 (1985), 276--289.


\bibitem{WeYo}
A. Weigandt and A. Yong, The prism tableau model for Schubert polynomials,
 J. Combin. Theory Ser. A 154 (2018), 551--582.


\bibitem{Wi}
R. Winkel, Diagram rules for the generation of Schubert polynomials, J. Combin. Theory
Ser. A. 86 (1999), 14--48.

\end{thebibliography}
\end{document}